\documentclass{article}
\usepackage{amsthm}
\usepackage{fancyhdr}
\usepackage{amsmath}
\usepackage{indentfirst}
\usepackage{titlesec}
\usepackage{amsfonts}
\usepackage{amssymb}

\newcommand{\cA}{\mathcal A}
\newcommand{\cB}{\mathcal B}

\newcommand{\cF}{\mathcal F}

\newcommand{\cH}{\mathcal H}
\newcommand{\cL}{\mathcal L}

\newcommand{\cR}{\mathcal R}

\newcommand{\N}{\mathbb N}
\newcommand{\R}{\mathbb R}
\newcommand{\Z}{\mathbb Z}
\newcommand{\cov}{\textrm{\rm Cov}}

\newcommand{\var}{\textrm{\rm Var}}

\newcommand{\zhfb}{\hfil\break}
\newcommand{\zeps}{\varepsilon}

\newcommand{\bbN}{\mathbb N}
\newcommand{\bbP}{\mathbb P}
\newcommand{\bbR}{\mathbb R}
\newcommand{\bbZ}{\mathbb Z}

 \newtheorem{theorem}{Theorem}[section]

 \newtheorem{lemma}[theorem]{Lemma}

 \newtheorem{corollary}[theorem]{Corollary}

 \theoremstyle{definition}
 \newtheorem{definition}[theorem]{Definition}

\newtheorem{remark}[theorem]{Remark}

 \newtheorem{notations}[theorem]{Notations}

 \newtheorem{background}[theorem]{Background}

 \numberwithin{equation}{section}

 \numberwithin{theorem}{section}

\begin{document}

% Full title of the paper (Capitalized)
\title{\textbf{On the central limit question for 
strictly stationary, reversible Markov chains}}

% Author Orchid ID: enter ID or remove command
\newcommand{\orcidauthorA}{0000-0000-000-000X} % Add \orcidA{} behind the author's name
%\newcommand{\orcidauthorB}{0000-0000-000-000X} % Add \orcidB{} behind the author's name

% Authors, for the paper (add full first names)
\author{Richard C.\ Bradley \\
Department of Mathematics, 
Indiana University, Bloomington, IN 47405, U.S.A.\\ bradleyr@iu.edu}

% Authors, for metadata in PDF
%\AuthorNames{Richard Bradley \\%}%Firstname Lastname, Firstname Lastname and Firstname Lastname}

\maketitle

\begin{abstract}
     This paper will provide several classes of strictly stationary, countable-state, irreducible, aperiodic Markov chains that are reversible and have finite second moments, such that the central limit theorem fails to hold.
The main purpose is to examine the extent to which, for the development of central limit theory for strictly stationary Markov chains (and functions of them) under the strong mixing ($\alpha$-mixing) and absolute regularity ($\beta$-mixing) conditions, the property of reversibility (if it holds) can provide extra leverage.
It is known, partly as a by-product of research done by Roberts, Rosenthal, and 
Tweedie in two papers in 1997 and 2001, that for the case of exponential mixing rates,  
reversibility provides notable extra leverage of that kind.  
In contrast, a class of counterexamples in a paper of Doukhan, Massart, and Rio in 1994
showed (implicitly) that for the case of power-type mixing rates, reversibility
apparently provides almost no such extra leverage.
Further perspective on that latter fact will be provided by some counterexamples 
in this paper.
Other counterexamples here will (indirectly) provide some tentative, uncertain 
evidence for the possibility that for mixing rates that are ``between'' power-type and 
exponential (for example, sub-exponential), reversibility may in fact provide some 
small but nontrivial extra leverage.    
\end{abstract}

%\maketitle
% Keywords
\hfil\break 
\null \hskip 1.2 in 
strong mixing; absolute regularity; central limit question
\medskip


%\maketitle

%%%%%%%%%%%%%%%%%%%%%%%%%%%%%%%%%%%%%%%%%%

%\setcounter{secnumdepth}{4}
%%%%%%%%%%%%%%%%%%%%%%%%%%%%%%%%%%%%%%%%%%

\section{Introduction}
\label{sc1}

   In this paper, for simplicity, every Markov chain
(i) is indexed by the set $\Z$ of all integers,
(ii) is strictly stationary, and
(iii) has as its state space either the set $\R$ of all real numbers 
(accompanied by its Borel $\sigma$-field
$\cR$) or a finite or countably infinite subset of $\R$.
\smallskip

   Over twenty years ago, in the papers of 
Roberts and Rosenthal \cite {ref-journal-RR} and 
Roberts and Tweedie \cite {ref-journal-RT} together,
some special theory for Markov chains was developed which included
 (among other things) some striking results for strictly stationary  
 Markov chains that are reversible.
 Of course the term ``reversible'' simply means that
the distribution of the Markov chain as a whole is invariant
under a reversal of the ``direction of time''.   
(More on that in Section \ref {sc2}.)\ \
It is well known that some key aspects of the theory in those two papers
\cite {ref-journal-RR} and \cite {ref-journal-RT} can be reformulated
in the terminology of the dependence coefficients associated with
certain strong mixing conditions; see e.g.\ 
 \cite {ref-journal-KM}, \cite {ref-journal-Longla14},
\cite {ref-journal-LP}, or the exposition in  
\cite {ref-journal-Bradley2021}. 
For example, it has long been well known (see e.g.\ 
[\cite {ref-journal-Bradley2007}, Theorem 21.19]) that for a given strictly
stationary (not necessarily reversible) Markov chain, the ``uniform ergodicity''
condition (which plays a prominent role in \cite {ref-journal-RR} and
\cite {ref-journal-RT}) is equivalent to $\beta$-mixing (absolute regularity) with a
mixing rate that is at least exponentially fast
(``exponential $\beta$-mixing'' for short).   
\medskip

     It has been well known that as a consequence of the aforementioned 
special theory developed in the papers 
\cite {ref-journal-RR} and \cite {ref-journal-RT},
together with classic central limit theorems 
(such as in Rosenblatt \cite {ref-journal-Rosenblatt1971}
or Lifshits \cite {ref-journal-Lifshits}) for (functions of) strictly 
stationary Markov chains that satisfy the $\rho$-mixing condition,
the following holds: 
For a (function of a) strictly stationary, reversible, ``exponentially $\beta$-mixing'' 
Markov chain, if the second moments are finite and the variances of the partial 
sums diverge to $\infty$, then a central limit theorem holds.
It is also known, from counterexamples in 
\cite {ref-journal-Bradley1983b}, \cite {ref-journal-CL2025}, 
\cite {ref-journal-Haggstrom}, and
(after a very careful examination of certain hidden details)
\cite {ref-journal-Herrndorf}, that if the assumption of reversibility were omitted 
altogether in that context, then the central limit theorem may fail to hold.
Thus for (functions of) strictly stationary, ``exponentially $\beta$-mixing'' 
Markov chains, when the second moments are finite, 
the extra property of reversibility, if satisfied, provides extra leverage for the
central limit theorem to hold.
\medskip

   Recently Cuny and Lin \cite {ref-journal-CL2025} examined that general theme in 
much detail in a broader context, with a special focus on the weaker mixing assumption of 
$\alpha$-mixing (Rosenblatt strong mixing), with some emphasis on an exponential mixing rate
or absence thereof, 
and with the condition of reversibility replaced by certain related weaker conditions.
Among other things, they extended the central limit theorem alluded to above to 
one involving that broader context.          
Some recent related work by the author \cite {ref-journal-Bradley2025} 
sheds some further light on that general theme in a different way,
with a particular class of concrete counterexamples (to the central limit theorem) 
involving reversibility and $\beta$-mixing.  
More details on all of this will be provided in Theorem \ref{thm2.7} and Remarks \ref{rem2.8}
and \ref{rem2.9} in Section \ref{sc2}.
\medskip
 
     Building on the work in \cite {ref-journal-Bradley2025},  
this paper here is intended to provide some further insights into 
that general theme, with the construction of various other classes of counterexamples 
intended to illustrate in certain concrete ways ``what can go wrong'' with respect to the 
central limit theorem under the seemingly ``nicest'' (i.e.\ most favorable) assumptions.
Accordingly, we shall retain our focus here on reversibility and (in the 
counterexamples) on $\beta$-mixing (which is stronger than $\alpha$-mixing).
\medskip

     In this paper, the following question will be examined:  
Suppose that for a given strictly stationary Markov chain, 
(i) $\beta$-mixing (absolute regularity) holds with a mixing rate
that is slower than exponential, and 
(ii) for the Markov chain itself or a given function of it, the marginal distribution is 
assumed to satisfy a ``moment'' condition (or a more 
general ``tail'' condition) that is stronger than finite second moment.
To what extent does the extra property of reversibility, if it is satisfied,
provide extra leverage for the central limit theorem to hold?
\medskip

     In this paper, some counterexamples will be constructed --- each of them being
simply a strictly stationary, countable-state, reversible Markov chain --- leading to the 
following tentative assessments regarding that question:
\medskip

For the case where the random variables in the Markov chain (or a given function of it) 
are bounded or at least have finite absolute moments of some order higher than second moments, the answer is in some sense ``practically no extra leverage''.
That was shown to hold in some sense by a class of counterexamples in the paper 
of Doukhan et al.\ \cite {ref-journal-DMR1994}.  
Further nuances on that will be provided by counterexamples constructed here in this 
paper (specifically, as described in Theorem \ref{thm4.4} and Corollary \ref{cor5.7}).
\medskip 
 
For the case where the second moments are finite but absolute moments 
of order $2 + \delta$ are infinite for all $\delta > 0$, the answer 
remains uncertain, but from counterexamples constructed in this paper,
it seems plausible that the answer might be something like 
``some small but nontrivial extra leverage''.
That musing will be put into a more precise form in Remarks \ref{rem5.9} and \ref{rem5.10}
in Section \ref{sc5} in this paper, in connection with some counterexamples
described in that section (in particular, in Corollary \ref{cor5.8}).    
\medskip

     The rest of this paper is organized as follows:
Section \ref{sc2} will give a review of certain mixing conditions for Markov 
chains, along with further relevant technical details on the discussion above.
Sections \ref{sc3} and \ref{sc4} will describe counterexamples that 
respectively deal with two different aspects of the question above for the 
case of strictly stationary, reversible Markov chains that are bounded.
Section \ref{sc5} will describe some counterexamples 
that deal broadly with the question above for the case of
strictly stationary, reversible Markov chains that are not bounded. 
Section \ref{sc6} will provide some preliminary work for
subsequent proofs. 
Sections \ref{sc7}, \ref{sc8}, and \ref{sc9} will then give proofs of 
the properties of the counterexamples described in Sections 
\ref{sc3}, \ref{sc4}, and \ref{sc5} respectively.
Finally, Section \ref{sc10} will take a brief further look at the (bounded) counterexamples
constructed in Section \ref{sc8}, in connection with a particular very sharp assumption
in a related central limit theorem of Merlev\`ede and Peligrad \cite {ref-journal-MP2}.
\medskip

\section{Mixing conditions and Markov chains}
\label{sc2}

   In this section, we shall lay out some notations and review some 
background facts regarding mixing conditions and Markov chains.
Most of the exposition in this section will be taken practically 
verbatim from the introductions of recent papers by the 
author \cite {ref-journal-Bradley2024} \cite {ref-journal-Bradley2025}.

\begin{notations}
\label{nt2.1}
First, here are a few really basic notations:\hfil\break
$\N$ denotes the set of all positive integers, \hfil\break
$\Z$ denotes the set of all integers, \hfil\break
$\R$ denotes the set of all real numbers, and \hfil\break
$\cR$ denotes the Borel $\sigma$-field on $\R$.
\hfil\break
The usual notations such as $\bbR^\bbZ$ and $\cR^\bbZ$ will be used for 
Cartesian products of sets and for product $\sigma$-fields.
\smallskip

   For sequences $(a_1, a_2, a_3, \dots)$ and $(b_1, b_2, b_3, \dots)$ of 
positive numbers, the notations $a_n = O(b_n)$ and $a_n << b_n$
(as $n \to \infty$) have the same meaning
(the latter notation will be used primarily when the numbers 
$a_n$ and $b_n$ are complicated); 
the notation $a_n \asymp b_n$ (as $n \to \infty$) will mean that
$a_n=O(b_n)$ and $b_n = O(a_n)$;
and the notation
$a_n \sim b_n$ (as $n \to \infty$) will mean that
$\lim_{n \to \infty} (a_n/b_n) = 1$.
\smallskip
   
     In what follows, the setting will be a probability space
$(\Omega, \cF, P)$, rich enough to accommodate all random variables declared.
In this paper, all random variables will be real-valued (possibly discrete).
The $\sigma$-field (on $\Omega$) generated by a given family $(Y_i, i \in I)$ of  
random variables (with the index set $I$ being nonempty), will be denoted
$\sigma(Y_i, i \in I)$.
\smallskip   

    The phrase ``central limit theorem'' will henceforth be abbreviated simply as CLT.
\end{notations}

\begin{notations}
\label{nt2.2}
Suppose $\cA$ and $\cB$ are any 
two $\sigma$-fields $\subset \cF$.
Define the following two measures of dependence:
First,
\begin{equation}
\label{eq2.21}
\alpha(\cA,\cB) := \sup_{A\in \cA, B\in \cB} |P(A\cap B)-P(A) P(B)|.
\end{equation}
Second, 
\begin{equation}\label{eq2.22}
\beta(\cA,\cB):= 
\sup \frac{1}{2} \sum^I_{i=1}\sum^J_{j=1} |P(A_i\cap B_j) - P(A_i)P(B_j)|
\end{equation}
where the supremum is taken over all pairs of finite partitions 
$\{A_1,A_2,\dots,A_I\}$ and $\{B_1,B_2,\dots,B_J\}$ of 
$\Omega$ such that $A_i\in \cA$ for each $i$ and $B_j\in \cB$ for each $j$.
[The factor of $1/2$ in (\ref{eq2.22}) is not of special significance, 
but has become customary in order to make 
certain inequalities a little nicer.]\ \ 
It is elementary and well known that
\begin{equation}
\label{eq2.23}
0 \leq 2\alpha(\cA,\cB) \leq \beta(\cA,\cB) \leq 1.
\end{equation}
(See e.g.\  [\cite {ref-journal-Bradley2007}, Vol.\ 1, Proposition 3.11].)\ \ 
The quantities $\alpha(\cA,\cB)$ and $\beta(\cA,\cB)$ are 
equal to $0$ if the $\sigma$-fields $\cA$ and $\cB$ are independent, 
and are each positive otherwise.
\end{notations}

\begin{notations}
\label{nt2.3}
Suppose $X:= (X_k$, $k\in \Z)$ is a (not necessarily Markovian) 
strictly stationary sequence of (real-valued) random variables.
Refer to Notations \ref{nt2.1}.
For each integer $j$, define the notations $\cF^j_{-\infty} :=\sigma(X_k, k\le j)$ and 
$\cF^\infty_j:= \sigma(X_k$, $k\ge j)$.
\medskip

     For each positive integer $n$, define the following two ``dependence coefficients'':
\begin{align}
\label{eq2.31}
\alpha(n) &= \alpha_X(n) := \alpha(\cF^0_{-\infty}, \cF^\infty_n); \quad {\rm and}\\
\label{eq2.32}
\beta(n) &= \beta_X(n) := \beta(\cF^0_{-\infty}, \cF^\infty_n).
\end{align}
For each positive integer $n$, one has by strict stationarity that 
$\alpha(n) = \alpha(\cF^j_{-\infty},\cF^\infty_{j+n})$ for every integer 
$j$, and the analogous comment holds for $\beta(n)$ as well.
\smallskip

     One (trivially) has that $\alpha(1) \ge \alpha(2)\ge \alpha(3)\ge \dots\,$; and the analogous comment holds for the numbers $\beta(n)$.
The strictly stationary sequence $X$ is said to satisfy

\noindent ``$\alpha$-mixing'' [or (Rosenblatt) ``strong mixing''] if 
$\alpha_X(n) \to 0$ as $n\to\infty$;

\noindent ``$\beta$-mixing'' (or ``absolute regularity'') if $\beta_X(n) \to 0$ as $n\to\infty$.
\smallskip

The $\alpha$-mixing  condition is due to 
Rosenblatt \cite {ref-journal-Rosenblatt1956}.
The $\beta$-mixing condition was first studied by 
Volkonskii and Rozanov \cite {ref-journal-VolkRoz}, and was attributed there 
to Kolmogorov.
\smallskip 

    By (\ref{eq2.23}), one has that for each positive integer $n$,
\begin{equation}
\label{eq2.33}
0 \leq 2\alpha(n) \leq \beta(n) \leq 1.
\end{equation}
By (\ref{eq2.33}), $\beta$-mixing (absolute regularity) implies $\alpha$-mixing (strong mixing).  
\end{notations}

\begin{remark}
\label{rem2.4}
Now suppose that 
$X:=(X_k, k\in \Z)$ is a strictly stationary {\em Markov chain} (with the random 
variables $X_k$, $k \in \Z$ being real-valued, possibly discrete).
(No assumption of ``reversibility'' yet.)
\smallskip

      As a well known consequence of the Markov property, for each 
positive integer $n$, eqs.\ (\ref{eq2.31})-(\ref{eq2.32}) hold in the following 
augmented forms for the given (strictly stationary) Markov chain $X$:
\begin{align}
\label{eq2.41}
\alpha(n) &= \alpha_X(n) =  
\alpha(\cF^0_{-\infty}, \cF^\infty_n) = \alpha\bigl(\sigma(X_0),\sigma(X_n)\bigl);\\
\label{eq2.42}
\beta(n) &= \beta_X(n) = 
\beta(\cF^0_{-\infty},\cF^\infty_n) = \beta\bigl(\sigma(X_0), \sigma(X_n)\bigl).
\end{align}
(See e.g.\  [\cite {ref-journal-Bradley2007}, Vol.\ 1, Theorem 7.3].)\ \ 
By strict stationarity and (\ref{eq2.41})-(\ref{eq2.42}), one has that for 
any integer $j$ and any positive integer $n$, the (strictly stationary) 
Markov chain $X$ satisfies 
$\alpha(n) = \alpha(\sigma(X_j),\sigma(X_{j+n}))$ and
$\beta(n) = \beta(\sigma(X_j), \sigma(X_{j+n}))$.    
\end{remark}

     {\bf Note.}\ \ In the transcribing of work of Roberts, Rosenthal, and Tweedie in their papers 
\cite {ref-journal-RR} and \cite {ref-journal-RT} on Markov chains 
into the language of (the dependence coefficients for) strong mixing conditions, 
the ``$\rho$-mixing'' condition also comes into play in a significant way.
However, with regard to the results in this paper here, that condition can be tacitly
left in the background; its definition need not be spelled out here.
 
\begin{remark}
\label{rem2.5}
There exist (say real) strictly stationary,
countable-state Markov chains
$X := (X_k, k \in \bbZ)$ and functions
$h: \R \to \R$ satisfying
$E[\, [h(X_0)]^2] < \infty$ and $\var[h(X_0)] > 0$,
such that
$\beta_X(n) \to 0$ at least exponentially fast
as $n \to \infty$
(a condition equivalent to geometric ergodicity, as noted in Section \ref{sc1}),
such that no normalized versions of the partial sums 
$\sum_{k=1}^n h(X_k)$ converge in distribution to $N(0,1)$
as $n \to \infty$.
\smallskip

     That has been shown with constructions by the author 
\cite {ref-journal-Bradley1983b},  
by H\"aggstrom \cite {ref-journal-Haggstrom},
and by Cuny and Lin [\cite {ref-journal-CL2025}, Proposition 3.7].  
The construction by Herrndorf \cite {ref-journal-Herrndorf} also provides 
such an example, though that fact (in particular, the pertinent underlying 
Markov chain in that example) apparently was not brought to light until
it was explained in [\cite {ref-journal-Bradley2007}, Vol.\ 3, Theorem 31.6].
\smallskip

     In the examples in \cite {ref-journal-Bradley1983b} and
\cite {ref-journal-Herrndorf}, the mixing rate for absolute regularity
can be made arbitrarily fast (short of $m$-dependence).
\smallskip

     In the example in \cite {ref-journal-Bradley1983b}
(in Theorem 2 and page 95, Remark 2.1 there),
(i) the partial sums, suitably normalized, converge in distribution along
subsequences to all infinitely divisible laws, and 
(ii) (with a suitable choice of certain parameters) the variances of those 
partial sums can blow up at a rate that is arbitrarily close to 
(but, due to ergodicity, necessarily not quite as fast as) $n^2$.
Further properties of the example in \cite {ref-journal-Bradley1983b} 
were added in [\cite {ref-journal-Bradley2007}, Vol.\ 3, Corollary 31.5].  
Also, an Appendix in \cite {ref-journal-Bradley2024arXiv} (the augmented version
of the paper \cite {ref-journal-Bradley2025} that was posted on arXiv)
provides for the example in \cite {ref-journal-Bradley1983b} as embellished in
[\cite {ref-journal-Bradley2007}, Vol.\ 3, Corollary 31.5],  one particular 
(somewhat long but gentle) explanation of the elementary fact that without any 
changes in other asserted properties, the function $h$ can be made to be one-to-one, 
making the example itself a (strictly stationary, countable-state) Markov chain.   
\end{remark}

\begin{notations}
\label{nt2.6}
(A) A given (say real) strictly stationary Markov chain 
$X := (X_k, k \in \Z)$
is said to be ``reversible'' if the distribution 
[on $(\R^\Z, \cR^\Z)$] of the ``time-reversed''
Markov chain $(X_{-k}, k \in \Z)$ is the same as that
of the sequence $X$ itself.
\smallskip

   (B) By a standard argument, a given (say real) strictly stationary Markov 
chain $X := (X_k, k \in \Z)$ is reversible if and only if
the distribution [on $(\R^2, \cR^2)$] of the random vector
$(X_1, X_0)$ is the same as that of the random vector 
$(X_0, X_1)$.  
\end{notations}

   Refer to the third and fourth paragraphs of Section \ref{sc1}.   
As a ``basic reference point'' for the material in this paper, we shall state a CLT that is part
of [\cite {ref-journal-CL2025}, Theorem 4.9]:

\begin{theorem} {\rm (Cuny and Lin)}.
\label{thm2.7} 
Suppose $X := (X_k, k \in \bbZ)$ is a real strictly stationary, 
{\bf reversible} Markov chain that satisfies $\alpha$-mixing with 
$\alpha_X(n) \to 0$ at least exponentially fast as $n \to \infty$.  
Suppose
$h: \R \to \R$  is a Borel function such that
\begin{equation}
\label {eq2.71}
0\ <\ E[\, [h(X_0)]^2]\ <\ \infty \indent  {\rm and} \indent E[h(X_0)]\ =\ 0.
\end{equation}
Then $\sigma^2 := E[\, [h(X_0)]^2] + 2 \cdot \sum_{k=1}^\infty E[h(X_0) \cdot h(X_k)]$
exists in $(0, \infty)$, with the sum converging absolutely (in fact with the summands
decaying at least exponentially fast); and
$n^{-1} \var(\sum_{k=1}^n h(X_k)) \to \sigma^2$ as $n \to \infty$.
Further, $(1/(n^{1/2}\sigma)) \sum_{k=1}^n h(X_k) \to N(0,1)$ in distribution as $n \to \infty$.  

   Thus, letting $\sigma_n(h) := \|\sum_{k=1}^n h(X_k)\|_2$ for each $n \in \bbN$, 
one has that $[\sigma_n(h)]^{-1} \sum_{k=1}^n h(X_k)$ converges in distribution
to $N(0,1)$ as $n \to \infty$.
\end{theorem}

\begin{remark}
\label{rem2.8}

     (A) There is more to [\cite {ref-journal-CL2025}, Theorem 4.9] than just a CLT as given here.
The last paragraph of Theorem \ref{thm2.7} is redundant,  
but it was included here in order to facilitate a comparison in 
Remark \ref{rem2.9}(B) below with another (``counterexample to CLT'') aspect of
 [\cite {ref-journal-CL2025}, Theorem 4.9].  
\smallskip

     (B) Theorem \ref{thm2.7} is a generalization of an earlier well known CLT, alluded to in the
third paragraph of Section \ref{sc1} --- essentially the same statement but with the stronger mixing assumption of $\beta$-mixing with mixing rate at least exponentially fast (equivalently, 
geometric ergodicity).
\smallskip

     (C) The CLT contained in [\cite {ref-journal-CL2025}, Theorem 4.9] is more general 
than the statement given here in Theorem \ref{thm2.7}.
As was explicitly covered in [\cite {ref-journal-CL2025}, Theorem 4.9], 
letting $P$ denote the Markov chain's one-step ${\cL}_2$ transition operator 
dealt with in \cite {ref-journal-CL2025}, one can replace the assumption of reversibility
in Theorem \ref{thm2.7} (that is, $P$ being self-adjoint) with the weaker assumption
that $P$ is normal in the functional-analysis sense; and one can weaken that assumption 
even further.   
But here we shall restrict the exposition to the case of (strictly stationary) reversible 
Markov chains, the focus of this paper here.
\smallskip

     (D) Here is a quick review of some relevant background material pertaining to
Theorem \ref{thm2.7}. 
By [\cite {ref-journal-Bradley2021}, Theorem 3.3] (the equivalence of conditions 
(R2), (R3), (R4), and (A1) there), a given strictly stationary, reversible Markov chain 
is $\rho$-mixing if and only if it satisfies $\alpha$-mixing with $\alpha(n) \to 0$ at least
exponentially fast as $n \to \infty$, and in that case the mixing rate for $\rho$-mixing
is exponential.
For a given strictly stationary, reversible Markov chain, by
[\cite {ref-journal-CL2025}, the Remark prior to Theorem 3.3], condition (i) in
[\cite {ref-journal-CL2025}, Theorem 4.9] is simply the condition that the Markov chain is
$\rho$-mixing.  
Hence, if a given strictly stationary, reversible Markov chain satisfies $\alpha$-mixing
with $\alpha(n) \to 0$ at least exponentially fast, then condition (i) in
[\cite {ref-journal-CL2025}, Theorem 4.9] is satisfied; and by the equivalence of all three 
conditions (i), (ii), and (iii) in that theorem (note the word ``non-degenerate''
in condition (iii)), one obtains Theorem \ref{thm2.7}.
\smallskip

   (E) For that last comment in (D), a quick word may be in order here with 
regard to how one might formally confirm Theorem \ref{thm2.7} in its entirety from
[\cite {ref-journal-CL2025}, Theorem 4.9].
The sentence right after (\ref{eq2.71}), with its phrase ``exists in $(0, \infty)$'' replaced by
``exists in $[0, \infty)$'', comes from a well known elementary calculation.
Condition (iii) in [\cite {ref-journal-CL2025}, Theorem 4.9] simultaneously yields 
(a) the last paragraph of Theorem \ref{thm2.7} 
and (b) the existence of a {\it positive\/} number $\tau$ such that
$(1/(n^{1/2}\tau)) \sum_{k=1}^n h(X_k) \to N(0, 1)$ in distribution as $n \to \infty$.
Those items (a) and (b) together yield that
$\sigma_n(h) \sim n^{1/2}\tau$ as $n \to \infty$ (see Notations \ref{nt2.1}). 
That and (the final piece of) the sentence right after (\ref{eq2.71}) together yield that
$\sigma^2 = \tau^2\ (> 0)$.
That and item (b) yield the second sentence after (\ref{eq2.71}).  
That completes the argument.      
\end{remark}

\begin{remark}
 \label{rem2.9}
    (A) For the spectrum of results in this paper here, an ``endpoint'' will be a 
known ``basic insight'', loosely stated as follows: 
\smallskip       
 
    {\it For (functions of) strictly stationary, {\bf reversible}, $\alpha$-mixing 
(or $\beta$-mixing) Markov chains, under the assumption of finite second moments,
the ``borderline'' mixing rates for the CLT are the exponential mixing rates.}
\smallskip

     That ``basic insight'' was expressed by 
Cuny and Lin [\cite {ref-journal-CL2025}, Theorem 4.9] 
in a different way, with much more detail,
in a broader context not restricted to the case of reversibility.
For just the case of reversibility, the focus of this paper, it comes from
Theorem \ref{thm2.7} above together with the following item (B):
\medskip

     (B) Suppose that in the statement of Theorem \ref{thm2.7}, one replaces the assumption
of ``exponential $\alpha$-mixing'' by $\alpha$-mixing with any mixing rate that is 
slower than exponential.
Then there will exist a Borel function $h: \R \to \R$ satisfying (\ref{eq2.71})
such that the final paragraph of Theorem \ref{thm2.7} fails to hold.
In light of the information reviewed in Remark \ref{rem2.8}(D) above, 
that was shown by Cuny and Lin [\cite {ref-journal-CL2025}, Theorem 4.9]
(the equivalence of conditions (i) and (ii) there).
With that result together with Theorem \ref{thm2.7} itself 
(in particular, its final paragraph), 
Cuny and Lin [\cite {ref-journal-CL2025}, Theorem 4.9] 
established the ``basic insight'' formulated in (A) above.
\medskip 

     (C) The author \cite {ref-journal-Bradley2025} provided some further insight
into the ``basic insight'' in (A) above, with 
a construction of a class of strictly stationary, countable-state, reversible 
Markov chains with finite second moment (and mean 0) which satisfy $\beta$-mixing 
with mixing rate arbitrarily close to (but not quite) exponential, such that the 
partial sums, suitably normalized, converge in distribution along a subsequence 
to a specific non-degenerate, non-normal law labeled as
$\mu_{P1sL}$ in Notations \ref{nt3.3}(B) in Section \ref{sc3} below.
(The construction there also satisfies a ``dissipation'' property in eq.\ (\ref{eq3.23})
in Section \ref{sc3}, which rules out what was implicitly referred to in 
 \cite {ref-journal-CL2025} as a ``degenerate annealed CLT''.) 
\medskip  
  
    (D) {\tt A redress of an oversight.}\ \ In the manner described in (B) above, 
 the ``basic insight'' in (A) above was
clearly established by Cuny and Lin [\cite {ref-journal-CL2025}, Theorem 4.9] 
at least as early as a preprint of that paper in the year 2023.
At a conference in Cincinnati in May 2024, Michael Lin \cite {ref-journal-Lin} gave a 
talk on that joint work with Christophe Cuny.
The author of this paper here was present at that talk, but failed to comprehend the 
fact that amid the various details, the ``basic insight'' in (A) above was present in that talk.  
In the paper \cite {ref-journal-Bradley2025} referred to in (C) above, starting with its 
first version in November 2024,  the author formulated that ``basic insight'', in the 
same manner as in (A) above, but without mentioning either the paper 
\cite {ref-journal-CL2025} or the talk \cite {ref-journal-Lin}.   
In the paper \cite {ref-journal-Bradley2025}, the references 
\cite {ref-journal-CL2025} (say a preprint of it) and \cite {ref-journal-Lin} 
needed to be cited as having {\it already earlier\/} established 
the ``basic insight'' in (A) above.
\end{remark}

\section{The bounded case --- fast growth of variances}
\label{sc3}

   This section is devoted mainly to the topic of strictly stationary, 
countable-state, irreducible, aperiodic Markov chains that are reversible 
and also bounded and have the following two ``bad'' properties:
(i) the variance of the $n^{\rm th}$ partial sum grows at 
almost the rate $n^2$ as $n \to \infty$; and 
(ii) the Markov chain fails to satisfy a CLT, in the sense that
regardless of what normalizations are used, the normalized partial sums 
fail to converge to a (non-degenerate) normal distribution.

\begin{background}
\label{bckg3.1}
It is well known that, as a result of ``long-run equilibrium'', any strictly stationary,
countable-state, irreducible, aperiodic Markov chain satisfies 
$\beta$-mixing (absolute regularity) and hence [see (\ref{eq2.33})] also $\alpha$-mixing.  
(See e.g.\ [\cite {ref-journal-Bradley2007}, Vol. 1, Theorem 7.7].)\ \ 
However, in that context, the mixing rate for $\beta$-mixing 
(or for $\alpha$-mixing) can be arbitrarily slow. 
(That is true even for the case where the Markov chain is reversible —
see e.g.\ the main result in \cite {ref-journal-Bradley2024}).
That will just be background information for the comments below.
\smallskip

   (A) Davydov [\cite {ref-journal-Davydov1973}, Example 2] 
constructed a class of examples of (non-reversible) real, bounded, 
strictly stationary, countable-state, irreducible, aperiodic
(and hence $\beta$-mixing) Markov chains $(X_k, k \in \Z)$
for which the CLT fails to hold.  
In those counterexamples, the partial sums, suitably normalized, 
converge in distribution to a (symmetric) stable law with index 
(``exponent'') strictly between 1 and 2.
By the Theorem of Types 
(see e.g.\ Billingsley [\cite {ref-journal-Bill}, Theorem 14.2]), there cannot be 
(in those examples) convergence of the partial sums to a (non-degenerate) 
normal law, even along a subsequence, regardless of what normalizing 
constants (multiplicative and additive) are employed.
\smallskip

    (B) With certain parameters in that class of (bounded) examples chosen 
appropriately, Davydov \cite {ref-journal-Davydov1973} 
showed that for any given $\zeps > 0$, such counterexamples exist
in which $\var(\, \sum_{k=1}^n X_k) \geq n^{2-\zeps}$
for all $n \in \N$ sufficiently large.
That is very sharp in that direction.
It is well known (from a simple truncation argument) that for any given strictly 
stationary ergodic random sequence (Markovian or not) with finite second 
moments, the variance of the $n^{th}$ partial sum is $o(n^2)$ as $n \to \infty$.
(Of course ergodicity is implied by, and is in fact considerably weaker than, 
$\beta$-mixing or even $\alpha$-mixing.)
\medskip

   (C) Later, the author [\cite {ref-journal-Bradley1983a}, Theorem 2]
constructed a class of examples of (non-reversible) 
strictly stationary, countable-state, irreducible, aperiodic
(and hence absolutely regular) Markov chains 
$X := (X_k, k \in \Z)$ such that for some
subset $D$ of the state space of $X$,
the partial sums $\sum_{k=1}^n I_D(X_k)$, with suitable normalization, 
converge in distribution along subsequences to all infinitely divisible laws.
[Here $I_D(.)$ denotes the indicator function.]\ \
That class of counterexamples was constructed in such a way that
the rate of growth of $\var [\, \sum_{k=1}^n I_D(X_k)]$
could be arbitrarily close to (though necessarily still strictly 
slower than) $n^2$.
\medskip
                     
   (D) Theorem \ref {thm3.4} below — the main purpose of
Section \ref{sc3} here — will describe a class of real 
strictly stationary, countable-state, irreducible, aperiodic 
(hence $\beta$-mixing) bounded Markov chains that are reversible, 
such that the CLT fails to hold and the variance of 
the $n^{th}$ partial sum grows at almost (but again not quite) the rate $n^2$.
That will extend Davydov's observation in (B) above to the case of
reversible Markov chains.
\medskip

First, in Remark \ref {rem3.2} and Notations \ref{nt3.3} below,
both taken directly from Section 3 of \cite {ref-journal-Bradley2025},
some pertinent comments and notations are given,
in a general form not restricted to bounded random variables.
Those comments and notations will be relevant to Theorem \ref{thm3.4} 
as well as to all other new results in 
(Sections \ref {sc4} and \ref{sc5} of) this paper.          
\end{background}

\begin{remark}
\label{rem3.2}
The new examples constructed in this paper
will be strictly stationary random sequences 
$(X_k, k \in \N)$ (reversible Markov chains) 
that satisfy the conditions
\begin{equation}
\label{eq3.21}
E[X_0^2] < \infty \indent {\rm and} \indent E[X_0] = 0,
\end{equation} 
as well as the condition
\begin{equation}
\label{eq3.22}
\lim_{n \to \infty} n^{-1}
\var\biggl(\, \sum_{k=1}^nX_k\biggl)\ 
=\ \infty
\end{equation}
and even the (stronger) condition
\begin{equation}
\label{eq3.23}
{\rm for\ every}\ c > 0, \quad 
\lim_{n \to \infty} \biggl[\, 
\sup_{r \in \R}\,
P\biggl(r-c\ <\ n^{-1/2}\sum_{k=1}^nX_k\ <\ r+ c\biggl) \biggl]\ =\ 0.
\end{equation}
Of course (\ref {eq3.23}) [in conjunction with (\ref {eq3.21})] implies (\ref {eq3.22}).
If (\ref {eq3.22}) were to fail to hold, then for some infinite set $T \subset \Z$,  
the numbers $n^{-1} \var(\, \sum_{k=1}^nX_k),\ n \in T$ would be bounded, and
the family of distributions of the random variables $n^{-1/2} \sum_{k=1}^nX_k$, $n \in T$
would therefore be tight by (\ref {eq3.21}) and Chebyshev's inequality; 
but that would contradict (\ref {eq3.23}).    
\smallskip

Eq.\ (\ref {eq3.23}) might be referred to as the ``complete dissipation'' of the 
random variables $n^{-1/2}\sum_{k=1}^nX_k$ as $n \to \infty$.  
For a given fixed positive integer $n$, the term in the brackets 
in (\ref {eq3.23}) is, as a function of $c > 0$, a well known 
``concentration function'' of the random variable $n^{-1/2}\sum_{k=1}^nX_k$.
\smallskip

   By a trivial argument, (\ref{eq3.23}) is equivalent to the 
(at first sight seemingly weaker) condition that there {\it exists\/} a 
positive number $c$ such that the equality in (\ref{eq3.23}) holds.
In particular, to verify (\ref{eq3.23}), it suffices to carry out 
the argument for (say) just the case $c = 1$, allowing one 
to slightly simplify the arithmetic.     
\end{remark} 

\begin{notations}
\label{nt3.3}
In the examples constructed in this paper, a certain non-degenerate, 
non-normal distribution, labeled below as ``$\mu_{\rm P1sL}$'', will 
play a role as a ``partial limiting distribution''.  
That distribution is formulated in part (B) below, employing a 
standard definition in part (A).
\medskip 

   (A)\ The Laplace (``double exponential'') distribution is the
absolutely continuous distribution (probability measure on
$(\R, \cR)$) with probability density function 
$f_{\rm Laplace}$ given by 
$f_{\rm Laplace}(x) = (1/2) e^{-|x|}$ for $x \in \R$.
\medskip

   (B)\ Let $\mu_{\rm P1sL}$ denote the distribution [probability 
measure on $(\R, \cR)$] of a random variable $Y$ of the form
\begin{equation}
\label{eq3.31}
Y(\omega)\ :=\ \sum_{k=1}^{N(\omega)} \eta_k(\omega)
\indent {\rm for}\ \omega \in \Omega 
\end{equation}       
where (i) $(\eta_1, \eta_2, \eta_3, \dots)$ is a sequence of independent, 
identically distributed random variables
with the Laplace distribution (see (A) above), and 
(ii) the random variable $N$ is independent of the sequence
$(\eta_1, \eta_2, \eta_3, \dots)$ and is nonnegative 
integer-valued and has the Poisson distribution with mean 1.
\smallskip

   In the subscript of the notation $\mu_{\rm P1sL}$, 
the symbols P, 1, s, and L stand for 
``{\bf P}oisson with mean {\bf 1}'' 
mixture of ``{\bf s}ums of {\bf L}aplace'' random variables.
In (\ref{eq3.31}), it is tacitly understood that for 
$\omega \in \Omega$ such that $N(\omega) = 0$, 
the value $Y(\omega)$ is defined to be 0, with the 
``empty sum'' in the right hand side of that equation
understood to take the value 0.
\end{notations}

   Now refer again to Remark \ref {rem3.2} and Notations \ref{nt3.3}.
Our first result is the following theorem:
\smallskip
 
\begin{theorem}
\label{thm3.4}
Suppose $(q_1, q_2, q_3, \dots)$ is a sequence of positive
numbers such that $q_n \to 0$ as $n \to \infty$.
Then there exists a (real) strictly stationary, countable-state,
irreducible, aperiodic (hence $\beta$-mixing) Markov chain 
$X := (X_k, k \in \Z)$ which is reversible and has
the following properties: \zhfb
(i) There exists a positive constant $C$ such that
$P(|X_0| \leq C) = 1$; \zhfb
(ii) there exists a positive integer $N$ such that 
$\var(\, \sum_{k=1}^n X_k) \geq q_n \cdot n^2$ for all 
$n \geq N$; \zhfb 
(iii) eqs.\ (\ref{eq3.21}), (\ref{eq3.22}), and (\ref{eq3.23})
hold; and \zhfb
(iv) there exists a strictly increasing sequence 
$\bigl(J(1), J(2), J(3), \dots\bigl)$ of positive integers, 
and a sequence $(b_1, b_2, b_3, \dots)$ of positive constants
satisfying $b_n \to \infty$ as $n \to \infty$, such that
\begin{equation}
\label{eq3.41}
b_n^{-1} \sum_{k=1}^{J(n)}X_k\ \to \mu_{\rm P1sL}\ \
{\rm in\ distribution\ as}\ n \to \infty.   
\end{equation}        
\end{theorem}

\begin{remark}
\label{rem3.5}
(A) The statement of Theorem \ref {thm3.4} has a lot of redundancy.
As noted in Remark \ref {rem3.2}, eq.\ (\ref {eq3.23})
[in conjunction with (\ref{eq3.21})] implies (\ref {eq3.22}).
Of course property (i) implies the ``first half'' of 
eq.\ (\ref {eq3.21}).
If $q_n \to 0$ sufficiently slowly, then
property (ii) implies (\ref {eq3.22}).
Eqs.\ (\ref {eq3.41}) and (\ref {eq3.23}) together
force the condition $b_n \to \infty$ (as $n \to \infty$)
in property (iv).
\smallskip
   
   The statements of the other new results in this paper
will have similar redundancies.
\smallskip

   (B) Conclusion (ii) --- in essence the ``almost quadratic''
growth of variances of the partial sums --- is the 
peculiar feature of Theorem \ref{thm3.4}.  
That theorem extends to reversible Markov chains the corresponding 
observation of Davydov in Background \ref {bckg3.1}(B) above.
\smallskip
  	
   (C) Conclusions (iii) and (iv) in Theorem \ref {thm3.4} are
``generic'' to, and will be cited in, the other new results in this paper.
By the Theorem of Types (see e.g.\ [\cite {ref-journal-Bill}, Theorem 14.2]),
conclusion (iv) rules out the possibility of convergence of the normalized partial sums
to a non-degenerate normal law, regardless of what normalizing constants 
(multiplicative and/or additive) are used.
Conclusion (iii) [specifically, eq.\ (\ref{eq3.23})] rules out what 
Cuny and Lin \cite{ref-journal-CL2025} refer to (at least implicitly) as a 
``degenerate annealed CLT'', where the $n^{\rm th}$ partial sum multiplied by 
$n^{-1/2}$ would converge in probability to 0 
[``convergence in distribution to $N(0,0)$''] ---
and that still holds even allowing also additive normalizing constants.
\smallskip

   (D) Eq.\ (\ref{eq3.41}) asserts that $\mu_{P1sL}$ is a ``partial limiting distribution'' 
for suitably normalized partial sums.  
As a consequence, for every positive number $\lambda$,
the partial limiting distributions include that of the sum of a random 
Poisson-mean-$\lambda$ number of independent, identically distributed 
Laplace random variables (with the Poisson random variable being independent 
of the family of Laplace random variables). 
That fact is given directly for $\lambda = 1$ by (\ref{eq3.41}).
If follows (as a consequence of $\alpha$-mixing) for positive integers $\lambda$ 
by basic calculations that go back to Cogburn \cite {ref-journal-Cogburn} and
Ibragimov \cite {ref-journal-Ibragimov1962} 
(see e.g.\ Theorem 12.2 in Volume 1 of  \cite {ref-journal-Bradley2007}).
It then follows for positive rational numbers $\lambda$ by a result of 
Cogburn \cite {ref-journal-Cogburn} 
(see e.g.\ Theorem 12.4 in Volume 1 of  \cite {ref-journal-Bradley2007}).
It then follows for positive real numbers $\lambda$ by an elementary limiting 
argument.
\smallskip

     Of course then by an ordinary CLT, it then follows that the (non-degenerate) normal laws are also ``partial limiting distributions''.
\smallskip

     (E) The Markov chain $X := (X_k, k \in \Z)$ constructed for Theorem 3.4 will in fact have, in addition to reversibility, another symmetry:  the Markov chain
$(-X_k, k \in \Z)$ has the same distribution on $(\R^\Z, \cR^\Z)$ as the Markov chain
$X$ itself.
\smallskip 
          
     Comments (C), (D), and (E) apply as well to the counterexamples described in
Sections \ref{sc4} and \ref{sc5} below as well as to the main result in 
\cite {ref-journal-Bradley2025} (Theorem 3.3 there). 
It would be interesting to know what other partial limiting distributions 
[besides the ones mentioned in (D)] are possible in such counterexamples 
to the CLT involving strictly stationary, reversible Markov chains.
(Compare to the first sentence in the final paragraph of Remark \ref{rem2.5} for the
``non-reversible'' case.)
\end{remark}

   Theorem \ref{thm3.4} will be proved in Section \ref{sc7}, after some 
preliminary work is done in Section \ref{sc6}.

\section{The bounded case again --- very sharp mixing rates}
\label{sc4}

As a ``first point of reference'', we state a classic CLT of Ibragimov 
\cite {ref-journal-Ibragimov1962}.
It is stated here essentially in the form given in 
Ibragimov and Linnik [\cite {ref-journal-IbrLin}, Theorem 18.5.4].
[Conclusion (2) here was not explicitly part of that statement, 
but comes directly from the proof there.] 
 
\begin{theorem} {\rm (Ibragimov)}.    
\label {thm4.1}
Suppose $X := (X_k, k \in \bbZ)$ is a 
(not necessarily Markovian)
strictly stationary sequence of random variables such that
the following two conditions (a), (b) hold: \zhfb
(a) $P(|X_0| \leq C) = 1$ for some 
positive constant $C$, and $E[X_0] = 0$; and \zhfb
(b) $\sum_{n=1}^\infty \alpha_X(n) < \infty$. \zhfb
Then the following three statements hold: \zhfb 
(1) $\sigma^2 := E[X_0^2] + 2\sum_{n=1}^\infty E[X_0X_n]$
exists in $[0, \infty)$, with the sum being absolutely convergent;
\zhfb
(2) $\lim_{n \to \infty} n^{-1} \var(\, \sum_{k=1}^n X_k) = \sigma^2$; and \zhfb
(3) if also $\sigma^2 > 0$, then
$[1/(\sqrt n\, \sigma)] \sum_{k=1}^n X_k$ converges in 
distribution to $N(0,1)$ as $n \to \infty$. 
\end{theorem}

Consider again the class of (bounded) examples of 
Davydov [\cite {ref-journal-Davydov1973}, Example 2]
cited in Background \ref{bckg3.1}(A)(B).
With certain choices of parameters in that construction,
those examples of Davydov show that for any given arbitrarily small fixed 
$\zeps > 0$, Theorem \ref{thm4.1}  would fail to hold if its 
mixing-rate assumption in (b) were replaced by the marginally 
slower rate $\alpha(n) = O(n^{-1+\zeps})$ as $n \to \infty$.   
In those examples, the partial sums, suitably normalized, converge in distribution 
to a (symmetric) stable law with index (``exponent'') strictly between 1 and 2.
Those examples, in conjunction with Theorem \ref {thm4.1}, established 
that for the central limit theorem for strictly stationary, bounded random 
sequences under the strong mixing condition, the ``borderline'' mixing rate is,
at least among power-type mixing rates, the mixing rate $O(n^{-1})$.
\smallskip

     In fact in those counterexamples, as in another set of examples in
[\cite {ref-journal-Davydov1973}, Example 1] in the same paper
that will be cited in Section \ref{sc5}, 
Davydov showed implicitly, with calculations on pages 327-328
of that paper, that the $\beta$-mixing condition 
holds with a mixing rate that is within a constant factor of being
the same rate as for $\alpha$-mixing.
\smallskip

     Later, a (nontrivial) strictly stationary, finite-state, non-Markovian counterexample
(to the CLT), satisfying $\beta$-mixing with that ``borderline'' mixing 
rate $O(n^{-1})$ itself, was constructed in \cite {ref-journal-Bradley1989}. 
\smallskip

     To get some further perspective on this, we state here a CLT of 
Merlev\`ede and Peligrad [\cite {ref-journal-MP} Corollary 1.1].
It is not quite the sharpest known CLT of this kind (more on that in Remark \ref{rem4.6}
below), but it is still very sharp; and because of the relative simplicity of its statement, 
it is well suited for comparison with a class of counterexamples described in
Theorem \ref{thm4.4} below.  

\begin{theorem} {\rm (Merlev\`ede and Peligrad)}.    
\label {thm4.2}
Suppose $X := (X_k, k \in \Z)$ is a 
(not necessarily Markovian)
strictly stationary sequence of random variables such that
the following three conditions (a), (b), (c) hold: \zhfb
(a) $P(|X_0| \leq C) = 1$ for some 
positive constant $C$, and $E[X_0] = 0$; \zhfb 
(b) $\liminf_{n \to \infty} n^{-1} 
\var(\, \sum_{k=1}^n X_k) > 0$; and \zhfb
(c) $\alpha_X(n) = o(n^{-1})$ as $n \to \infty$. \zhfb
Then letting $b_n := (\pi/2)^{1/2} E|\sum_{k=1}^nX_k|$
for each $n \in \N$, one has that 
$b_n \to \infty$ as $n \to \infty$, and   
$(1/b_n) \sum_{k=1}^n X_k$ converges in distribution to
$N(0,1)$ as $n \to \infty$.
\end{theorem}

     This is actually a special case of Theorem \ref{thm5.3} in the next
section.  For checking the literature, the comment right after the 
statement of that theorem is relevant.

\begin{remark}
\label{rem4.3}
[In comments (B) and (C) below, for checking the literature, the 
notation and comments in Notations \ref{nt5.1} in the next section 
will be relevant.]
\medskip

   (A) Theorem \ref {thm4.1} (in the case $\sigma^2 > 0$) and 
Theorem \ref {thm4.2} each involve strictly stationary sequences that satisfy 
$\liminf_{n \to \infty} n^{-1} \var(\, \sum_{k=1}^n X_k) > 0$.
The ``trade-off'' between those two theorems is essentially
as follows:   (i) In Theorem \ref {thm4.1} the normalizing constant 
(for convergence to a normal law) is of a 
``standard'' type (that is, $\asymp 1/\sqrt n$\ ), 
while in Theorem \ref {thm4.2} that is not necessarily the case.
(ii) On the other hand, Theorem \ref {thm4.2} involves
a marginally slower mixing rate assumption than in
Theorem \ref {thm4.1}. 
\medskip    

   (B) In its own way, Theorem \ref {thm4.1} is very sharp. 
If its (summable) mixing rate assumption were weakened 
at all (without any change in the other assumptions), 
then its conclusion as stated could fail to hold. 
That was shown with a class of counterexamples 
--- strictly stationary, non-Markovian sequences --- 
in [\cite {ref-journal-Bradley1997}, Theorem 4 
with $q(u) \equiv 1$ there],
satisfying the ``complete dissipation'' property
(\ref {eq3.23}) here in Section \ref {sc3}. 
\smallskip

   Those counterexamples violate the conclusion --- and hence 
also the mixing-rate assumption --- of Theorem \ref{thm4.1}. 
However, some of them  --- for 
example, one with $q(u) \equiv 1$ and 
$h(n) = n^{-1}(\log n)^{-1}$ (for $n \geq 2$) there --- still
satisfy the hypothesis of Theorem \ref{thm4.2},
and their partial sums will therefore, when normalized as in 
Theorem \ref{thm4.2}, converge to $N(0,1)$
in distribution. 
\medskip 

   (C) In its own way,
Theorem \ref {thm4.2} is very sharp. 
That was shown with a class of counterexamples 
--- again, strictly stationary, non-Markovian sequences --- 
in [\cite {ref-journal-Bradley1997}, Theorem 3, 
with $q(u) \equiv 1$ and $h(n) = n^{-1}$ there],
satisfying (a) and (b) in Theorem \ref{thm4.2} as well as
$\alpha(n) \leq \beta(n) \leq n^{-1}$ for all
$n \in \N$, and having the weird property
(adapted from the construction of Herrndorf
\cite {ref-journal-Herrndorf}) that 
$\inf_{n \in \N} P(\, \sum_{k=1}^nX_k = 0) > 0$.
(But see also Remark \ref{rem4.6} below.)
\medskip

   (D) Theorem \ref {thm4.4} below will show that 
even for strictly stationary, countable-state, reversible 
Markov chains, Theorem \ref {thm4.2} is still quite sharp.   
The mixing rate property declared in Theorem \ref {thm4.4} below 
is very close to the rate $O(n^{-1})$ fulfilled in the counterexamples 
alluded to in comment (C) above
[since the numbers $g_n,\, n \in \N$ in Theorem \ref{thm4.4} can be 
taken to diverge to $\infty$ arbitrarily slowly].\ \   
Remark \ref {rem3.5}(A) regarding redundancy applies here.      
\end{remark}

\begin{theorem}
\label{thm4.4}
Suppose $(g_1, g_2, g_3, \dots)$ is a sequence of positive
numbers such that $g_n \to \infty$ as $n \to \infty$.
Then there exists a (real) strictly stationary, countable-state,
irreducible, aperiodic Markov chain 
$X := (X_k, k \in \Z)$ which is reversible and has
the following properties: \zhfb
(i) There exists a positive constant $C$ such that
$P(|X_0| \leq C) = 1$; \zhfb
(ii) $\alpha(n) \leq \beta(n) << g_n \cdot n^{-1}$ as $n \to \infty$ 
(see Notations \ref{nt2.1});  \zhfb
(iii) conclusions (iii) and (iv) in Theorem \ref {thm3.4} hold; and \zhfb
(iv) defining ${\bf h}(n) := n^{-1} \var(\sum_{k=1}^n X_k)$ for every $n \in \N$, one has that
\begin{align}
\label{eq4.41}
&{\bf h}(n)\ \leq\ g_n\ \ {\rm for\ every}\ n \in \N;\ \ \ {\rm and} \\
\label{eq4.42}
&{\rm for\ every}\ N \in \N,\ \ \ 
\limsup_{n \to \infty} \frac {{\bf h}(N n)} {{\bf h}(n)}\ =\ N. 
\end{align}       
\end{theorem}
  
   Theorem \ref{thm4.4} will be proved in Section \ref{sc8}.
   
\begin{remark} 
\label{rem4.5}
    As described in Remark \ref{rem3.5}(C), conclusion (iii) in Theorem \ref{thm4.4}
rules out a possible CLT (even a ``degenerate annealed CLT'').
In Theorem \ref{thm4.4}, the purpose of its conclusion (iv) [in particular, eq.\ (\ref{eq4.42})]
is to provide some ``supporting insight'' with regard to how it is (in this counterexample)
that a CLT fails to hold.
Note that in Theorem \ref{thm3.4} in the previous section, one can trivially ensure that
eq.\ (\ref{eq4.42}) holds [but with no useful analog of (\ref{eq4.41})], 
by simply choosing the numbers $(q_n, n \in \N)$ in that theorem to 
converge to 0 sufficiently slowly [for example, such that $1/(\log n) << q_n$ as $n \to \infty$].
Here are some more comments on (\ref {eq4.42}) [together with (\ref{eq4.41})]
in connection with the rest of Theorem \ref{thm4.4} here.
(These comments will apply to the counterexamples described in Section \ref{sc5} as well.)   
\smallskip      

     (A) In (iv) the boldface letter $\bf h$ is employed in order to avoid confusion with the use of
the ordinary letter $h$ as a mathematical symbol elsewhere in this paper.
(A similar use of a few other boldface letters will appear elsewhere in this paper.)
\smallskip

     (B) As a complement to (\ref{eq4.41}), because of eq.\ (\ref{eq3.22})
[from conclusion (iii) in Theorem \ref{thm3.4}, 
as part of conclusion (iii) in Theorem \ref{thm4.4}], 
the numbers ${\bf h}(n)$ are bounded below by a positive constant.
\smallskip

     (C) Of course by the assumptions on the sequence 
$(g_n,\ n \in N)$, one has that $\inf_{n \in N} g_n > 0$.
With an elementary argument, without loss of generality, one can assume 
that the sequence $(g_n)$ is monotonically increasing and slowly varying.   
The issue of slowly varying sequences will come up again in (D) below 
and in Remark \ref{rem4.6}. 
\smallskip
  
     (D) In eq.\ (\ref{eq4.42}), for any given $N \in \N$, $N$ itself is the maximum possible
value for the right hand side, by e.g.\ a simple application of Minkowski's inequality.
Even if (\ref {eq4.41}) holds with the sequence ($g_n$) being slowly varying,
that does not prevent the sequence (${\bf h}(n)$) from severely {\it failing\/} to be 
slowly varying, a failure manifested in (\ref{eq4.42}).
\smallskip

     (E) The peculiar behavior of the variances of the partial sums $\sum_{k=1}^n X_k$ 
manifested in (\ref{eq4.42}) is a fundamental part of what makes the Markov chain $X$ fail to satisfy a CLT [as described in Remark \ref{rem3.5}(C)].
Roughly speaking, in Theorem \ref{thm4.4} (as well as in Theorem \ref{thm3.4} and in the 
related counterexamples described in Section \ref{sc5}),
the Markov chain $X$ is constructed from a countably infinite family 
of ``building blocks'' ---  each a strictly stationary, reversible, 
three-state Markov chain with a particular structure.  
(Those ``building blocks'' are independent of each other.)\ \ 
 Eq.\ (\ref{eq4.42}) is (indirectly) a fundamental part of a scheme that (as $n \to \infty$) 
 allows each of these ``building block'' Markov chains, in turn, to have 
 its ``day in the sun'', when 
(i) with regard to partial sums, it ``dominates'' all of the other ``building block'' Markov chains together, and 
(ii) its own partial sums are not yet close to being normally distributed, but are (under suitable normalization) close to having the (non-degenerate, non-normal) distribution 
$\mu_{P1sL}$ in Notations \ref{nt3.3}.
\end{remark}
   
\begin{remark}
\label{rem4.6} (A) Theorem \ref{thm4.2} above is a special 
case of Theorem \ref{thm5.3} in the next section.
Theorem \ref{thm4.2} served as a convenient ``basic reference point'' for the 
counterexamples described in Theorem \ref{thm4.4}.
Likewise, in the next section, 
Theorem \ref{thm5.3} will serve as a convenient ``basic reference point'' for the 
counterexamples described in Theorem \ref{thm5.5} and in a couple of corollaries. 
However, several years after the original publication of Theorems \ref{thm4.2} and \ref{thm5.3}
in \cite{ref-journal-MP}, those theorems were
slightly sharpened with a CLT by the same two authors, 
Merlev\`ede and Peligrad, in [\cite {ref-journal-MP2}, Theorem 4].
(Their result there was a functional CLT.)\ \ 
A natural question is what is the comparison between
that latter CLT and the counterexamples constructed here in this paper.
In the comments below, we shall try to crystalize that question a bit. 
\smallskip

      (B) In fact the assumptions (b) and (c) in Theorem \ref{thm5.3} in the next section are,
in the special case of bounded random variables, 
(equivalent to) assumptions (b) and (c) in Theorem \ref{thm4.2} above.  
(Indeed the assumptions (b) in both places are identical.)\ \  
In both the general context and the bounded context, the pair of assumptions (b) and (c)
is replaced in the newer result in [\cite {ref-journal-MP2}, Theorem 4] by
a single assumption --- eq.\ (3.2) in that theorem --- which is implied by, 
and is strictly weaker (hence more general and more versatile) than, 
the pair of assumptions (b) and (c).
That refined assumption [\cite {ref-journal-MP2}, eq.\ (3.2)] is intricate and involves 
the quantiles in Notations \ref{nt5.1} in the next section.
An explicit formulation of it will be postponed to Section \ref{sc10}, 
where it will be examined more closely.
But a few more comments on both [\cite {ref-journal-MP2} Theorem 4] 
and its refined assumption (3.2) are in order here. 
\smallskip

       (C) That assumption [\cite {ref-journal-MP2}, eq.\ (3.2)] involves the rate of growth of
the variances of the partial sums [in a more flexible way than in assumption (b) in Theorem
\ref{thm4.2}] together with mixing rates (for a variant of $\alpha$-mixing
employed in \cite {ref-journal-MP2}) and the ``moment'' or (more generally) ``tail'' properties 
of the random variables. 
It applies to some cases where the mixing-rate 
assumption (c) in Theorem \ref{thm4.2} is (possibly) not satisfied.
For example, in their paper on CLTs for (functions of) reversible Markov chains, 
Zhao et al.\ \cite {ref-journal-ZWV} pointed out how, in their Example 1, involving an explicit bounded functional of a particular reversible Markov chain that 
satisfies $\beta(n) \asymp 1/n$, one could derive the asymptotic 
normality of the partial sums by applying that CLT of 
Merlev\`ede and Peligrad [\cite {ref-journal-MP2}, Theorem 4].
\smallskip

   (D) Theorem \ref{thm4.4} and Corollary \ref{cor5.7} (in Section \ref{sc5}) 
seem to be the two results in this paper here that involve counterexamples in which 
the ``borderline'' of the CLT --- involving ``essentially sharpest possible'' combinations of properties of the marginal distribution together with mixing rates for $\alpha$-mixing ---
seems to be very nearly the same with reversibility as it is without reversibility.
The question posed in (A) above seems to be most pertinent to those two results, and
it seems to be a question of    
how the refined assumption [\cite {ref-journal-MP2}, eq.\ (3.2)]
compares to related properties of the constructions for those two results.
In Section \ref{sc10}, we shall give (essentially) the formulation of the refined assumption
[\cite {ref-journal-MP2}, eq.\ (3.2)], and then illustrate directly how in the construction in 
Section \ref{sc8} for Theorem \ref{thm4.4}, involving the relatively simple case of 
bounded random variables, it fails to hold.
(We shall not bother to carry out such an illustration for Corollary \ref{cor5.7}.)
\smallskip    

   (E) The question posed in (A) and (D) is connected to a more basic issue: 
a stark contrast between the behavior of the variances of the partial sums in the CLT 
in [\cite {ref-journal-MP2}, Theorem 4] and the corresponding behavior 
in the counterexamples constructed in this paper here. 
Both of the papers \cite{ref-journal-MP} and \cite{ref-journal-MP2} involved the following
condition for the numbers ${\bf h}(n):= n^{-1} \var (\sum_{k=1}^n X_k),\ n \in \N$ defined in conclusion (iv) of Theorem \ref{thm4.4}:
\begin{equation}
\label{eq4.61}
{\bf h}(n)\ {\rm is\ slowly\ varying\ (in\ the\ strong\ sense)\ as}\ n \to \infty.  
\end{equation}
That condition is common in central limit theory for strictly stationary sequences
under strong mixing conditions.
In the proof in \cite{ref-journal-MP} of Theorem \ref{thm5.3} (and as a special case,
Theorem \ref{thm4.2}), eq.\ (\ref{eq4.61}) here was verified as an intermediate step
in the proof.
Similarly, as part of the proof of [\cite{ref-journal-MP2}, Theorem 4],  
its eq.\ (3.2) was employed in order to verify eq.\ (\ref{eq4.61}) here as an intermediate step.
Eq.\ (\ref{eq4.61}) is in stark contrast to eq.\  (\ref{eq4.42}) in Theorem \ref{thm4.4} 
and (implicitly) in Theorem \ref{thm5.5} (and its Corollaries \ref{cor5.7} and \ref{cor5.8})
in Section \ref{sc5}.
For the significance of that, recall Remark \ref{rem4.5}(E).
Recall also the third sentence of Remark \ref{rem4.5}, in connection with
Theorem \ref{thm3.4}.         
\end{remark}

\section{The unbounded case --- with quantiles}
\label{sc5}

   Now let us turn to the case of random variables that are not bounded.
\smallskip

    In addition to Theorem \ref{thm4.1} (which involved bounded random variables), 
Ibragimov \cite {ref-journal-Ibragimov1962} proved an analog for some 
classes of unbounded random variables: namely, 
a CLT for strictly stationary sequences $(X_k, k \in \Z)$ 
such that for some $\delta > 0$,
$E[|X_0|^{2 + \delta}] < \infty$ and
$\sum_{n=1}^\infty [\alpha(n)]^{\delta / (2 + \delta)} < \infty$.  
A statement of that theorem can also be found in Ibragimov 
and Linnik [\cite {ref-journal-IbrLin}, Theorem 18.5.3].
Davydov [\cite {ref-journal-Davydov1973}, Example 1]
constructed a class of examples that showed that that result of 
Ibragimov is very sharp, particularly with regard to power-type mixing rates,
even under $\beta$-mixing (see the second paragraph
after Theorem \ref{thm4.1} in the preceding section).
\smallskip

   Later, with the use of quantiles instead of moments,
Doukhan et al.\ [\cite {ref-journal-DMR1994}, Theorem 1] obtained a 
CLT that is slightly sharper --- in fact in a certain sense as sharp as possible.
It will be stated here in Theorem \ref {thm5.2} below,
after some relevant terminology is reviewed in Notations \ref {nt5.1}.

\begin{notations}
\label{nt5.1}
Here it will suffice to restrict to nonnegative random variables.
For any nonnegative random variable $W$, define the
``upper tail'' quantile function $Q_W: (0,1) \to [0, \infty)$
as follows:  For each $u \in (0,1)$,
\begin{equation}
\label{eq5.11}
Q_W(u)\ :=\ \inf\bigl\{t \geq 0: P(W > t) \leq u\bigl\}.
\end{equation}
(This is different from, but closely related to, another standard
definition of quantile function.)\ \
It is easy to see that on the entire open unit interval $(0,1)$, this function 
$Q_W$ is nonnegative, nonincreasing, and right-continuous 
(continuous from the right).
\smallskip 
  
   By an elementary (if slightly tricky) well known argument,
if $U$ is a random variable uniformly distributed on the unit interval $(0,1)$, 
then the random variable $Q_W(U)$ has the same distribution as the 
random variable $W$ itself, and hence for (say) any real nonnegative 
Borel function $f$ on $[0, \infty)$, one has that
\begin{equation}
\label{eq5.12}
E[f(W)]\ =\ \int_0^1 f\bigl(Q_W(u)\bigl)\, du. 
\end{equation}
(The common value can be $\infty$.)\ \ 
Eq.\ (\ref{eq5.12}) with the function $f(x) := x^2$ 
will be particularly pertinent
(say with $W := |X|$ for a given real-valued random variable $X$). 
\smallskip

   In the rest of this paper, for a given random variable $X_0$
(in, say, a strictly stationary sequence $(X_k, k \in \Z)$),   
the quantile function of the nonnegative random variable
$|X_0|$ will be written as $Q_{|X(0)|}$ for typographical convenience.   
\end{notations}

   The CLT of Doukhan et al.\ [\cite {ref-journal-DMR1994}, Theorem 1] 
alluded to above is as follows:   

\begin{theorem} {\rm (Doukhan, Massart, and Rio)}.    
\label {thm5.2}
Suppose $X := (X_k, k \in \Z)$ is a (not necessarily Markovian)
strictly stationary sequence of random variables such that
the following conditions (a), (b) hold: \zhfb
(a) $E[X_0^2] < \infty$ and $E[X_0] = 0$; and \zhfb
(b) $\sum_{n=1}^\infty 
\int_0^{\alpha(n)} Q_{|X(0)|}^2(u)\, du < \infty$. \zhfb
Then the following three statements hold: \zhfb 
(1) $\sigma^2 := E[X_0] + 2\sum_{n=1}^\infty E[X_0X_n]$
exists in $[0, \infty)$, with the sum being absolutely convergent;
\zhfb
(2) $\lim_{n \to \infty} n^{-1} \var(\, \sum_{k=1}^n X_k) = \sigma^2$; and \zhfb
(3) if also $\sigma^2 > 0$, then
$[1/(\sqrt n\, \sigma)] \sum_{k=1}^n X_k$ converges in 
distribution to $N(0,1)$ as $n \to \infty$. 
\end{theorem}

   Conclusions (1) and (2) are due to Rio 
[\cite {ref-journal-Rio1993}, Theorem 1.2], and are based
on his covariance inequality (in that paper) involving quantiles 
and the measure of dependence $\alpha(\cdot \, ,\, \cdot)$.
Conclusion (3) is due to Doukhan et al.\
[\cite {ref-journal-DMR1994}, Theorem 1].
Their result was actually a functional CLT.
\smallskip

   The formulation of condition (b) here in Theorem \ref{thm5.2} is
different from, but equivalent to, the formulation of it in 
\cite {ref-journal-DMR1994}.
(See \cite {ref-journal-DMR1994} or \cite {ref-journal-Rio1993} or
the next-to-last paragraph on page 327 of Vol.\ 1 of
\cite {ref-journal-Bradley2007}.)
\smallskip 

   The following CLT is due to Merlev\`ede and Peligrad
[\cite {ref-journal-MP} Theorem 1.3].
(It includes Theorem \ref{thm4.2} in the previous section as a special case.  
Their result was actually a functional CLT.)\ \ 
But again recall Remark \ref{rem4.6}. 

\begin{theorem} {\rm (Merlev\`ede and Peligrad)}.    
\label {thm5.3}
Suppose $X := (X_k, k \in \bbZ)$ is a 
(not necessarily Markovian)
strictly stationary sequence of random variables such that
the following three conditions (a), (b), (c) hold: \zhfb
(a) $E[X_0^2] < \infty$ and $E[X_0] = 0$; \zhfb
(b) $\liminf_{n \to \infty} n^{-1} 
\var(\, \sum_{k=1}^n X_k) > 0$; and \zhfb
(c) $\int_0^{\alpha(n)} Q_{|X(0)|}^2(u)\, du 
= o(n^{-1})$ as $n \to \infty$. \zhfb
Then letting $b_n := (\pi/2)^{1/2} E|\sum_{k=1}^nX_k|$
for each $n \in \N$, one has that
$b_n \to \infty$ as $n \to \infty$, and 
$(1/b_n) \sum_{k=1}^n X_k$ converges in distribution to
$N(0,1)$ as $n \to \infty$.
\end{theorem}

   (In that paper \cite {ref-journal-MP}, the fact that $b_n \to \infty$ as $n \to \infty$ 
was not explicitly built into the statement of that theorem; but with a careful 
examination, one can see that it was shown in the proof there.
Alternatively, for a slightly more generously detailed presentation of their 
theorem and their proof (including an explicit highlighting of the fact that 
$b_n \to \infty$ as $n \to \infty$), see e.g.\ 
[\cite {ref-journal-Bradley2007}, Vol.\ 2, Theorem 16.18].)

\begin{remark}
\label{rem5.4}
(A) Theorem \ref {thm5.2} (in the case $\sigma^2 > 0$) and 
Theorem \ref {thm5.3} each involve strictly stationary sequences that satisfy 
$\liminf_{n \to \infty} n^{-1} \var(\, \sum_{k=1}^n X_k) > 0$.
The ``trade-off'' between those two theorems is essentially as follows:   
(i) In Theorem \ref {thm5.2} the normalizing constant (for convergence to a 
normal law) is of a ``standard'' type (that is, $\asymp 1/\sqrt n$\ ), 
while in Theorem \ref {thm5.3} that is not necessarily the case.
(ii) On the other hand, Theorem \ref {thm5.3} involves
a marginally slower mixing rate assumption than in
Theorem \ref {thm5.2}. 
\medskip    

   (B) In its own way, Theorem \ref {thm5.2} is very sharp.
If its mixing rate assumption were weakened 
 at all (without any other change in the assumptions), 
then its conclusion as stated could fail to hold.
\smallskip

For mixing rates of the ``power type'' form
$O(n^{-a})$ where $a > 1$, that was shown with a class
of counterexamples --- strictly stationary, reversible
Markov chains --- constructed by 
Doukhan et al.\ [\cite {ref-journal-DMR1994},
Theorems 2 and 5 and Corollary 1].
(The property of ``reversibility'' was apparently not mentioned 
explicitly there.)
\smallskip
  
   For more general mixing rates, that was also shown with a class of 
counterexamples --- strictly stationary, non-Markovian sequences --- 
in [\cite {ref-journal-Bradley1997}, Theorem 4], 
with the ``complete dissipation'' property
(\ref {eq3.23}) in Section \ref {sc3} explicitly stated there. 
Those counterexamples violate the conclusion --- and hence 
also the mixing-rate assumption --- of Theorem \ref{thm5.2}. 
However, some of them still satisfy the hypothesis of 
Theorem \ref{thm5.3}, and their partial sums will therefore, when 
normalized as in Theorem \ref{thm5.3}, converge to $N(0,1)$ in distribution. 
\medskip 

   (C) In its own way, Theorem \ref {thm5.3} is very sharp.
That was shown with a class of counterexamples 
--- again, strictly stationary, non-Markovian sequences --- 
in [\cite {ref-journal-Bradley1997}, Theorem 3] that
satisfy conditions (a) and (b) in Theorem \ref{thm5.3} as well as
$\int_0^{\alpha(n)} Q_{|X(0)|}^2(u)\, du \leq n^{-1}$ for all $n \in \N$,
and have (as in the context of Remark \ref{rem4.3}(C))
the weird property (again adapted from \cite {ref-journal-Herrndorf}) that 
$\inf_{n \in \N} P(\, \sum_{k=1}^nX_k = 0) > 0$.
But again, see Remark \ref{rem4.6}, regarding the refinement of Theorem \ref{thm5.3}
that was given in [\cite {ref-journal-MP2}, Theorem 4].
\medskip

   (D) The next theorem, primarily involving strictly stationary, reversible Markov 
chains whose random variables are unbounded, is intended in essence as a 
counterpart to Theorem \ref{thm4.4} (which involved the case of bounded 
random variables).    
\end{remark}

\begin{theorem}
\label{thm5.5}
   Suppose $f$ is a real function on the closed half line 
$[1, \infty)$ such that the following conditions are
satisfied: \zhfb
(a) $0 < f(x) \leq 1$ for all $x \in [1, \infty)$; \zhfb  
(b) the function $f$ is continuous and strictly decreasing on $[1, \infty)$,
and $\lim_{x \to \infty} f(x) = 0$; \zhfb
(c) the function $f$ is twice continuously differentiable on the
open half line $(1, \infty)$; \zhfb
(d) for every positive number $c$, one has that
$\lim_{x \to \infty} [e^{cx} f(x)] = \infty$; \zhfb
(e) there exists ${\bf w} \in (1,\infty)$ such that  
the function $x \mapsto x \cdot f(x)$ for $x \in [1, \infty)$ 
is nonincreasing on $[{\bf w}, \infty)$; \zhfb 
(f) the function $x \mapsto \log f(x)$ for 
$x \in [1, \infty)$ is convex. 
\smallskip 
   
   Suppose $g:[1, \infty) \to [1, \infty)$ is a function such that
$\lim_{x \to \infty} g(x) = \infty$.
\smallskip 
  
   Then there exists a (real) strictly stationary, countable-state,
irreducible, aperiodic Markov chain 
$X := (X_k, k \in \Z)$ which is reversible,
such that $E[X_0^2] < \infty$ and $E[X_0] = 0$ and
(see Notations \ref{nt2.1} and \ref{nt5.1})
\begin{align}
\label{eq5.51}   
    &\int_0^{f(x)} Q_{|X(0)|}^2(u)\, du\  <<\ 
\frac {-f'(x)} {f(x)} \cdot g(x)\ \ 
{\rm as}\ x \to \infty, \quad  {\rm and}\\
\label{eq5.52}   
    &\alpha_X(n)\ \leq\ \beta_X(n)\ \leq\ f(n)\ \
    {\rm for\ all}\ n \in \N,  
\end{align}
and conclusions (iii) and (iv) in Theorem \ref{thm3.4} hold,
and conclusion (iv) in Theorem \ref{thm4.4} holds as well
[with $g_n = g(n)$].     
\end{theorem}

\begin{remark}
\label{rem5.6}
Theorem \ref{thm5.5} will be proved in Section \ref{sc9}.
Its primary interest is with unbounded random variables.
(Theorem \ref{thm4.4} in Section \ref{sc4} gives in essence 
an analog for bounded random variables.)\ \
The properties asserted in Theorem \ref{thm5.5} 
include (indirectly) eqs.\ (\ref{eq3.23}) and (\ref{eq3.41}), 
which rule out a CLT, as described in Remark \ref{rem3.5}(C),
as well as eqs.\ (\ref{eq4.41}) [with $g_n = g(n)$] and (\ref{eq4.42}), which provide
some ``supporting insight'' regarding the failure of the CLT. 
Here are some other comments on Theorem \ref{thm5.5}:
\smallskip

(A) Of course the positive function $g(x),\ x \in [1, \infty)$ can diverge to $\infty$ 
arbitrarily slowly (as $x \to \infty$).
\smallskip
 
(B) In connection with the integral in (\ref{eq5.51}) and the assumption
here that $0 < f(x) \leq 1$ for all $x \in [1, \infty)$,  
one has that 
\begin{equation}
\label{eq5.61} 
\int_0^1 Q_{|X(0)|}^2(u)\, du\ =\ E(X_0^2)\ < \ \infty,
\end{equation}
where the equality comes from (\ref{eq5.12}) and the inequality
$E(X_0^2) < \infty$ was asserted just before (\ref{eq5.51}).
\smallskip

(C) Refer to hypotheses (a) and (b) in the theorem. 
One has that $\lim_{x \to \infty} \log(f(x)) = -\infty$, and also that
$\log(f(x))$ is strictly decreasing on $[1,\infty)$.
 Hence by hypotheses (c) and (f), one has that at each
$x$ in the open half line $(1, \infty)$, the derivative of
of $\log(f(x))$ [namely $f'(x)/f(x)$] is negative
and hence [trivially by hypothesis (a)] $f'(x)$ is negative.
Thus the right hand side of (\ref{eq5.51}) is positive for every $x \in (1, \infty)$.
\smallskip

(D) From hypotheses (b) and (d), the ``bound-on-mixing-rate'' function $f(\, .\, \,)$
itself approaches 0 more slowly than any exponential decay.
\smallskip 

(E) Refer to Remark \ref{rem5.4}(B) and the construction in 
\cite {ref-journal-DMR1994} referred to there. 
Theorem \ref{thm5.5} is an attempt to get 
further information about the ``borderline'' of the 
CLT for strictly stationary, {\it reversible\/} Markov chains,
in terms of combinations of 
(i) mixing rates for $\alpha$-mixing and/or $\beta$-mixing, and 
(ii) ``moment'' assumptions --- or rather, more general ``tail'' assumptions, 
expressed through quantiles --- on the marginal distribution.
More on that in Corollaries \ref{cor5.7} and \ref{cor5.8} and 
Remarks 5.9 and \ref{rem5.10} below.     
\smallskip

(F) In order to simplify some arguments in the proof (in Section \ref{sc9})
of Theorem \ref{thm5.5}, the assumptions in that theorem are stronger than 
necessary and have some redundancies.
However, those assumptions are compatible with ordinary mixing-rate 
functions such as the power-type functions $1/x^p$ for $p \geq 1$,
or the sub-exponential rates such as $\exp(-x^q)$ for $0 < q < 1$.
(See Corollaries \ref{cor5.7} and \ref{cor5.8} below.)  
\smallskip

(G) One obvious shortcoming in Theorem \ref{thm5.5} is
that a quantile function $Q_{|X(0)|}$ --- or alternatively
a function $q: (0,1) \to [0, \infty)$ with appropriate 
properties (depending on the function $f$) such that
$Q_{|X(0)|}(u) \leq q(u)$ for all $u \in (0,1)$ ---
has not been specified {\it beforehand\/}.
The effective specification beforehand of such functions seems quite feasible for 
mixing rates that are of power type or sufficiently close to that, but possibly 
rather cumbersome or ungainly for faster 
mixing rates such as the sub-exponential rates in Corollary \ref{cor5.8} below.
That issue will not be pursued further here.   
\end{remark}

   Here are the statements of Corollaries \ref{cor5.7} and \ref{cor5.8}.
Corollary \ref{cor5.7} holds as a direct application of Theorem \ref{thm5.5} 
with the function $f(x) := (1/x)^p$ (with $p > 1$) for $x \in [1,\infty)$.
Corollary \ref{cor5.8} will perhaps warrant a slight explanation.

\begin{corollary}
\label{cor5.7}
Suppose $1 < p < \infty$.

Suppose $g:[1, \infty) \to [1, \infty)$ is a function such that
$\lim_{x \to \infty} g(x) = \infty$. 
  
   Then there exists a (real) strictly stationary, countable-state,
irreducible, aperiodic Markov chain 
$X := (X_k, k \in \Z)$ which is reversible, such that
$E[X_0^2]< \infty$ and $E[X_0] = 0$ and 
\begin{align}
\label{eq5.71}   
    &\int_0^{(1/x)^p} Q_{|X(0)|}^2(u)\, du\ <<\ 
(1/x) \cdot g(x)\ \ {\rm as}\ x \to \infty, 
\quad  {\rm and}\\
\label{eq5.72}   
    &\alpha_X(n)\ \leq\ \beta_X(n)\ 
    \leq\ (1/n)^p\ \
    {\rm for\ all}\ n \in \N,  
\end{align}
and conclusions (iii) and (iv) in Theorem \ref{thm3.4} hold,
and conclusion (iv) in Theorem \ref{thm4.4} holds as well [with $g_n = g(n)$]. 
\end{corollary}

\begin{corollary}
\label{cor5.8}
Suppose $0 < q < 1$.

Suppose $g:[1, \infty) \to [1, \infty)$ is a function such that
$\lim_{x \to \infty} g(x) = \infty$. 
  
   Then there exists a (real) strictly stationary, countable-state,
irreducible, aperiodic Markov chain 
$X := (X_k, k \in \Z)$ which is reversible, such that
$E[X_0^2] < \infty$ and $E[X_0] = 0$ and 
\begin{align}
\label{eq5.81}   
    &\int_0^{\exp(-x^q)} Q_{|X(0)|}^2(u)\, du\  <<\ 
(1/x)^{1-q} \cdot g(x)\ \ 
{\rm as}\ x \to \infty, \quad  {\rm and}\\
\label{eq5.82}   
    &\alpha_X(n)\ \leq\ \beta_X(n)\ \leq\ \exp(-n^q) \ \
    {\rm for\ all}\ n \in \N,  
\end{align}
and conclusions (iii) and (iv) in Theorem \ref{thm3.4} hold,
and conclusion (iv) in Theorem \ref{thm4.4} holds as well [with $g_n = g(n)$].  
\end{corollary}

   To derive this corollary, apply Theorem \ref{thm5.5}
with $f(x) := \exp (-x^q)$ for $x \in [1, \infty)$.
Here we note only that
verifying condition (e) in Theorem \ref{thm5.5} is equivalent to verifying that
$\log(x \cdot f(x)) = (\log x) - x^q$ is nonincreasing (as $x$ increases)
for $x \in (1, \infty)$ sufficiently large.
For that it suffices to note that the derivative satisfies the inequality
$(1/x) -qx^{q-1} < 0$ for all $x \in (1, \infty)$ sufficiently large.

\begin{remark}
\label{rem5.9}
Here are a few further comments on Corollaries \ref{cor5.7} and \ref{cor5.8}
and on Theorem \ref{thm5.5}.
\medskip

   (A) Corollary \ref {cor5.7} shows that for power-type mixing rates, 
Theorem \ref{thm5.3} is very sharp, even for the class of 
strictly stationary, countable-state Markov chains that are reversible.   
The combination of properties (\ref {eq5.71}) and (\ref {eq5.72}) 
(with $g(x)$ allowed to blow up arbitrarily slowly) 
is just marginally weaker than the assumption (c) in Theorem \ref {thm5.3}.
But recall Remark \ref{rem4.6}, regarding the refinement in
[\cite{ref-journal-MP2}, Theorem 4].
\smallskip

     In some sense, Corollary \ref{cor5.7}  
complements the construction (involving reversible Markov chains 
and power-type mixing rates) of 
Doukhan et al.\ [\cite {ref-journal-DMR1994}, Theorems 2 and 5 and Corollary 1],
which showed that for power-type mixing rates, Theorem \ref{thm5.2} is in a certain
sense as sharp as possible [as was pointed out in Remark \ref{rem5.4}(B) above]. 
\smallskip

     In essence, Corollary \ref{cor5.7} deals with that part of Theorem \ref{thm5.5} 
where the random variables $X_k$ have a finite absolute moment of some order 
greater than 2.
(See again the two paragraphs right before Notations \ref{nt5.1}.) 
\medskip
 
   (B) Now refer again to the ``basic insight'' of Cuny and Lin in Remark \ref{rem2.9}(A).
Refer to Theorem \ref{thm2.7},
and to the counterexamples (under say ``barely slower than 
exponential'' mixing rates for $\alpha$-mixing or $\beta$-mixing) 
alluded to in Remark \ref{rem2.9}(B)(C).
\smallskip   
    
    Corollary \ref{cor5.8} seems to fill much of the ``chasm''  
between Corollary 5.7 and the ``basic insight'' in Remark \ref{rem2.9}(A).   
In the context of Corollary \ref{cor5.8} (with $g(x) \to \infty$ arbitrarily slowly),
one has for example, 
$\int_0^{\alpha(n)} Q_{|X(0)|}^2(u)\, du$ ``almost'' $ <<\ n^{q - 1}$ as $n \to \infty$.
Compare that to the (non-Markovian) counterexamples from [\cite {ref-journal-Bradley1997}, Theorem 3] alluded to in Remark \ref{rem5.4}(C) (in connection
with Theorem \ref{thm5.3}), where (say)
$\int_0^{\alpha(n)} Q_{|X(0)|}^2(u)\, du <<\ n^{-1}$ as $n \to \infty$.
\smallskip

     It is an open question whether, for mixing rates distinctly faster than
power-type (but, say, slower than exponential), 
Theorem \ref{thm5.5} is essentially sharp.
If that turns out to be the case, then the comparison in the preceding 
paragraph would indicate that reversibility (if it is satisfied) would provide,
for mixing rates that are distinctly faster than power-type 
(but, say, slower than exponential), 
some small but nontrivial extra leverage in the development of 
central limit theory for (functions of) strictly stationary Markov chains.
A natural open question (for mixing rates that are distinctly faster
than power-type but slower than exponential) is whether or not a CLT would hold if
in the context of Theorem \ref{thm5.5}, without changes of any real 
significance in the other conditions, one replaces eq.\ (\ref{eq5.51}) 
by the barely stronger condition
$\int_0^{f(x)} Q_{|X(0)|}^2(u)\, du\  =\ o\bigl(-f'(x)/f(x)\bigl)$ 
[or even just $O(\dots)$?] as $x \to \infty$.
\medskip

     (C) Certain assumptions in Theorem \ref{thm5.5} appear to be
of only ``secondary'' value, and could reasonably be 
weakened or omitted altogether.
However, such an effort to seek ``maximal generality'' in Theorem \ref{thm5.5}
seems pointless unless and until the question in (B) above is 
answered in the affirmative.    
\end{remark}

\begin{remark}
\label{rem5.10}
   The following comments may provide some additional perspective on the
question posed in Remark \ref{rem5.9}(B).
\smallskip

   Herrndorf [\cite {ref-journal-Herrndorf1985}. p.\ 543, Corollary to Theorem 2]
gave an analog of Theorem \ref{thm4.1} 
for strictly stationary (not necessarily Markovian) sequences
$X := (X_k, k \in \Z)$ such that 
$\alpha(n) \to 0$ at least exponentially fast
and for some $\gamma > 1$,
\begin{equation}
\label{eq5.1001}
E[X_0^2(\log^+|X_0|)^\gamma]\ < \infty.
\end{equation}
Here of course $\log^+ a = \log a$ resp.\ 0 if $a > 1$ resp.\ $0 \leq a \leq 1$.
Later, from covariance equalities of Rio [\cite {ref-journal-Rio1993}, pp.\ 592-593]
and the CLT of Doukhan et al.\ \cite {ref-journal-DMR1994},
it was established that Herrndorf's CLT (with exponential $\alpha$-mixing) 
would still hold under the ``moment'' assumption (\ref{eq5.1001}) with
$\gamma = 1$.
From non-Markovian counterexamples (such as in \cite {ref-journal-Bradley1997}) 
it has been known that such a CLT (with exponential $\beta$-mixing) 
would not hold under the assumption that (\ref{eq5.1001}) is satisfied 
for every $\gamma \in (0,1)$.
\smallskip

     Cuny and Lin [\cite {ref-journal-CL2025}, Proposition 3.7]
described an example of a (say real) strictly stationary, countable-state,
non-reversible Markov chain $Y := (Y_k, k \in \Z)$
satisfying exponential $\beta$-mixing, 
and a Borel function $h: \R \to \R$, such that
the random sequence $(X_k, k \in \Z)$ defined by $X_k := h(Y_k)$ for $k \in \Z$
satisfies $E[X_0] = 0$ and $E[X_0^2] > 0$ and  
(\ref{eq5.1001}) for every $\gamma \in (0,1)$ but fails to satisfy the CLT.
Recall from Theorem \ref{thm2.7} that if
if the Markov chain $Y$ here were also reversible, then
(\ref{eq5.1001}) with $\gamma = 0$ [for $X_0 = h(Y_0)$] 
--- that is, $E[X_0^2] < \infty$ --- would suffice for a CLT. 
\smallskip

That comparison suggests that perhaps, for the development of the CLT
for (functions of) strictly stationary Markov chains that satisfy $\alpha$-mixing 
or $\beta$-mixing with a mixing rate that is distinctly faster than power-type 
but slower than exponential (for example, sub-exponential), 
the amount of extra leverage provided by 
reversibility (if satisfied) might be something akin to allowing a ``moment assumption''
(such as in (\ref{eq5.1001})) to be ``weakened by a small power of $\log^+|X_0|$''.
\end{remark}

\section {Preliminaries}
\label{sc6} 

     This section will lay out assorted preliminary information that will be employed in
Sections \ref{sc7}, \ref{sc8}, and \ref{sc9} in the proofs of Theorems 
\ref{thm3.4}, \ref{thm4.4}, and \ref{thm5.5}. 
This section will be divided into three disjoint pieces: 

      Sub-section 6A will present a class of probability functions on $\Z$ that will 
play a key role in connection with the partial limiting distributions in those three
theorems.

     Sub-section 6B will present a class of strictly stationary, 3-state, reversible Markov chains that will serve as ``building blocks'' for the constructions of the Markov chains in
the three theorems.

     Sub-section 6C will present a review of various assorted miscellaneous pieces of 
mathematical information that will be employed in the proofs of the three theorems.

     A large portion of the material here in Section \ref{sc6} will be taken verbatim from
Sections 4-6 of the paper \cite {ref-journal-Bradley2025}.   
\hfil\break

     {\bf Sub-section 6A:\ \ A class of probability functions on $\Z$.}
In the convergence in distribution described in the first paragraph of 
Notations \ref{nt3.3}, a key role will be played by the following 
class of probability functions for discrete random variables.
\medskip

\begin{notations}\label{nt6.1}    
     For each pair of real numbers $a,p \in (0,1)$ (the open unit interval), 
define the function ${\bf g}_{a,p}: \Z \to [0,1]$ as follows:
\begin{equation}
\label{eq6.11}
{\bf g}_{a,p}(0)\ :=\ 1- a, 
\indent {\rm and\ for\ each}\ n \in \N, \indent
{\bf g}_{a,p}(n)\ :=\ {\bf g}_{a,p}(-n)\ :=\ (a/2) p(1-p)^{n-1}\ . 
\end{equation}
For any given pair of numbers $a,p \in (0,1)$, the
numbers ${\bf g}_{a,p}(k)$, $k \in \Z$ are positive and
(by a simple calculation) add up to 1.
Thus ${\bf g}_{a,p}$ is a probability function on the set $\Z$ of all integers.
For typographical convenience, the probability function ${\bf g}_{a,p}$ will 
sometimes be denoted as ${\bf g}[a,p]$; that is,
${\bf g}[a,p](k) = {\bf g}_{a,p}(k)$ for $k \in \Z$. 
\medskip
\end{notations}

\begin{lemma}
\label{lem6.2}
Suppose that for every ordered pair
$(a,p) \in (0,1)^2$, $J(a,p)$ is a positive integer; 
and suppose that 
\begin{equation}
\label{eq6.21}
\lim_{a \to 0+,\thinspace p \to 0+} a \cdot J(a,p)\ =\ 1\ .
\end{equation}
Suppose that for every ordered pair $(a,p) \in (0,1)^2$,
$\zeta^{(a,p)}_1, \zeta^{(a,p)}_2, \dots,
\zeta^{(a,p)}_{J(a,p)}$ are independent, identically distributed
integer-valued random variables, each with the probability function 
${\bf g}_{a,p}$ in eq.\ (\ref{eq6.11}) (with the same parameters $a$ and $p$).
Then 
\begin{equation}
\label{eq6.22}
p \cdot \sum_{k=1}^{J(a,p)} \zeta^{(a,p)}_k\
\to\ \mu_{P1sL}\ \ {\rm in\ distribution\ as}\ \  
a \to 0+,\  p \to 0 +.
\end{equation}
\end{lemma}

    That is [\cite {ref-journal-Bradley2025}, Lemma 4.2]; its proof is given there.
\smallskip 

    The probability functions ${\bf g}_{a,p}$ will appear again in 
Lemma \ref {lem6.6} later on in this section.  
Lemmas \ref {lem6.2} and \ref{lem6.6} will together play a key role in the proofs, 
in Sections \ref{sc7}-\ref{sc9}, of
Theorems \ref{thm3.4}, \ref{thm4.4}, and \ref{thm5.5}.
\hfil\break

     {\bf Sub-section 6B:\ \ The ``building block'' Markov chains.}\ \
This is basically an abbreviated repeat of Section 5 of \cite {ref-journal-Bradley2025},
with the proofs there omitted here. 
In Definition \ref{def6.4} below, we shall define a class of 
fairly simple strictly stationary 3-state, reversible Markov chains that will 
be used in Sections \ref{sc7}, \ref{sc8}, and \ref{sc9} as 
``building blocks'' in the construction of the Markov chains 
for Theorems \ref{thm3.4}, \ref{thm4.4}, and \ref{thm5.5}.
The ``building block'' Markov chains here will be the same as 
in the construction in the paper \cite {ref-journal-Bradley2025}
(and related to, though distinctly different from, the ``building block''
Markov chains in the construction in the paper
\cite {ref-journal-Bradley2024}).   
\medskip

   In Notations \ref{nt6.3}, we shall formulate  
some key properties of a few items that will 
then be, explicitly or implicitly, key components 
in Definition \ref{def6.4}: \zhfb
the (invariant) probability vectors, in eq.\ (\ref{eq6.331});
\zhfb
the one-step transition probability matrices, in
eq.\ (\ref{eq6.341}); and \zhfb
the  {\it joint\/} probability functions of two consecutive
random variables, in (\ref{eq6.321}).     
\smallskip

   In Lemmas \ref{lem6.5} and \ref{lem6.6},
certain key properties (including reversibility) will be listed for 
the Markov chains in Definition \ref{def6.4}.

\begin{notations}
\label{nt6.3}
 In the ``building block'' Markov chains in Definition \ref{def6.4}, the three states
will be $-1$, $0$, and $1$.
Accordingly, the relevant $3 \times 3$ matrices will
be indexed by $\{-1,0,1\} \times \{-1,0,1\}$ (instead of, say, $\{1,2,3\} \times \{1,2,3\}$),
and the $1 \times 3$ row vectors will be indexed by $\{-1,0,1\}$.
\medskip
 
   {\bf Part 1.  The parameters.}\ \ In the notations that follow, the two main 
parameters will be ``small'' positive numbers 
$\zeps$ and $\theta$.  
As an ongoing reminder that they are intended to be
``small'', we shall impose the conditions
\begin{equation}
\label{eq6.311}
0 < \zeps \leq 1/9 \quad {\rm and} \quad
0 < \theta \leq 1/9\ .                    
\end{equation}
As a function of those two parameters 
satisfying (\ref{eq6.311}), we shall define another,
closely related parameter $\theta^*$ as follows:
\begin{equation}
\label{eq6.312}
\theta^*\ =\ \theta^*(\zeps, \theta)\ :=\
{\theta \over {1 - \zeps}}\ .
\end{equation}
Of course under (\ref{eq6.311}) and (\ref{eq6.312}),
one has (with redundancy here), for convenient later 
reference, that 
\begin{equation}
\label{eq6.313}
(1 - \zeps)\thinspace \theta^* = \theta, \quad 
{\rm and} \quad
0 < \theta < \theta^* \leq 1/8 \quad {\rm and} \quad
1\ >\ 1 - \theta\ >\ 1 - \theta^*\ \geq\ 7/8\ >\ 0\ . 
\end{equation}

   {\bf Part 2.  Matrices for joint probabilities.}\ \ The following  $3 \times 3$ 
matrices will be used for {\it joint\/} (not conditional) probabilities
(for consecutive random variables in the Markov chains defined in
Definition \ref{def6.4}).
 For any $\zeps$ and any $\theta$ satisfying (\ref{eq6.311}), 
 define the $3 \times 3$ matrix
$\Lambda^{(\zeps, \theta)} = 
(\lambda^{(\zeps, \theta)}_{i,j},\ i,j\in \{-1,0,1\})$ 
as follows:
\begin{align}
\label{eq6.321} 
\lambda^{(\zeps, \theta)}_{0,0}  
&=\ 1 - \zeps - \theta\zeps; \nonumber\\
\lambda^{(\zeps, \theta)}_{-1,-1}\ 
=\ \lambda^{(\zeps, \theta)}_{1,1}   
&=\ (1 - \theta)\zeps/2 ; \nonumber\\
\lambda^{(\zeps, \theta)}_{-1,1}\ 
=\ \lambda^{(\zeps, \theta)}_{1,-1}\ 
&=\ 0;\ \ {\rm and} \nonumber \\
\lambda^{(\zeps, \theta)}_{0,-1}\ 
=\ \lambda^{(\zeps, \theta)}_{0,1}\       
= \lambda^{(\zeps, \theta)}_{-1,0}\ 
=\ \lambda^{(\zeps, \theta)}_{1,0}\
&=\ \theta\zeps/2.  
\end{align}      
Note that the nine entries in this matrix are all 
nonnegative, and they add up to 1.
Also, note that this matrix is symmetric:
$\lambda^{(\zeps, \theta)}_{i,j} =
\lambda^{(\zeps, \theta)}_{j,i}$
for $i,j \in \{-1,0,1\}$. 
(There is further obvious symmetry as well.) 
 \medskip
 
     {\bf Part 3.  Vectors for marginal distributions.}\ \ 
For each $\zeps \in (0, 1/9]$ [as in eq.\ (\ref{eq6.311})],
define the $1 \times 3$ row vector 
$\pi^{(\zeps)} := 
[\pi^{(\zeps)}_{-1}, \pi^{(\zeps)}_0, \pi^{(\zeps)}_1]$
as follows:
\begin{equation}
\label{eq6.331}
\pi^{(\zeps)}_0\ =\ 1 - \zeps \indent {\rm and} \indent
\pi^{(\zeps)}_{-1}\ =\ \pi^{(\zeps)}_1\ =\ \zeps/2\ .
\end{equation}
By simple arithmetic, assuming (\ref{eq6.311}) as usual,  
if $(X,Y)$ is a random vector
taking its values in the set $\{-1,0,1\}^2$ with joint
probability function $\Lambda^{(\zeps, \theta)}$
in (\ref{eq6.321}), then for each of the random 
variables $X$ and $Y$ separately, the marginal probability function is 
$\pi^{(\zeps)}$.
\medskip  

     {\bf Part 4.  Matrices for one-step transition probabilities.}\ \ 
The matrices in eq.\ (\ref{eq6.341}) below will be employed later on 
as one-step transition probability matrices 
for the ``building block'' Markov chains in Definition \ref{def6.4}.
\medskip 
     
     For any $\zeps$ and $\theta$ satisfying
(\ref{eq6.311}), define the $3 \times 3$ matrix
$\bbP^{(\zeps, \theta)} := (p^{(\zeps, \theta)}_{ij},\ 
i,j \in \{-1,0,1\})$ as follows 
[employing (\ref{eq6.312}) 
in the first two rows of this display]:
\begin{align}
\label{eq6.341}
p^{(\zeps, \theta)}_{0,0}\ &=\ 1 - \theta^*\zeps;
\nonumber \\
p^{(\zeps, \theta)}_{0,-1}\ =\ p^{(\zeps, \theta)}_{0,1}\ 
&=\ \theta^* \zeps/2;
\nonumber \\  
p^{(\zeps, \theta)}_{-1,-1}\ =\ p^{(\zeps, \theta)}_{1,1}\  
&=\ 1 - \theta; \nonumber \\     
p^{(\zeps, \theta)}_{-1,0}\ =\ p^{(\zeps, \theta)}_{1,0}\ 
&=\ \theta; \nonumber \\ 
p^{(\zeps, \theta)}_{-1,1}\ =\ p^{(\zeps, \theta)}_{1,-1}\ 
&=\ 0. 
\end{align}

By (\ref{eq6.321}), (\ref{eq6.331}), (\ref{eq6.341}), and the initial
equality in (\ref{eq6.313}),  
one has [under (\ref{eq6.311})--(\ref{eq6.312})] that
\begin{equation}
\label{eq6.342}
{\rm for\ every}\ (i,j) \in \{-1,0,1\}^2,\ \ \
\pi^{(\zeps)}_i \cdot p^{(\zeps, \theta)}_{i,j}\ 
=\ \lambda^{(\zeps, \theta)}_{i,j}\ ; 
\end{equation}
and also one obtains the matrix product equality
\begin{equation}
\label{eq6.343}
\pi^{(\zeps)} \bbP ^{(\zeps, \theta)} = \pi^{(\zeps)}\ .    
\end{equation}
\end{notations}

\begin{definition}
\label{def6.4}
     Suppose $\zeps$ and $\theta$ satisfy (\ref{eq6.311}).
Refer to (\ref{eq6.312}).
A Markov chain $X := (X_k, k \in \Z)$ 
is said to satisfy ``Condition ${\cH}(\zeps, \theta)$'' if 
$X$ is a strictly stationary Markov chain such that \hfil\break
(i) its state space is $\{-1, 0, 1\}$, \hfil\break
(ii) its (invariant) marginal probability vector is the vector 
$\pi^{(\zeps)}$ in (\ref{eq6.331}), and \hfil\break
(iii) its one-step transition probability matrix is the matrix
$\bbP^{(\zeps, \theta)}$ in (\ref{eq6.341}).     
\smallskip

   That is [under Condition ${\cH}(\zeps, \theta)$],
for $i,j \in \{-1, 0, 1\}$, 
$P(X_0 = i) = \pi_i^{(\zeps)}$ from (\ref{eq6.331}), 
and $P(X_1=j|X_0 = i) = p^{(\zeps, \theta)}_{i,j}$ from
(\ref{eq6.341}).  
The requirements (ii) and (iii) in this definition are
compatible by (\ref{eq6.343}).
In fact Condition ${\cH}(\zeps, \theta)$ uniquely 
determines the distribution of the entire Markov chain $X$
on $(\R^\Z, \cR^\Z)$.  
\end{definition}

   As a quick summary, [under Condition 
${\cH}(\zeps, \theta)$ for $\zeps$ and $\theta$ 
satisfying (\ref{eq6.311})], the Markov chain $X$ is 
strictly stationary, finite-state (in fact 3-state),
and [as an elementary by-product of (\ref{eq6.341})]
irreducible and aperiodic.
As noted in Background \ref{bckg3.1}, the Markov chain $X$ 
therefore satisfies $\beta$-mixing.
It is therefore also ergodic
(see e.g.\ [\cite {ref-journal-Bradley2007}, Vol.\ 1, 
Theorem 7.7 and Chart 5.22]). 
\medskip

\begin{lemma}
\label{lem6.5} 
Suppose $\zeps$ and $\theta$ satisfy (\ref{eq6.311}). 
Suppose $X := (X_k, k \in \Z)$ is a (strictly stationary) Markov chain 
that satisfies Condition $\cH (\zeps, \theta)$.
Then the state space of this Markov chain $X$ is
$\{-1,0,1\}$; and the following statements hold: \zhfb
\noindent (1) The (joint) probability function of the random 
vector $(X_0, X_1)$ is the matrix $\Lambda^{(\zeps, \theta)}$ 
from (\ref{eq6.321}).
\zhfb
\noindent (2) This Markov chain $X$ is reversible. \zhfb
\noindent (3a) One has that $P(X_0 = 0) = 1 - \zeps$ and
$P(X_0 = -1) = P(X_0 = 1) = \zeps/2$. \zhfb
\noindent (3b) For each $s \in \{-1,0,1\}$, 
$P(\{X_0 = 0\} \cap \{X_1 = s\}) > 0$; in fact
$P(\{X_0 = 0\} \cap \{X_1 = 0\}) \geq 1 - 2\zeps$. \zhfb 
\noindent (4) One has that $E(X_0) = 0$ and 
$E|X_0| =  E(X_0^2) = \var(X_0) = \zeps > 0$. \zhfb
\noindent (5) The sequence $Y := (-X_k, k \in \Z)$
is a strictly stationary Markov chain that satisfies 
Condition ${\cal H}(\zeps, \theta)$. \zhfb 
\noindent (6) For each $n \in \N$, 
$\beta_X(n) \leq 6 \zeps (1 - \theta)^n 
\leq 6 \zeps \exp(-\theta n)$. \zhfb 
\noindent (7) For each $n \in \N$,
$P(X_0 = X_n = 1) = 
(\zeps/2)(1 - \theta)^n + P(\{X_0 = 1\} \cap \{X_n = -1\})$.
\zhfb
\noindent (8a) For each $n \in \N$,
$E(X_0X_n) = \cov(X_0, X_n) = \zeps (1 - \theta)^n > 0$. \zhfb
\noindent (8b) One has that 
$E(X_0^2) + 2 \sum_{n=1}^\infty E(X_0X_n) 
= \zeps [ (2/\theta) -1]$.
\zhfb 
\noindent (9) Letting 
$u_n := (1/n) E[(X_1 + X_2 + \dots + X_n)^2]$ for
each $n \in \N$, one has that
\begin{equation}
\label{eq6.51}
u_1 \leq u_2 \leq u_3 \leq \dots \quad {\rm and} \quad
\lim_{n \to \infty} u_n = \zeps [(2 / \theta) - 1].  
\end{equation}
\noindent (10) As $n \to \infty$, the normalized partial sum
$n^{-1/2} (X_1 + X_2 + \dots + X_n)$
converges in distribution to the normal law with 
mean 0 and variance $\zeps [(2 / \theta) - 1]$. 
\end{lemma}

     That lemma is [\cite {ref-journal-Bradley2025}, Lemma 5.3]; 
its proof is given there.
Note that the Markov chain in Lemma \ref{lem6.5} is ``positively correlated''
in the sense of property (8a).
That property yields other nice properties, such as the inequalities in the
``first half'' of eq.\ (\ref{eq6.51}), which in turn will help establish certain key properties 
of the variances of partial sums in the counterexamples constructed in this paper. 
(In fact it will be easy to see that the counterexamples constructed in this paper, 
as well as the one constructed for [\cite {ref-journal-Bradley2025}, Theorem 3.3], are all
``positively correlated'' in this sense.)

\begin{lemma}
\label{lem6.6} 
Suppose $\zeps$ and $\theta$ satisfy (\ref{eq6.311}).
Suppose $X := (X_k, k \in \Z)$ is a strictly stationary
Markov chain that satisfies Condition $\cH(\zeps, \theta)$
(see Definition \ref{def6.4}).
Refer to (\ref{eq6.312}) and (\ref{eq6.313}).
Then there exists a sequence $(W_1, W_2, W_3, \dots)$
of independent, identically distributed, integer-valued
random variables with the following two properties: \zhfb
(1) The probability function of each $W_m$ ($m \in \N$)
is ${\bf g}[\theta^*\zeps, \theta]$ (recall Notations \ref{nt6.1},
including its last sentence).
\zhfb      
(2) Letting $I$ be the positive integer defined by 
\begin{equation}
\label{eq6.61}
I\ :=\ \biggl[ {1 \over {\theta^*\zeps}} \biggl]. 
\end{equation}
[where in this particular equation, the tall brackets indicate
that the right hand side is the greatest integer that is
$\leq 1/(\theta^*\zeps)$], one has that
\begin{equation}
\label{eq6.62}
P \biggl( \sum_{k=1}^I X_k\ \neq\ \sum_{j=1}^I W_j \biggl)\
\leq\ 3\zeps.
\end{equation}
\end{lemma}

   That lemma is [\cite{ref-journal-Bradley2025}, Lemma 5.4];
its proof is given there.
\hfil\break

     {\bf Sub-section 6C:\ \ Other Preliminaries.}\ \ 
In this sub-section, we shall collect together several miscellaneous technical 
background facts in analysis and probability theory that will be employed in 
Sections \ref{sc7}-\ref{sc9} in the proofs of Theorems 
\ref{thm3.4}, \ref{thm4.4}, and \ref{thm5.5}.
\vskip 0.1 in    

  We start here with a standard technical lemma. 
A proof for it can be found e.g.\ in
[\cite {ref-journal-Bradley2007}, Vol.\ 1, Theorem 6.2(III)].
It will be employed in the arguments in 
Sections \ref{sc8} and \ref{sc9}.
\smallskip 

\begin{lemma}
\label{lem6.7}
Suppose $\cA_1, \cA_2, \cA_3, \dots$ and
$\cB_1, \cB_2, \cB_3, \dots$ are 
$\sigma$-fields [on the given probability space
$(\Omega, \cF, P)$], and the
$\sigma$-fields $\cA_j \vee \cB_j$, $j \in \N$
are independent.
Then
\begin{equation}
\nonumber
\beta\biggl(\, \bigvee_{j=1}^\infty \cA_j,\ 
\bigvee_{j=1}^\infty \cB_j \biggl)\
\leq\ \sum_{j=1}^\infty \beta(\cA_j, \cB_j). 
\end{equation}  
\end{lemma}

   The following lemma is of course well known. 
It is stated here for convenient reference for arguments
in Sections \ref{sc7}-\ref{sc9}.
\medskip

\begin{lemma}
\label{lem6.8}
     Suppose $\zeta_1, \zeta_2, \zeta_3, \dots$ is a sequence
of independent square-integrable random variables such that
(i) $E[\zeta_j] = 0$ for each $j \in \N$, and
(ii) $\sum_{j=1}^\infty E[\zeta_j^2] < \infty$.
     Then the sum $\sum_{j=1}^\infty \zeta_j$ converges almost surely and 
in $\cL^2$ to a square-integrable random variable $Y$ such that $EY = 0$ 
and $E(Y^2) = \sum_{j=1}^\infty E(\zeta_j^2)$, 
\end{lemma}

   The next lemma is just an elementary technical fact that 
will be employed in order to show that the constructions in 
Sections \ref{sc7}-{\ref{sc9} are in fact Markov chains ---
not just functions of Markov chains.

\begin{lemma}
\label{lem6.9}
    Suppose $y_1, y_2, y_3, \dots$ and $z_1, z_2, z_3, \dots$
are each a sequence of elements of the set $\{-1,0,1\}$.
Suppose that (i) there exists $m \in \N$ such that
$y_k = z_k = 0$ for all $k > m$, and
(ii) there exists $k \in \N$ such that $y_k \neq z_k$.

   Suppose $r_1, r_2, r_3, \dots$ is a sequence of positive numbers 
such that one of the following two conditions holds:
\zhfb
(a) $r_{k+1}/r_k \geq 3$ for every $k \in \N$; or
\zhfb
(b) $r_{k+1}/r_k \leq 1/3$ for every $k \in \N$.

   Then each of the sums $\sum_{k=1}^\infty r_ky_k$ and 
$\sum_{k=1}^\infty r_kz_k$ has at most finitely many non-zero summands 
(and therefore trivially converges absolutely); and they satisfy 
$\sum_{k=1}^\infty r_ky_k \neq \sum_{k=1}^\infty r_kz_k$.
\end{lemma} 

     This lemma is elementary and well known.
Its statement and proof under hypothesis (a) were spelled out in
[\cite {ref-journal-Bradley2024}, Lemma 6.3].
Here we shall spell out the (similar) argument under hypothesis (b). 
\smallskip

     Suppose the numbers $y_k$, $z_k$, and $r_k$ for $k \in \N$ are as in
the statement of this lemma, with hypothesis (b) satisfied.
Let $n$ denote the least positive integer such that $y_n \neq z_n$. 
It suffices to show that 
$\sum_{k=n}^\infty r_ky_k \neq \sum_{k=n}^\infty r_k z_k$.
\smallskip

      Now $|r_ny_n - r_nz_n| = r_n |y_n - z_n| \geq r_n \cdot 1 = r_n$.
Also, 
\begin{equation}
\nonumber
\biggl|\sum_{k=n+1}^\infty  r_ky_k - \sum_{k=n+1}^\infty r_kz_k\biggl|\ 
\leq\ \sum_{k=n+1}^\infty r_k |y_k - z_k|\ 
<\  \sum_{k=n+1}^\infty 2r_k\
\leq\ \sum_{k=n+1}^\infty 2 \cdot 3^{-(k-n)} r_n \ 
=\ r_n,          
\end{equation}
where the strict inequality $\ldots < \dots$ in that display comes from the fact
that $y_k - z_k = 0$ for some $k \geq n+1$ 
(in fact all sufficiently large $k$). 
These two facts together force the inequality at the end of the preceding 
paragraph to hold.
That completes the argument.

\begin{notations}
\label{nt6.10}
   In the next lemma, and in arguments in 
Sections \ref{sc8} and \ref{sc9},    
the term ``affine function $L: \R \to \R$'' 
of course just means a function (``line'') given by 
$L(x) = q + rx$ for $x \in \R$, 
where $q$ and $r$ are each a real number.
For such an affine function $L$, the notation
``$[{\rm slope\ of}\ L]$'' will simply mean its (constant) 
derivative (the parameter $r$ in the preceding sentence).

   (Those notations and the next two lemmas will later help 
facilitate applications of Lemma \ref{lem6.5}(6) together with
Lemma \ref{lem6.7}.)
\end{notations}  

   The next lemma involves some elementary
analysis for convex functions.
It will be applied in both of Sections \ref{sc8} and \ref{sc9}.
Its formulation here is unnecessarily restrictive (and also has 
redundancy); but it will be fine for those applications.
 
\begin{lemma}
\label{lem6.11}
Suppose $\phi: [1,\infty) \to (-\infty, 0]$ is a function
that satisfies the following conditions: \zhfb
(a) $\phi$ is continuous, strictly decreasing, and convex 
on $[1, \infty)$; \zhfb
(b) $\phi$ is continuously differentiable on the open 
half line $(1, \infty)$; and \zhfb 
(c) one has that
\begin{equation}
\label{eq6.111}
\lim_{x \to \infty} \phi(x)\ =\ -\infty 
\indent {\rm and} \indent
\lim_{x \to \infty} \frac {\phi(x)} {x}\ =\ 0.
\end{equation}

   For each $y \in (1, \infty)$ (the open half line), let 
$L^{(y)}: \R \to \R$ denote the affine function whose graph is tangent to that 
of $\phi$ at the point $(y, \phi(y))$ --- that is, the affine function that satisfies
\begin{equation}
\label{eq6.112}
L^{(y)}(y)\ =\ \phi(y) \indent {\rm and} \indent 
L^{(y)}(x)\ \leq\ \phi(x)\ \ {\rm for\ all}\ x \in [1, \infty). 
\end{equation}
 
   Then the following statements (I), (II), and (III) hold:
\medskip

   (I) The following statements hold: \zhfb
(i) For each $x \in (1, \infty)$,
$\phi(x) < \phi(1) \leq 0$ and $\phi'(x) < 0$.  \zhfb  
(ii) The function $x \mapsto \phi'(x)$ for $x \in (1, \infty)$, 
is nondecreasing. \zhfb  
(iii) One has that $\lim_{x \to \infty} \phi'(x) = 0$. \zhfb
(iv) For each $y \in (1, \infty)$,  
$[{\rm slope\ of}\ L^{(y)}] = \phi'(y) < 0$. \zhfb 
(v) For each $y \in (1, \infty)$, 
$L^{(y)}(0) = \phi(y) + y \cdot (-\phi'(y)) > \phi(y)$.
\medskip

   (II) For any negative number $D$ and any positive number 
$s$, there exists a number    
$T = T(\phi, D, s) > 1$ with the following property:

   For each $y \geq T$, the affine function 
$L^{(y)}$ in (the entire sentence of) eq.\ (\ref{eq6.112})
satisfies 
\begin{equation}
\label{eq6.113}
-s\ \leq\ [{\rm slope\ of}\ L^{(y)}]\ =\ \phi'(y)\ <\ 0 
\indent {\rm and} \indent
L^{(y)}(0)\ \leq\ D. 
\end{equation}
  
   (III) Suppose also that $\psi: [1, \infty) \to [0, \infty)$ is a  
function such that $\lim_{x \to \infty} \psi(x) = \infty$.
 
Then for any negative number $D$ and any positive numbers
$B$ and $s$, there there exists a real number 
$T^* := T^*(\phi, \psi, B, D, s) > 1$ with the following 
property:
\smallskip

For each $y \geq T^*$, the affine function 
$L^{(y)}$ in (the entire sentence of) eq.\ (\ref{eq6.112})
satisfies (\ref{eq6.113}) and also 
[strengthening part of (\ref{eq6.113})] satisfies  
\begin{equation}
\label{eq6.114}  
L^{(y)}(0) + B\ \leq\ D \indent {\rm and} \indent
L^{(y)}(x) + B\ \leq\ \phi(x) + \psi(x)\ \ 
{\rm for\ all}\ x \geq 1.    
\end{equation}
\end{lemma}

   Apart from a couple of trivial extra details listed in (I) here, the (elementary)
proofs of (I) and (II) were spelled out in 
[\cite {ref-journal-Bradley2025}, proof of Lemma 6.5].
[Note that --- because of the ``extra'' assumption (b) here in Lemma \ref{lem6.11}
(not present in [\cite {ref-journal-Bradley2025}, Lemma 6.5]) ---
the (``bad'') set $\Gamma$ in the context there 
in \cite {ref-journal-Bradley2025} is, in our context here,  
empty and therefore omitted.]
\smallskip

   The proof of Statement (III) will be given here:
\hfil\break 

{\bf Proof of Lemma \ref {lem6.11}(III).}\ \
Suppose $D$, $B$, $s$, and the function $\psi$
are as in statement (III). 
To help keep track of the various threads in the argument, we
shall label each paragraph here in order: (P1), (P2), etc.
\smallskip
 
(P1) By the hypothesis of statement (III), one has that 
$D < 0 \leq \psi(x)$ for all $x \in [1, \infty)$.

(P2) Referring again to the assumptions on $\psi$,
let $q \in (1, \infty)$ be such that
$\psi(x) \geq B$ for all $x \geq q$.

(P3) By the assumptions in the lemma here
[including the assumption that
$\phi(x) \leq 0$ for all $x \in [1, \infty)$],
one has that $D - B + \phi(q) < 0$.
Referring to the notation $T(\dots)$ in statement (II),
define the positive number
$T^*(\phi, \psi, B, D, s)$ by 
\begin{equation}
\label{eq6.11P1}
T^*(\phi, \psi, B, D, s)\ :=\ 
T \Bigl(\phi,\ D - B + \phi(q),\ s \Bigl).
\end{equation}

(P4) Now suppose that 
$y \geq T^*(\phi, \psi, B, D, s)$.
Refer to the affine function $L^{(y)}$ in 
(the entire sentence of) eq.\ (\ref{eq6.112}).
To complete the proof of statement (III), our task is
to show that (\ref{eq6.113}) and (\ref{eq6.114}) hold.

(P5) By (P4) and (\ref{eq6.11P1}),\ \   
$y\ \geq\ T(\phi,\ D - B + \phi(q),\ s)$. 
Hence from statement (II), one has that
$-s < [{\rm slope\ of}\ L^{(y)}] = \phi'(y) < 0$ and that
$L^{(y)}(0) \leq D - B + \phi(q) < D - B < D$
[since $B > 0$ by hypothesis; and $q>1$ by (P2) 
and hence $\phi(q) < \phi(1) \leq 0$ 
by the assumptions on the function $\phi$].
Thus (\ref{eq6.113}) and in fact the ``first half'' of 
(\ref{eq6.114}) hold.
To complete the proof of statement (III), all that remains
is to verify the ``second half'' of (\ref{eq6.114}) --- that is,
to show that for every $x \geq 1$, 
\begin{equation}
\label{eq6.11P2}  
L^{(y)}(x) + B\ \leq\ \phi(x)\, +\, \psi(x)\ .  
\end{equation}

(P6) For any $x \in [q, \infty)$ one has by (\ref{eq6.112}) and (P2)
that $L^{(y)}(x) + B \leq \phi(x) + B \leq \phi(x) + \psi(x)$. 
Thus (\ref{eq6.11P2}) holds for all $x \geq q$.

(P7) For reasons that will be spelled out in order below, one has that 
for any $x \in [1,q]$, the following hold:
\begin{equation}
\label{eq6.11P3}
L^{(y)}(x)\ <\ L^{(y)}(0)\ \leq\ D-B+\phi(q)\ \leq\ D-B+\phi(x)\
<\ \psi(x) - B + \phi(x).
\end{equation}
Here, by the inequality $[{\rm slope\ of}\ L^{(y)}] < 0$
in (P5), the first inequality in (\ref{eq6.11P3}) holds. 
Then because of the inequality 
$y\, \geq\, T(\phi,\, D - B + \phi(q),\, s)$ in (P5),
the second inequality in (\ref{eq6.11P3}) comes simply from 
an application of statement (II).   
Next, since (by hypothesis) the function $\phi$ is strictly decreasing 
on $[1,\infty)$, the third inequality in (\ref{eq6.11P3}) holds 
(under the stipulation here that $x \in [1,q]$).
The final inequality in (\ref{eq6.11P3}) holds by (P1). 
Thus (\ref{eq6.11P2}) holds for all $x \in [1,q]$.

(P8) By the final sentence of each of paragraphs (P6) and (P7),
one has that (\ref{eq6.11P2}) holds for all $x \geq 1$.
By the final sentence of (P5), the proof of
Lemma \ref{lem6.11}(III) is complete.
\medskip

   The next lemma involves quantiles [as in eq.\ (\ref{eq5.11}) in 
Section \ref{sc5}], and will be applied in Section {\sc9}.

\begin{lemma}
\label{lem6.12}
   Suppose $h_1, h_2, h_3, \dots$ is a sequence of positive
numbers such that $h_j \leq (1/2) h_{j+1}$ for each
$j \in \N$.

   Suppose $A_1, A_2, A_3, \dots$ is a sequence of events
[on the probability space $(\Omega, \cF, P)$] such that
$P(A_j) > 0$ for each $j \in \N$, and $P(A_1) < 1/2$ and
$P(A_j) \leq (1/2)P(A_{j-1})$ for each $j \geq 2$.
For each $j \in \N$, define the positive number 
$a_j := 2 \cdot P(A_j)$.
Then
\begin{align}
\label{eq6.121}
&0\ <\ \dots\ <\ a_3\ <\ a_2\ <\ a_1\ <\ 1\ \quad
{\rm and} \quad \lim_{j \to \infty} a_j = 0;\  \indent {\rm and} \\
\label{eq6.122}
&P(A_j\ {\rm occurs\ for\ infinitely\ many}\ j \in \N)\ =\ 0. 
\end{align}

   Referring to (\ref{eq6.122}), let $Y$ denote the 
nonnegative random variable defined ($a.s.$) by
\begin{equation}
\label{eq6.123}
Y\ :=\ \sum_{j=1}^\infty h_j I(A_j)
\end{equation}
[where $I(\dots)$ refers to the indicator function].

   Then referring to (\ref{eq5.11}) and (\ref{eq6.121}), one has that
\begin{equation}
\label{eq6.124}
Q_Y(u)\ =\ 0\ \ {\rm for\ every}\ u \in [a_1,\, 1)\ ;
\end{equation}
and for each $j \geq 2$, one has that
\begin{equation}
\label{eq6.125}
Q_Y(u)\ \leq\ 2h_{j-1}\ \ {\rm for\ every}\ 
u \in [a_j,\, a_{j-1})\ .
\end{equation}
\end{lemma}

   {\bf Proof.}\ \ 
Eq.\ (\ref{eq6.121}) holds as a trivial consequence of the assumptions on the 
events $A_j$ and the definition of the numbers $a_j$.
\smallskip

   For any given $j \in \N$, one has by induction that
$P(A_i) \leq (1/2)^{i-j} P(A_j)$ for every $i \geq j$.
Hence for every $j \in \N$,
\begin{equation}
\label{eq6.126}
P\biggl(\ \bigcup_{i=j}^\infty A_i\biggl)\
\leq\ \sum_{i=j}^\infty P(A_i)\
\leq\ \sum_{i=j}^\infty (1/2)^{i-j}P(A_j)\
=\ 2P(A_j)\ =\ a_j\ .
\end{equation}

   By (\ref{eq6.126}) (starting with its second term) 
with (say) $j=1$, together with
(\ref{eq6.121}) and the Borel-Cantelli Lemma, one has that (\ref{eq6.122}) holds.
\medskip

What remains is to prove eq.\ (\ref{eq6.124}) and 
(for $j \geq 2$) eq.\ (\ref{eq6.125}).
But first we shall carry out one further preliminary observation:
For each $j \geq 2$, one has by induction that
$h_i \leq (1/2)^{j-i} h_j$ for every $i \in \{1,2,\dots, j\}$.
Hence for each $j \in \N$ (trivially also including $j=1$), 
\begin{equation}
\label{eq6.127}
\sum_{i=1}^j h_i\ \leq\ \sum_{i=1}^j (1/2)^{j-i}h_j\
=\ \sum_{\ell = 0}^{j-1} (1/2)^\ell h_j\
<\ \sum_{\ell = 0}^\infty (1/2)^\ell h_j\ =\ 2h_j\ .
\end{equation}

   {\bf Proof of (\ref{eq6.124}).}\ \  
By (\ref{eq6.123}), one has the equality of events 
$\{Y>0\} = \bigcup_{i=1}^\infty A_i$.
For each $u \in [a_1, 1)$ [see (\ref{eq6.121})]
one has by (\ref{eq6.126}) that 
$P(Y > 0) = P(\bigcup_{i=1}^\infty A_i) \leq a_1 \leq u$. 
Hence (see (\ref{eq5.11})) for each $u \in [a_1,1)$,
$Q_Y(u) = \inf\{t \geq 0: P(Y > t) \leq u\} = 0$
(since by the preceding sentence the number $t=0$ itself is in that set, 
and is the least possible element therein).
Thus (\ref{eq6.124}) holds. 
\medskip

     {\bf Proof of (\ref{eq6.125}).}\ \ Suppose $j \geq 2$.  
For each $\omega \in \Omega$ such that 
$\omega \in \{\bigcup_{i=j}^\infty A_i\}^c$
(where the superscript $c$ denotes complement in $\Omega$),
one has by (\ref{eq6.123}) and (\ref{eq6.127}) that
$Y(\omega) = \sum_{i=1}^{j-1} h_iI(A_i)(\omega)
\leq \sum_{i=1}^{j-1} h_i \leq 2h_{j-1}$.  
Hence $\{Y > 2h_{j-1}\} \subset \bigcup_{i=j}^\infty A_i$.
Hence for each $u \in [a_j, a_{j-1})$ [see (\ref{eq6.121})],
one has by (\ref{eq6.126}) that
$P(Y > 2h_{j-1}) \leq P(\bigcup_{i=j}^\infty A_i) 
\leq a_j \leq u$. 
Hence for each $u \in [a_j, a_{j-1})$, 
$Q_Y(u) = \inf \{t \geq 0: P(Y > t) \leq u\} \leq 2h_{j-1}$
(since the number $t = 2h_{j-1}$ itself is in that set). 
Thus (\ref{eq6.125}) holds for the (arbitrary) given $j \geq 2$.
That completes the proof of Lemma \ref{lem6.12}.

\begin{lemma}
\label{lem6.13}
     Suppose $f: (1,\infty) \to (0,\infty)$ is a (positive) continuous nonincreasing
function [on the open half line $(1,\infty)$]  
such that $\lim_{x \to \infty} f(x) = 0$.

     Suppose $g:[1,\infty) \to [1, \infty)$ is a function such that
$\lim_{x \to \infty} g(x) = \infty$.

      Then there exists a continuous, nondecreasing function 
$h: [1, \infty) \to [1, \infty)$ with the following properties: \zhfb
(i) $1 \leq h(x) \leq g(x)$ for all $x \in [1, \infty)$; \zhfb
(ii) $\lim_{x \to \infty} h(x) = \infty$; and \zhfb
(iii) the function $x \mapsto h(x) \cdot f(x)$ for $x \in (1,\infty)$
is nonincreasing.
\end{lemma}

   In this lemma, the function $g$ is not assumed to be continuous or monotonic.
\medskip

   \indent {\bf Proof.}\ \ 
For the purpose of keeping track of the threads in the argument, the paragraphs
here will be labeled in order as (P1), (P2), etc.
\smallskip

   (P1) Recursively define an alternating sequence
$t_1 < u_1 < t_2 < u_2 < t_3 < u_3 < \dots$ of positive numbers 
with the following properties:
(i) $t_1 := 1$; 
(ii) for each $j \in \N$, $u_j > t_j + 1$ and $g(x) \geq j+1$ for all $x \geq u_j$; 
(iii) for each $j \in \N$, $t_{j+1} > u_j$ and $f(t_{j+1}) = (j/(j+1)) \cdot f(u_j)$.
(For that last equality, we are using the Intermediate Value Theorem in calculus,
together with the assumptions on the function $f$ here.)
\smallskip

   (P2) By (P1)(i), (P1)(ii) (its first inequality), and (P1)(iii), one has that
\begin{equation}
\label{eq6.131}
1 = t_1 < u_1 < t_2 < u_2 < t_3 < u_3 < \dots \uparrow \infty.
\end{equation}
With that in mind, define the function $h: [1,\infty) \to (0,1]$ as follows: \zhfb
(i) for every $j \in \N$, define $h(x) := j$ for all $x \in [t_j, u_j]$; 
in particular, $h(1) = h(t_1) := 1$; \zhfb
(ii) for every $j \in \N$, define
$h(x) := j \cdot f(u_j)/f(x)$ for all $x \in (u_j, t_{j+1})$. \zhfb
By (\ref{eq6.131}), that completes the definition of the function $h$ on $[1, \infty)$.
Now let us verify the properties of $h$ asserted in the lemma.
First review the assumptions on the function $f$.  
\smallskip

   (P3) {\it The function $h$ is continuous on $[1, \infty)$}.
Refer to (\ref{eq6.131}).
The function $h$ is constant on each of the closed intervals
$[1, u_1],\  [t_2, u_2],\  [t_3, u_3],\  [t_4, u_4],\  \dots$, and continuous 
on each of the open intervals $(u_1,t_2),\  (u_2, t_3),\  (u_3,t_4),\  \dots$. 
For each $j \in \N$, $h(x) \to j = h(u_j)$ as $x \to u_j+$ by (P2)(ii) and (P2)(i); and 
$h(x) \to j \cdot f(u_j)/f(t_{j+1}) = j+1 =  h(t_{j+1})$ as $x \to t_{j+1}-$
by (P2)(ii), (P1)(iii), and (P2)(i).
Hence $h$ is continuous on $[1, \infty)$.  
\smallskip

   (P4) {\it The function $h$ is nondecreasing on $[1,\infty)$}.
The function $h$ is constant on each of the closed intervals
$[1, u_1],\  [t_2, u_2],\  [t_3, u_3],\  \dots$, and nondecreasing 
on each of the open intervals $(u_1,t_2),\  (u_2, t_3),\  (u_3,t_4),\  \dots$
[by (P2)(ii) and the assumption that the (positive) function $f$ is nonincreasing].
Hence by (P3) (continuity), $h$ is nondecreasing on $[1, \infty)$. 
\smallskip

   (P5) {\it The product $h(x) \cdot f(x)$ is nonincreasing on the 
open half line $(1,\infty)$}.
Again refer to (\ref{eq6.131}). 
Separately on the half-open interval $(1, u_1]$ and on each of the close intervals 
$[t_2, u_2],\  [t_3, u_3],\  [t_4, u_4],\  \dots$,
the product $h(x) \cdot f(x)$ is nonincreasing since $h$ is constant and $f$ is nonincreasing.
For each $j \in \N$, $h(x) \cdot f(x) = j \cdot f(u_j)$, a constant, for $x \in (u_j, t_{j+1})$
by (P2)(ii).
Also, the product $f(x) \cdot h(x)$ is continuous on $(1,\infty)$ [by (P3) and the
continuity of $f$].
It follows that that product is nonincreasing on $(1, \infty)$.
\smallskip
 
   (P6) {\it One has that $\lim_{x \to \infty} h(x) = \infty$}.  This is a trivial
consequence of (P4) and (P2)(i).
\smallskip

   (P7) {\it One has that $1 \leq h(x) \leq g(x)$ for all $x \in [1, \infty)$}. 
Again refer to (\ref{eq6.131}).     
One has that $h(x) \geq h(1) = 1$ 
for all $x \in [1, \infty)$ by (P4) and (P2)(i).
For each $x \in [1, u_1]$, $h(x) = 1 \leq g(x)$ by (P2)(i) and the 
assumptions on $g$.
For each $j \geq 2$ and each $x \in [t_j, u_j]$, 
$h(x) = j \leq g(x)$ by (P2)(i) and the fact from (P1)(ii)  
that $g(x) \geq j$ for all $x \geq u_{j-1}$.
For each $j \in \N$ and each $x \in (u_j, t_{j+1})$,
$h(x) \leq h(t_{j+1}) = j+1 \leq g(x)$ by (P4), (P2)(i), and (P1)(ii).
Thus $1 \leq h(x) \leq g(x)$ for all $x \in [1, \infty)$.
\smallskip

   (P8) All properties asserted for the function $h$ in the lemma have been
verified in (P3)--(P7) above.
That completes the proof. 
    
\begin{lemma}
\label{lem6.14}
Suppose $(E_1, E_2, E_3, \dots)$ is a sequence of independent events
on a probability space $(\Omega, \cF, P)$, such that 
$\sum_{j=1}^\infty P(E_j^c) < \infty$ (where the superscript 
$c$ denotes complement).

   Suppose that on the same probability space, $(F_1, F_2, F_3, \dots)$ is a
sequence of independent events such that
(i) $P(F_j) > 0$ for each $j\in \N$, and 
(ii) $F_j = E_j$ (an equality of events) for all except at most finitely
many indices $j \in \N$.
Then $P(\bigcap_{j=1}^\infty F_j) = \prod_{j=1}^\infty P(F_j) > 0$.  
\end{lemma}

   Of course the final equality holds by independence.
The subsequent final inequality ($\dots > 0$) is simply an application of the 
basic fact in analysis that if $(s_1, s_2, s_3, \dots)$ is a sequence of 
numbers in $[0,1)$ such that $\sum_{j=1}^\infty s_i < \infty$,
then $\prod_{j=1}^\infty (1 - s_j) > 0$.
[Simply take $s_j = P(F_j^c)$.]

\section{Proof of Theorem \ref{thm3.4}}
\label{sc7}

\begin{remark}
\label{rem7.1}
    In the rest of this paper, in order to simplify certain mathematical
arguments, we shall sometimes build unnecessary redundancies into the
choices of parameters, and also we shall sometimes settle on adjectives of
``less than full strength'' (for example, stating that a sequence of real numbers
is ``monotonically decreasing'' when in fact it is strictly decreasing) when
``full strength'' is unnecessary.       
\end{remark}

The proof of Theorem \ref{thm3.4} will be carried out in a series of  ``steps''.
\medskip

   {\bf Step 1.  The parameters for the construction.}\ \
This ``step'' will be divided into several ``sub-steps''.
\vskip 0.1 in

     {\bf Sub-step 1A.}\ \ For each positive integer $j$, define the 
positive number $h_j$ by
\begin{equation}
\label{eq7.1A1}
h_j\ :=\ 1/3^j.
\end{equation}

     Recursively choose a sequence $(\zeps_1, \zeps_2, \zeps_3, \dots)$ 
of positive numbers in the interval $(0, 1/9]$ such that
for every $j \geq 2$,  
\begin{equation}
\label{eq7.1A2}
\zeps_j\ \leq\ \min\{ 9^{-j}, \zeps_{j-1} \} \indent {\rm and} \indent
\frac {h_j^2 \zeps_j} {h_\ell^2 \zeps_\ell^2}\ \leq\ 2^{-j}\
{\rm for\ every}\ \ell\ {\rm such\ that}\ 1 \leq \ell \leq j-1.   
\end{equation} 
(Note the absence of an exponent 2 for $\zeps_j$
in the numerator of the fraction.)\ \ 
\medskip     

{\bf Sub-step 1B.}\ \ 
Refer to the sequence of numbers $(q_1, q_2, q_3, \dots)$ in
the statement of Theorem \ref{thm3.4}. 
Recursively choose positive integers 
$M_j$ for $j \in \N$ such that
\begin{equation}
\label{eq7.1B1}
1 \leq M_1 < M_2 < M_3 < \dots;  \quad 
{\rm and\ for\ each}\ j \in \N,\ {\rm one\ has\ that}\
q_n \leq h_j^2 \zeps_j/2\ \ {\rm for\ all}\ n \geq M_j.
\end{equation}

{\bf Sub-step 1C.}\ \ 
Recursively choose a sequence $(\theta_1, \theta_2, \theta_3,\dots)$ of 
positive numbers in the interval $(0,1/9]$ such that 
$\theta_1 M_2 \leq 1/2$ and for each $j \geq 2$,  
\begin{equation}
\label{eq7.1C1}
\theta_j\ \leq\  \min\{9^{-j}, \theta_{j-1}\} \indent {\rm and} \indent  
\theta_j M_{j+1} \leq 1/2 \indent {\rm and} \indent
h_j^2 \zeps_j/\theta_j\ \geq\ 
j \cdot \biggl(1 + \sum_{u=1}^{j-1} 
(h_u^2 \zeps_u/\theta_u)\biggl). \quad
\end{equation}

   Recall from above that for each $j \in \N$, 
$\zeps_j \in (0, 1/9]$ and $\theta_j \in (0, 1/9]$, which
imply $0 < \theta_j/(1 - \zeps_j) \leq 1/8$.
Accordingly, for each $j \in \N$, define the positive number $\theta^*_j$
[\`a la (\ref{eq6.312})] and the positive integer $I_j$ by
\begin{equation}
\label{eq7.1C2}
\theta^*_j\ :=\ \frac {\theta_j} {1 - \zeps_j}
\indent {\rm and} \indent 
 I_j\ :=\ \biggl[\frac {1} {\theta_j^* \zeps_j} \biggl]
\end{equation}    
[where in the latter equality the brackets indicate the
greatest integer that is $\leq 1/(\theta_j^* \zeps_j)$].
\medskip

    {\bf Sub-step 1D.\ \ A summary of some key properties of the parameters:}  \zhfb
(a) For each $j \in \N$, $\zeps_j \in (0, 1/9]$, $\theta_j \in (0, 1/9]$,
$\theta^*_j \in (0, 1/8]$,  $I_j \geq 72$, and $h_j > 0$. \zhfb 
(b) As $j \to \infty$, $\zeps_j \to 0$ monotonically
and $\theta_j \to 0$ monotonically. \zhfb
(c)  As $j \to \infty$, $\theta^*_j \to 0$ monotonically
and $I_j \to \infty$ monotonically. \zhfb 
(d)  For each $j \in \N$, one has that $\theta_j/\theta^*_j < 1$ and
$\theta_j^*\zeps_j I_j \leq 1$.  
\zhfb
(e) As $j \to \infty$, one has that $\theta_j/\theta^*_j \to 1$ and 
$\theta_j^*\zeps_j I_j \to 1$ and $\theta_j/h_j \to 0$. \zhfb
(f)  One has that 
$\sum_{j=1}^\infty h_j^2 \zeps_j  < \sum_{j=1}^\infty \zeps_j \leq 1/8$. 
\zhfb 
(g) For each $j \geq 2$, one has that
$h_j^2 \zeps_j / \theta_j \geq j \cdot [1+ \sum_{u=1}^{j-1} (h_u^2 \zeps_u / \theta_u)]
> j$. 
\medskip

   Of the five assertions in (a), the first two are built into
the choices of the numbers $\zeps_j$ and $\theta_j$ in (\ref{eq7.1A2}) and
(\ref{eq7.1C1}), and the other three then follow from (\ref{eq7.1C2})
or hold by (\ref {eq7.1A1}).
Property (b) holds by (\ref{eq7.1A2}) and (\ref{eq7.1C1}),
and property (c) then follows from (\ref{eq7.1C2}).
Property (d) holds by (\ref{eq7.1C2}) and property (a). 
In property (e), the first two limits hold by (\ref{eq7.1C2}) and  
properties (a)(b)(c), and the third limit holds by (\ref{eq7.1A1}) and 
the ``first third'' of (\ref{eq7.1C1}).
Property (f) holds by (\ref{eq7.1A1}) and (\ref{eq7.1A2}) and property (a).  
Property (g) is simply a repeat of the ``final third'' of (\ref{eq7.1C1})
(and then a trivial application of property (a)).       
\medskip
 
   {\bf Step 2.  The construction.}\ \
This step will repeat practically verbatim Step 3 in Section 7 of
\cite {ref-journal-Bradley2025} (in the proof of the main result of that paper).
For the reader's convenience, we shall start by essentially repeating here
(with superficial changes in equation numbers, etc.) 
the first three ``sub-steps'' of Step 3 there.   
\medskip

   {\bf Sub-step 2A.}\ \ For each positive integer $j$, 
referring to the parameters
$\zeps_j$ and $\theta_j$ defined in Sub-steps 1A and 1C
[see Sub-step 1D(a)], let
$X^{(j)} := (X^{(j)}_k, k \in \Z)$ be a strictly stationary
Markov chain that satisfies Condition ${\cal H}(\zeps_j, \theta_j)$
of Definition \ref{def6.4}.
Let this construction be carried out in such a way that these
(strictly stationary) Markov chains 
$X^{(1)},\ X^{(2)},\ X^{(3)},\ \dots$
are independent of each other.
\medskip

  {\bf Sub-step 2B.}\ \ By Definition \ref{def6.4} and
Lemma \ref{lem6.5}, for each $j \in \N$ and each $k \in \Z$,
the random variable $X^{(j)}_k$ takes its values in the set
$\{-1,0,1\}$, with probabilities
$P(X^{(j)}_k = 0) = 1 - \zeps_j$ and
$P(X^{(j)}_k = -1) = P(X^{(j)}_k = 1) = \zeps_j/2$, and 
hence $P(X^{(j)}_k \neq 0) = \zeps_j$.
Hence by Sub-step 1D(f), one has that
\begin{equation}
\label{eq7.2B1}
{\rm for\ each}\ k \in \Z, \quad \sum_{j=1}^\infty P(X^{(j)}_k \neq 0)\ 
\leq\ 1/8\ <\ \infty.
\end{equation}
Hence for each $k \in \Z$, by the Borel-Cantelli Lemma, one has that
\begin{equation}
\label{eq7.2B2}
P\Bigl(X^{(j)}_k \neq 0\ {\rm for\ infinitely\ many}\
j \in \N \Bigl)\ =\ 0.
\end{equation}

For each $k \in \Z$, let $G_k$ denote the event inside the parentheses
in (\ref{eq7.2B2}).
Then by (\ref{eq7.2B2}) itself, for each $k \in \Z$, $P(G_k) = 0$.
Define the event $G := \bigcup_{k \in \Z}G_k$.
Since the set $\Z$ is countable, one has that $P(G) = 0$.
For each $\omega \in G$, each $j \in \N$, and each $k \in \Z$, redefine
$X^{(j)}_k(\omega) := 0$.
Then without loss of generality, i.e.\ without having changed the joint distribution of the
random variables $(X^{(j)}_k,\ j \in \N,\ k \in \Z)$ (and hence without having changed 
any of the ``independence'' or ``Markov'' properties), we now have the following:
\smallskip

For every $\omega \in \Omega$ and every $k \in \Z$, 
$X^{(j)}_k(\omega) \neq 0$ for at most finitely many $j \in \N$ (possibly no $j \in \N$ at all).
   
\medskip          

 {\bf Sub-step 2C.}\ \ 
Referring to (\ref{eq7.1A1}) and the last sentence of Sub-step 2B, 
define the sequence $X := (X_k, k \in \Z)$ of random variables as follows:   
For each $k \in \Z$ and each $\omega \in \Omega$, 
\begin{equation}
\label{eq7.2C1}
X_k(\omega)\ :=\ \sum_{j=1}^\infty h_j X^{(j)}_k(\omega)\ . 
\end{equation}  
For each $k \in \Z$ and each $\omega \in \Omega$, 
by the last sentence of Sub-step 2B, 
the sum in (\ref{eq7.2C1}) has at most finitely
many non-zero terms and is therefore well defined.
\medskip

     {\bf Sub-step 2D.}\ \ 
Next we shall repeat essentially verbatim Sub-steps 3D, 3E, and 3F
in Section 7 of \cite {ref-journal-Bradley2025}.
The arguments there refer to certain equation numbers, etc.\ from {\it earlier\/}
in that paper.
The following chart gives the match-ups between 
(i) the relevant equation numbers, etc.\ from {\it earlier\/} in that paper, and
(ii) the corresponding equation numbers, etc.\ from {\it earlier\/} in this paper here.  
\begin{tabbing}
{\bf In Sub-steps 3D/3E/3F in Section 7 of \cite {ref-journal-Bradley2025}}  \hskip 1 in                    
\= {\bf In the context here in this paper} \\ 
Sub-steps 3A, 3B, and 3C    \> Sub-steps 2A, 2B, and 2C \\
(7.23) and Lemma 6.3  \>  (\ref{eq7.1A1}) and Lemma \ref{lem6.9}\\  
(7.27) \> (\ref{eq7.2C1})\\
Lemma 6.6 \> Lemma \ref{lem6.14} \\ 
(7.25) \> (\ref{eq7.2B1}) \\
Lemma 5.3(2)(3b) \>  Lemma \ref{lem6.5}(2)(3b) \\
Notations 2.8(B) \> Notations \ref{nt2.6}(B) \\
Sub-step 2D(f) \> Sub-step 1D(f) 
\end{tabbing}

     Then applying the first paragraph of Sub-step 3G of Section 7 
of \cite {ref-journal-Bradley2025}, we obtain our analog of the second paragraph
of Sub-step 3G there:
\medskip

   {\it The random sequence $X$ defined in (\ref {eq7.2C1})
is a strictly stationary, countable-state, irreducible, aperiodic,
reversible Markov chain [with $P(X_0 = s) > 0$ for every state $s$]\/.} 
\medskip

    Referring again to (\ref{eq7.2C1}), let us summarize here a key portion of that material
carried over (explicitly or implicitly) 
from [\cite {ref-journal-Bradley2025}, Section 7, Sub-steps 3D-3G].
Refer to (\ref{eq7.2C1}), Lemma \ref{lem6.9}, and Sub-steps 2A and 2B. 
\smallskip
 
First, for any $k \in \Z$, one has that $\sigma(X_k) = \bigvee_{j=1}^\infty \sigma(X^{(j)}_k)$.
(For a given $k \in \Z$, one has that $X_k$ uniquely determines $X^{(j)}_k$,\ $j \in \N$,
 and vice versa.) 
\smallskip

   Second, the state space of the Markov chain $X$ is the (countably infinite) set of all real numbers of the form $s = \sum_{j=1}^\infty h_j z_j$ for which $(z_1, z_2, z_3, \dots)$
is a sequence of elements of $\{-1, 0, 1\}$ such that $z_j \neq 0$ 
for at most finitely many $j \in \N$ (possibly none).
For any such state $s = \sum_{j=1}^\infty h_j z_j$, one has that
its (invariant) marginal probability satisfies $P(X_0 = s) > 0$ --- in fact
\begin{equation}
\label{eq7.2D1}
P (X_0 = s)\ 
=\ P \biggl(\, \bigcap_{j=1}^\infty \{X^{(j)}_0 = z_j\} \biggl)\ 
=\ \prod_{j=1}^\infty P(X^{(j)}_0 = z_j)\ >\ 0. 
\end{equation}

\medskip

   {\bf Step 3.\ \ Proof of Statement (i)} [in Theorem
\ref {thm3.4}]. 
Statement (i) holds (with $C= 1/2$) by (\ref {eq7.2C1}), (\ref{eq7.1A1}), 
and the first sentence of Sub-step 2B.
\medskip

   {\bf Step 4.}\ \ Refer to Sub-step 2A and the first sentence of Sub-step 2B.
Of course for any given $j \in \N$ and any given 
$k \in \Z$, one has that $E[X^{(j)}_0] = 0$ and $\var [X^{(j)}_0] = \zeps_j$
and hence also $E[h_jX^{(j)}_k] = 0$ and 
$E[(h_jX^{(j)}_k)^2] = \var[h_jX^{(j)}_k] = h_j^2 \zeps_j$.
\smallskip
 
From that and stationarity and Minkowski's Inequality,
for any given $j \in \N$ and any given $n \in \N$,
one trivially has $E(\sum_{k=1}^n h_jX^{(j)}_k) = 0$ and   
$\var (\sum_{k=1}^n h_jX^{(j)}_k) \leq n^2 h_j^2 \zeps_j$.    
Also, trivially by (\ref{eq7.1A1}) and Sub-step 1D(a), 
$\sum_{j=1}^\infty  h_j^2 \zeps_j 
< \sum_{j=1}^\infty  h_j^2 < \infty$. 
Hence by Lemma \ref{lem6.8} and the second sentence of Sub-step 2A,
for any (finite or countably infinite) 
nonempty set $S \subset \N$ and any $n \in \N$,
the random variable 
$\sum_{j \in S} \sum_{k=1}^n h_j X^{(j)}_k$
is square-integrable [in fact it is bounded by (\ref{eq7.1A1})] 
and has mean 0 and satisfies 
\begin{equation}
\label{eq7.41}
E \biggl[\biggl(\,\sum_{j \in S} \sum_{k=1}^n h_j X^{(j)}_k\biggl)^2 \, \biggl]\
=\
\sum_{j \in S} \var \biggl(\, \sum_{k=1}^n h_j X^{(j)}_k\biggl)\ 
=\ \sum_{j \in S} \biggl[h_j^2\, 
\var \biggl(\, \sum_{k=1}^n X^{(j)}_k\biggl)\biggl].  
\end{equation}
As a special case, with $S = \N$ and $n = 1$, one obtains (say) eq.\ (\ref{eq3.21}) —
$E[X_0^2] < \infty$ (which of course holds trivially by Statement (i) in the theorem
— recall Step 3 above) and $E[X_0] = 0$.
\medskip

     {\bf Step 5.  Proof of Statement (ii)} [in Theorem \ref {thm3.4}].
In this argument, we shall let the integer $N$ in 
Statement (ii) be the integer $M_1$ from (\ref {eq7.1B1}).
Now let $m$ be an arbitrary fixed positive integer such that
$m \geq M_1$.  
To complete the proof of Statement (ii), it suffices to prove for this $m$ that
$\var(\, \sum_{k=1}^m X_k) \geq q_m \cdot m^2$.

\medskip

   Referring to (\ref{eq7.1B1}), let $J \in \N$ be such that
$M_J \leq m \leq M_{J+1}$.
\smallskip

   By Sub-step 2A and Lemma \ref{lem6.5}(4)(8a), and then (\ref{eq7.1C1})
(its entire sentence) and (\ref{eq7.1B1}), one has that
\begin{align}
\label{eq7.51}
\var&\biggl(\, \sum_{k=1}^m X^{(J)}_k\biggl)\
=\ E\biggl[\biggl(\, \sum_{k=1}^m 
X^{(J)}_k\biggl)^2\, \biggl]\
=\ \sum_{k=1}^m \sum_{\ell = 1}^m 
E\Bigl(X^{(J)}_k X^{(J)}_\ell\Bigl)\  
=\  \sum_{k=1}^m \sum_{\ell = 1}^m 
\zeps_J (1 - \theta_J)^{| k - \ell |}\ \nonumber \\
&\geq\ m^2 \zeps_J(1 - \theta_J)^m\
\geq\ m^2 \zeps_J(1 - m\theta_J)\
\geq\ m^2 \zeps_J(1 - M_{J+1}\theta_J)\
\geq\ m^2 \zeps_J \cdot 1/2\
\geq\ m^2 \cdot q_m/h_J^2 .       
\end{align}
   
   Now by (\ref {eq7.2C1}) and the first equality in (\ref{eq7.41}) with $S = \N$, 
and then (\ref{eq7.51}), one has that
\begin{equation}
\nonumber
\var\biggl(\, \sum_{k=1}^m X_k\biggl)\ 
=\ \sum_{j \in \N}
\var\biggl(\, \sum_{k=1}^m h_j X^{(j)}_k\biggl)\
\geq \ \var\biggl(\, \sum_{k=1}^m h_J X^{(J)}_k\biggl)\
\geq\ m^2 \cdot q_m,
\end{equation}   
giving the desired inequality.
That completes the proof of Statement (ii).
\medskip

   {\bf Step 6.  Proof of Statement (iii)} [in Theorem 3.4].\ \
Eq.\ (\ref{eq3.21}) in Statement (ii) 
was confirmed right after (\ref{eq7.41}). 
As noted right after eq.\ (\ref{eq3.23}) itself,
eq.\ (\ref {eq3.22}) then follows from (\ref {eq3.23}).
To complete the proof of Statement (iii), it now suffices 
to verify (\ref {eq3.23}).
\smallskip

     To verify (\ref{eq3.23}), we shall repeat essentially verbatim the argument
(for the same equation) in Step 6 (starting with its second paragraph)
in Section 7 of \cite {ref-journal-Bradley2025}.
The argument there refers to certain equation numbers, etc.\ from {\it earlier\/}
in that paper.
The following chart gives the match-ups between 
(i) the relevant equation numbers, etc.\ from {\it earlier\/} in that paper, and
(ii) the corresponding equation numbers, etc.\ from \ {\it earlier\/} in this paper here.  
\begin{tabbing}
{\bf In Step 6 in Section 7 of \cite {ref-journal-Bradley2025}}  
\hskip 1 in    \= {\bf In the context here in this paper} \\
Remark 3.1 \> Remark \ref {rem3.2} \\ 
(3.3)  \>  (\ref {eq3.23}) \\
Sub-step 2D(g)  \>  Sub-step 1D(g) \\  
Lemma 5.3(10) \> Lemma \ref{lem6.5}(10) \\
Sub-steps 3A and 3B \> Sub-steps 2A and 2B \\ 
(7.27) \> (\ref{eq7.2C1})
\end{tabbing}

   That completes the proof of (\ref {eq3.23}),
and of statement (iii) in Theorem \ref{thm4.4}.
 \medskip

    {\bf Step 7.  Proof of Statement (iv).}\ \
The argument will be divided into four ``Sub-steps'', as was the corresponding 
argument in Step 7 of Section 7 of the paper \cite {ref-journal-Bradley2025}.
Sub-steps 7A, 7B, and 7D here will repeat essentially verbatim the corresponding
Sub-steps there in that paper.
Only Sub-step 7C here will differ from the corresponding
Sub-step there.
For the reader's convenience, we shall spell out the whole argument here.
 
\medskip

    {\bf Sub-step 7A.}\ \ 
Refer to Notations \ref{nt6.1}, including its last sentence.
Refer to Sub-step 1D(a) and the first sentence of Sub-step 2A.   
For each $j \in \N$, applying Lemma \ref{lem6.6}, let
$W^{(j)} := (W^{(j)}_1, W^{(j)}_2, W^{(j)}_3, \dots)$
be a sequence of independent, identically distributed,
discrete random variables, each having the probability function 
${\bf g}[\theta^*_j \zeps_j, \theta_j]$, 
such that
\begin{equation}
\label{eq7.7A1}
P \biggl(\, \sum_{k=1}^{I(j)} X^{(j)}_k\ \neq\ 
\sum_{u=1}^{I(j)} W^{(j)}_u \biggl)\
\leq\ 3\zeps_j,
\end{equation}
where here and below, the positive integer $I_j$ from
(\ref{eq7.1C2}) is also written as $I(j)$ for typographical convenience.

\smallskip

Recall from Sub-step 1D(e) that
$\theta^*_j \zeps_j \cdot I_j \to 1$ as $j \to \infty$.
By Lemma \ref{lem6.2}, 
$\theta_j\sum_{u=1}^{I(j)} W^{(j)}_u \to \mu_{P1sL}$ in 
distribution as $j \to \infty$.            
As a consequence of (\ref{eq7.7A1}) and Sub-step 1D(b), 
$(\theta_j \sum_{k=1}^{I(j)} X^{(j)}_k) 
- (\theta_j \sum_{u=1}^{I(j)} W^{(j)}_u)\, \to 0$
in probability as $j \to \infty$.
Hence by Slutsky's Theorem,
$\theta_j \sum_{k=1}^{I(j)} X^{(j)}_k\ \to\ \mu_{P1sL}$
in distribution as $j \to \infty$.
That is,
\begin{equation}
\label{eq7.7A2}
\frac {\theta_j} {h_j} \sum_{k=1}^{I(j)} h_j X^{(j)}_k\ \to\ \mu_{P1sL}\ 
{\rm in\ distribution\ as}\ j \to \infty. 
\end{equation}

     {\bf Sub-step 7B.}\ \ 
 By Sub-step 2A and Lemma \ref{lem6.5}(9),
for each $u \in \N$ and each $n \in \N$, 
$\var(\sum_{k=1}^n X^{(u)}_k) < 
n \cdot 2\zeps_u/\theta_u$. 
Hence for each $j \geq 2$, by Lemma \ref{lem6.8},
both inequalities in Sub-step 1D(d), and then Sub-step 1D(g),
\begin{align}
\nonumber
E& \biggl[ \biggl(\frac {\theta_j} {h_j} 
\sum_{u=1}^{j-1} \sum_{k=1}^{I(j)} h_u X^{(u)}_k
\biggl)^2 \, \biggl]\
=\ \frac {\theta_j^2} {h_j^2} 
\sum_{u=1}^{j-1} \biggl[\, h_u^2\, 
\var \biggl(\, \sum_{k=1}^{I(j)} X^{(u)}_k \biggl) \biggl]\
\leq\ \frac {\theta_j^2} {h_j^2} 
\sum_{u=1}^{j-1} \Bigl(h_u^2 \cdot 
I_j \cdot 2\zeps_u/\theta_u \Bigl) \nonumber \\
\nonumber
&=\ \frac {\theta_j^2} {h_j^2} I_j 
\sum_{u=1}^{j-1} \Bigl(h_u^2 \cdot 
2\zeps_u/\theta_u \Bigl)\  
\leq\ \frac {\theta_j^2} {h_j^2} 
\frac {1} {\theta^*_j \zeps_j} 
\sum_{u=1}^{j-1} \Bigl(h_u^2 \cdot 
2\zeps_u/\theta_u \Bigl)\
\leq\ \frac {\theta_j} {h_j^2 \zeps_j}  
\sum_{u=1}^{j-1} \Bigl(h_u^2 \cdot 
2\zeps_u/\theta_u \Bigl)\
\leq\ \frac {2} {j}\ .  
\end{align}
As a consequence,
\begin{equation}
\label{eq7.7B1}
\frac {\theta_j} {h_j} 
\sum_{u=1}^{j-1} \sum_{k=1}^{I(j)} h_u X^{(u)}_k\
\to\ 0\ \, {\rm in\ probability\ as}\ j \to \infty.
\end{equation}
 
   {\bf Sub-step 7C.}\ \ 
For each $u \in \N$ and each $n \in \N$, by strict stationarity, Lemma \ref{lem6.5}(4),
and a trivial application of Minkowski's inequality, one has that
$E [(\sum_{k=1}^n X^{(u)}_k)^2] \leq n^2 \zeps_u$. 
Hence for each $j \in \N$, by Lemma \ref{lem6.8},
both inequalities in Sub-step 1D(d), and
then (\ref{eq7.1A2}),  
\begin{align}
\nonumber
E&\biggl[ \biggl(\frac {\theta_j} {h_j} 
\sum_{u=j+1}^{\infty} \sum_{k=1}^{I(j)} h_u X^{(u)}_k
\biggl)^2 \, \biggl]\
=\ \frac {\theta_j^2} {h_j^2} 
\sum_{u=j+1}^{\infty} \biggl[\, h_u^2\, 
\var \biggl(\, \sum_{k=1}^{I(j)} X^{(u)}_k \biggl) \biggl]\
\leq\ \frac {\theta_j^2} {h_j^2} 
\sum_{u=j+1}^{\infty} \bigl( h_u^2\, 
I_j^2 \zeps_u \bigl)\\
\nonumber
&=\ \frac {\theta_j^2} {h_j^2} I_j^2 
\sum_{u=j+1}^{\infty} \bigl(\, h_u^2 
\zeps_u \bigl)\  
\leq\ \frac {\theta_j^2} {h_j^2} 
\frac {1} {(\theta^*_j \zeps_j)^2} 
\sum_{u=j+1}^{\infty} \bigl(\, h_u^2 
\zeps_u \bigl)\
\leq\ \frac {1} {h_j^2 \zeps_j^2} 
\sum_{u=j+1}^{\infty} \bigl(\, h_u^2 
\zeps_u \bigl)\
\leq\ \sum_{u = j+1}^\infty 2^{-u}\ 
=\ 2^{-j}.    
\end{align}
As a consequence,
\begin{equation}
\label{eq7.7C1}
\frac {\theta_j} {h_j} 
\sum_{u=j+1}^{\infty} \sum_{k=1}^{I(j)} h_u X^{(u)}_k\
\to\ 0\ \, {\rm in\ probability\ as}\ j \to \infty.
\end{equation}

   {\bf Sub-step 7D.}\ \    
By (\ref{eq7.7A2}), (\ref{eq7.7B1}), (\ref{eq7.7C1}), and 
Slutsky's Theorem, and then (\ref{eq7.2C1}),  
$(\theta_j/h_j) \sum_{k=1}^{I(j)} X_k = 
(\theta_j/h_j) \sum_{u=1}^\infty \sum_{k=1}^{I(j)} h_u X_k^{(u)}$
converges in distribution to $\mu_{P1sL}$ as $j \to \infty$.
Also, from Sub-step 1D(e), $\theta_j/h_j \to 0$ as $j \to \infty$.
Finally [recall Sub-step 1D(c)], one can pass to a suitable subsequence of the $I_j$'s if necessary
in order to have the final sequence of such integers be strictly increasing.
Thus property (iv) in Theorem 3.4 holds.
That completes the proof of Theorem 3.4.

\section{Proof of Theorem \ref{thm4.4}}
\label{sc8}

The proof will be somewhat similar to that of
Theorem \ref{thm3.4}, with some stretches of that 
earlier argument taken verbatim here.   
It will be carried out in a  series of ``steps'', starting with ``Step 0''.
\medskip

   {\bf Step 0.  Preliminary material.}\ \ 
Define the function $\phi: [1, \infty) \to (-\infty, 0]$
by $\phi(x) := \log(1/x) = -\log x$ for $x \in [1, \infty)$.
Then this function $\phi$ satisfies the hypothesis of 
(that is, all of the first paragraph of) Lemma \ref{lem6.11}. 
Of course its derivative on the open half line $(1,\infty)$ is given by
\begin{equation}
\label{eq8.01}
\phi'(x)\ =\ -1/x \quad {\rm for\ all}\ x \in (1,\infty)
\end{equation}

   Refer to the sequence of positive numbers
$(g_1, g_2, g_3, \dots)$ in the statement of 
Theorem \ref{thm4.4}.
Without loss of generality, assume that $g_n \geq 1$ for all $n \in \N$.
Let $\psi: [1, \infty) \to [0,\infty)$ be the function such that
$\psi(n) = \log g_n$ for each $n \in \N$
and the function $\psi$ is affine (``linear'') on each of the
closed intervals $[n, n+1]$, $n \in \N$.
Then (see the hypothesis  of Theorem {\ref {thm4.4}) 
this function $\psi$ satisfies the hypothesis of 
Lemma \ref{lem6.11}(III), that is, $\psi$ is a nonnegative function and
$\lim_{x \to \infty} \psi(x) = \infty$.
\medskip

   {\bf Step 1.  The parameters for the construction.}\ \
In this step, we shall recursively define for each 
positive integer $j$ ten parameters
(one of those ``parameters'' being an affine function). 
They are listed here, along with their most basic restrictions: 
\medskip

    {\bf Key List:} \zhfb
(i) a positive number $B_j$; \zhfb
(ii) a negative number $D_j$; \zhfb
(iii) a positive number ${\bf c}_j$; \zhfb
(iv) a number $t_j \in (1, \infty)$; \zhfb
(v) an affine function ${\bf L}_j: \R \to \R$; \zhfb
(vi) a number $\zeps_j \in (0, 1/9]$; \zhfb
(vii) a number $\theta_j \in (0, 1/9]$; \zhfb
(viii) a number $\theta^*_j \in (0, 1/8]$; \zhfb  
(ix) a positive integer $I_j$; and \zhfb
(x) a positive number $h_j$.
\medskip

   Such parameters $\zeps_j$, $\theta_j$, $I_j$, and $h_j$ 
in (vi)-(vii) and (ix)-(x) will also be defined for $j=0$.
Those four parameters for $j=0$ will provide a convenient ``springboard'' 
to start the recursion (but they will not otherwise play a significant
role later on in the construction itself for Theorem \ref{thm4.4}). 
\medskip
   
    The definition of the positive numbers $h_j$, $j \geq 0$ (item (x) in the 
Key List above) is straightforward and need not be part of the recursion.
It is as follows:
\begin{equation}
\label{eq8.101}
{\rm For\ each}\ j \geq 0, \quad h_j\ :=\ 1/3^j.
\end{equation}
The basic requirement in item (x) in the Key List above is obviously
satisfied.
Trivially $1 = h_0 > h_1 > h_2 > h_3 > \dots \downarrow 0$. 
\medskip

   {\bf Sub-step 1A.\ \  The initial step.}\ \ Now we start the recursion, starting 
with $j=0$, for just a few parameters.  
Define the numbers $\zeps_0$, $\theta_0$, and $I_0$ by
\begin{equation}
\label{eq8.1A1}
\zeps_0\ =\ \theta_0\ := 1/9 \indent {\rm and} \indent I_0\ :=\ 1.
\end{equation}
The basic requirements in items (vi), (vii), and (ix) in the Key List
above are obviously satisfied.
Note (for later reference) the trivial fact that by (\ref{eq8.101}) and (\ref{eq8.1A1}), 
\begin{equation}
\label{eq8.1A2}
\frac {h_0^2 \zeps_0} {\theta_0}\ =\ 1.
\end{equation}

   {\bf Sub-step 1B.\ \  The recursion step.}\ \ 
Now suppose $j \geq 1$ is an integer, and that
for each integer $u$ such that $0 \leq u \leq j-1$,
the numbers $\zeps_u, \theta_u \in (0, 1/9]$ 
and the positive integer $I_u$ have already been defined.
\smallskip

     Referring to (\ref{eq8.101}) and ({\ref{eq8.1A2}), define the parameters
$B_j$ and $D_j$ by
\begin {equation}
\label{eq8.1B1}
B_j\ :=\ (j+2)\ +\ \log\biggl(j \cdot h_j^{-2} \cdot 
\sum_{u=0}^{j-1} \frac {h_u^2 \zeps_u} {\theta_u} \biggl)
\indent {\rm and} \indent 
D_j\ :=\ \log \biggl( \min \biggl\{ \frac {9^{-j}} {I_{j-1}}\, ,\, \frac {\zeps_{j-1}} {2^j} \biggl\} \biggl).    
\end{equation}
Of course by (\ref{eq8.101}) and (\ref{eq8.1A2}), one trivially has 
that $B_j > 0$ (in fact both terms in the definition of $B_j$ are positive) and $D_j < 0$.
\smallskip

     Refer to the assumptions on the sequence $(g_n,\ n \in \N)$ in the statement of
Theorem \ref{thm4.4} [and to the second sentence after (\ref{eq8.01})].
Define the real number ${\bf c}_j$ by
\begin{equation}
\label{eq8.1B1a}
  {\bf c}_j\ :=\ \min_{n \in {\bf U}(j)} \bigl( 2^{-j} \cdot (g_n /n) \cdot \exp(-B_j)\bigl)
\end{equation}
where ${\bf U}(j) := \{1\} \cup \{u \in \N: g_u < 2^{j+1} \exp(B_j) \}$.
That set ${\bf U}(j)$ is nonempty (as ensured by the ``$\{1\} \cup \dots$''), and it is finite; 
the number ${\bf c}_j$ in (\ref{eq8.1B1a}) is well defined, and it is positive.
\smallskip  
 
Refer again to the sentence after (\ref{eq8.1B1}).
Accordingly, in the terminology of Lemma \ref{lem6.11}(III),
using the convex function $\phi(\,.\,)$ and the function $\psi(\,.\,)$ in Step 0, 
define the number $t_j > 1$ by
\begin{equation}
\label{eq8.1B2}    
t_j\ :=\ T^*(\phi, \psi, B_j, D_j, \min\{9^{-j},\, \theta_{j-1}/2,\, {\bf c}_j\}). 
\end{equation}
From (\ref{eq8.01}) and Lemma \ref{lem6.11}(III), one has that
\begin{equation}
\label{eq8.1B3}
-\min\{9^{-j},\, \theta_{j-1}/2,\, {\bf c}_j\}\ \leq\ \phi'(t_j)\ =\ -1/t_j\ <\ 0.   
\end{equation}
In superscripts, the number $t_j$ will be written as $t(j)$ for 
typographical convenience.
In the notations of the second paragraph of Lemma \ref{lem6.11},
for the convex function $\phi$ in Step 0 here above, 
define the affine function ${\bf L}_j: \R \to \R$ by ${\bf L}_j := L^{(t(j))}$, that is,
${\bf L}_j(x) := L^{(t(j))}(x)$ for all $x \in \R$.
Referring to (\ref{eq6.112}) and 
Step 0 [including (\ref{eq8.01})], 
one has that (with a lot of redundancy)
\begin{align}
\label{eq8.1B4a}
&{\bf L}_j(t_j)\, =\, \phi(t_j)\, =\, -\log t_j \indent {\rm and} \indent
{\bf L}_j(x)\, \leq\, \phi(x)\, =\, -\log x\ \ 
{\rm for\ all}\ x \in [1, \infty); \indent {\rm and} \\
\label{eq8.1B4b}
& [ {\rm slope\ of}\ {\bf L}_j]\, =\, \phi'(t_j)\, =\, -1/t_j \indent {\rm and} \indent
{\bf L}_j(x)\, =\, {\bf L}_j(0)\, +\, (-1/t_j)x\ \ {\rm for\ all}\ x \in \R; \indent {\rm and} \\
\label{eq8.1B4c}
&{\bf L}_j(0)\ =\ 1 + {\bf L}_j(t_j)\ =\ 1 + \phi(t_j)\ =\ 1 - \log t_j
\end{align}
[where (\ref{eq8.1B4c}) simply comes from the ``second half'' of (\ref{eq8.1B4b})
and the ``first half'' of (\ref{eq8.1B4a})].
\smallskip 

   Define the positive numbers $\zeps_j$ and $\theta_j$ by
\begin{equation}
\label{eq8.1B5}
\zeps_j\ :=\ \exp\Bigl({\bf L}_j(0) + B_j - (j +2) \Bigl) \indent 
{\rm and} \indent \theta_j\ :=\ 1/t_j.
\end{equation}
By (\ref{eq8.1B2}), the above definition ${\bf L}_j:= L^{(t(j))}$, 
and Lemma \ref{lem6.11}(III), 
${\bf L}_j(0) + B_j \leq D_j$; and hence by (\ref{eq8.1B5}) and  
(\ref{eq8.1B1}) (its ``second half''),
$0 < \zeps_j \leq \exp(D_j) \leq 9^{-j}/I_{j-1} \leq 9^{-j}$.
Also, trivially $0  < \theta_j \leq 9^{-j}$ by 
(\ref{eq8.1B5}) and (\ref{eq8.1B3}); and hence also
$0 < \theta_j/(1 - \zeps_j) \leq 1/8$.
Similarly $0  < \theta_j \leq {\bf c}_j$ by (\ref{eq8.1B5}) and (\ref{eq8.1B3}). 
\smallskip

     Accordingly, define the positive number $\theta^*_j$
[\`a la (\ref{eq6.312})] and the positive integer $I_j$ by
\begin{equation}
\label{eq8.1B6}
\theta^*_j\ :=\ \frac {\theta_j} {1 - \zeps_j}
\indent {\rm and} \indent 
 I_j\ :=\ \biggl[\frac {1} {\theta_j^* \zeps_j} \biggl]
\end{equation}    
[where in the latter equality the brackets indicate the
greatest integer that is $\leq 1/(\theta_j^* \zeps_j)$].
By (\ref{eq8.1B6}) and the two sentences right after 
(\ref{eq8.1B5}), for the given $j \geq 1$, one has that
$0 < \theta^*_j \leq 1/8$, hence $\theta^*_j \zeps_j \leq 1/72$,
and hence $I_j \geq 72$.
\smallskip

   As has been noted at various places above, the basic requirements 
in the Key List of the parameters are all satisfied for the given $j \geq 1$.  
That completes the recursion step.
\smallskip

   That completes the recursive definition of the parameters.
\medskip 

   {\bf Sub-step 1C.  Miscellaneous information.}\ \
By {(\ref{eq8.1B5}) (and the sentence after it) and (\ref{eq8.1B1}),    
$0 < \zeps_j \leq \exp(D_j) \leq \zeps_{j-1}/2^j$ for each $j \geq 1$.  
Also, by (\ref{eq8.1B5}) and (\ref{eq8.1B3}), 
$0 < \theta_j \leq \theta_{j-1}/2$ for each $j \geq 1$.
\smallskip

    Also, by (\ref{eq8.1B5}), (\ref{eq8.1B4c}), (\ref{eq8.1B1}), 
(\ref{eq8.1B5}) again, and (\ref{eq8.1A2}),
for each $j \geq 2$ (say), 
\begin{align}
\nonumber
\zeps_j\ &=\ \exp\bigl( {\bf L}_j(0) \bigl)\, \cdot\, \exp \bigl(B_j - (j+2)\bigl)\
>\ \exp(-\log t_j) \cdot  \exp \bigl(B_j - (j+2)\bigl)\\
\label{eq8.1C1}
&=\ \exp\bigl(\log(1/t_j)\bigl)\, \cdot\, 
\biggl[ j h_j^{-2} \cdot \sum_{u=0}^{j-1} \frac {h_u^2 \zeps_u} {\theta_u} \biggl]\
=\ \theta_j \cdot jh_j^{-2} \cdot \biggl[ 1\,  +\, 
\sum_{u=1}^{j-1} \frac {h_u^2 \zeps_u} {\theta_u} \biggl]\ .     
\end{align}

    {\bf Sub-step 1D.\ \ A summary of some key properties 
of the parameters:}  \zhfb
(a) For each $j \in \N$, $\zeps_j \in (0, 1/9]$, $\theta_j \in (0, 1/9]$,
$\theta^*_j \in (0, 1/8]$, $I_j \geq 72$, and $h_j > 0$.  \zhfb
(b) As $j \to \infty$, $\zeps_j \to 0$ monotonically
and $\theta_j \to 0$ monotonically. \zhfb
(c)  As $j \to \infty$, $\theta^*_j \to 0$ monotonically
and $I_j \to \infty$ monotonically. \zhfb     
(d) For each $j \in \N$, one has that $\theta_j/\theta^*_j < 1$ and
$\theta_j^*\zeps_j I_j \leq 1$.   
\zhfb
(e) As $j \to \infty$, one has that $\theta_j/\theta^*_j \to 1$ and 
$\theta_j^*\zeps_j I_j \to 1$ and $\theta_j/h_j \to 0$. \zhfb
(f) One has that $\sum_{j=1}^\infty h_j^2 \zeps_j < 
\sum_{j=1}^\infty \zeps_j \leq 1/8$.
\zhfb
(g) For each $j \geq 2$, one has that $\zeps_j \leq 9^{-j}/I_{j-1}$ and that
$h_j^2 \zeps_j/\theta_j \geq j \cdot [1 + \sum_{u=1}^{j-1} (h_u^2 \zeps_u/\theta_u)] > j$.  
\medskip

     Everything in property (a) was noted in the paragraphs of 
eqs.\ (\ref{eq8.1B5}) and (\ref{eq8.1B6}) or comes from eq.\ (\ref{eq8.101}).
Properties (b) holds by the first two sentences of Sub-step 1C; and 
property (c) then follows by (\ref{eq8.1B6}).     
Property (d) holds by (\ref{eq8.1B6}) and property (a). 
In property (e), the first two limits hold by (\ref{eq8.1B6}) and properties (b)(c),
and the third limit holds by (\ref{eq8.101}) and the second sentence after
(\ref{eq8.1B5}).
Property (f) holds by (\ref{eq8.101}), 
the first sentence after (\ref{eq8.1B5}), 
and a simple calculation.
In property (g), the first inequality holds by the first sentence after
(\ref{eq8.1B5}), and the final pair of inequalities holds by
(\ref{eq8.1C1}) and property (a).
\medskip 
   
   {\bf Step 2.}\ \ Here we repeat essentially verbatim the 
entire ``Step 2'' in the proof given in Section \ref{sc7} for 
Theorem \ref{thm3.4}.
Where (\ref{eq7.1A1}) and Sub-step 1D(f) are cited in that context there, 
one simply cites the corresponding facts in (\ref{eq8.101}) and Sub-step 1D(f) here. 
\medskip

      {\bf Step 3.\ \ Proof of Statement (i)} [in Theorem
\ref {thm4.4}].
This is the same trivial argument as in Step 3 in the
proof (in Section \ref{sc7}) of Theorem \ref{thm3.4}.
[Again, (\ref{eq7.1A1}), cited there, is for $j \in \N$ the same as (\ref{eq8.101}) here.]
\medskip

   {\bf Step 4.}\ \ Again, this is identical to Step 4 in the 
proof (in Section {\sc 7}) of Theorem \ref{thm3.4}.    
Again, eq.\ (\ref{eq7.1A1}) and Sub-step 1D(a) cited there, are (for $j \in \N$) identical
to (\ref{eq8.101}) and Sub-step 1D(a) here. 
Of course (\ref{eq7.41}) holds in our context here.
Simply as noted in the final sentence of Step 4 there,
one has in our context here that 
$E[X_0^2] < \infty$ and $E[X_0] = 0$.  
\medskip   
       
   {\bf Step 5.  Proof of Statement (ii)}  
[in Theorem \ref{thm4.4}].
Let $n \in \N$ be arbitrary but fixed.
By (\ref{eq2.33}), it suffices to prove for this $n$ that
\begin{equation}
\label{eq8.51}
\beta_X(n)\ \leq\ (1/n) g_n\ .
\end{equation}

     To accomplish this, we shall employ for each $j \in \N$, 
in essentially the order here, the following: \zhfb 
{\bf (a)} eq.\ (\ref{eq2.42}) and Lemma \ref{lem6.5}(6); \zhfb 
{\bf (b)} the inequalities $6 < e^2$ and $0 < 1 - \theta_j < \exp(-\theta_j) < 1$; \zhfb 
{\bf (c)} from ({\ref {eq8.1B5}) the equality 
$\zeps_j = \exp({\bf L}_j(0)) \cdot \exp(B_j) \cdot e^{-j} \cdot e^{-2}$; \zhfb 
{\bf (d)} from (\ref {eq8.1B5}) and the last equality in (\ref{eq8.1B4b}),
the equality ${\bf L}_j(n) = {\bf L}_j(0) - \theta_j n$; \zhfb
{\bf (e)} from the second sentence after (\ref{eq8.1B3}),
the equality ${\bf L}_j(n) + B_j = L^{(t(j))}(n) + B_j$; \zhfb 
{\bf (f)} from eq.\ (\ref{eq6.114}) in Lemma \ref{lem6.11}, together with 
(\ref{eq8.1B2}), the inequality $L^{(t(j))}(n) + B_j \leq \phi(n) + \psi(n)$; \zhfb
{\bf (g)} the definitions of the functions $\phi$ and $\psi$ in Step 0. \zhfb    
Thereby, for each $j \in \N$, one has that
\begin{align}
\nonumber
\beta  &\Bigl(\sigma(X^{(j)}_0),\, \sigma(X^{(j)}_n) \Bigl)\
=\ \beta_{X(j)}(n)\
\leq\ 6 \zeps_j \cdot (1 - \theta_j)^n\
\leq\ e^2 \zeps_j \cdot (\exp(-\theta_j))^n \nonumber\\
&=\ \exp({\bf L}_j(0)) \cdot \exp(B_j) \cdot e^{-j} \cdot 
\exp(-\theta_j n)\ 
=\ e^{-j} \cdot \exp\bigl({\bf L}_j(0) -\theta_j n + B_j\bigl)\ 
=\ e^{-j} \cdot \exp \bigl({\bf L}_j(n) + B_j\bigl)\
\nonumber \\
&=\ e^{-j} \cdot \exp \bigl(L^{(t(j))}(n) + B_j\bigl)\
\leq\ e^{-j} \cdot \exp \bigl(\phi(n) + \psi(n) \bigl)\
=\ e^{-j} \cdot (1/n) \cdot g_n. \nonumber      
\end{align}

   Hence by Lemma \ref{lem6.7} and Sub-steps 2A and 2C 
(carried over directly from Section \ref{sc7}),
\begin{equation}
\nonumber 
\beta_X(n)\ =\ \beta\bigl(\sigma(X_0), \sigma(X_n)\bigl)\
\leq\ \sum_{j=1}^\infty
\beta\bigl(\sigma(X^{(j)}_0), \sigma(X^{(j)}_n)\bigl)\
\leq\ (1/n) g_n \cdot \sum_{j=1}^\infty e^{-j}\
<\ (1/n) g_n \cdot 1.   
\end{equation}
Thus (\ref{eq8.51}) holds.
 That completes the proof of Statement (ii). 
\medskip

     Our next task is to prove Statement (iii) in Theorem \ref{thm4.4} — that is,
to prove in our context here the conclusions (iii) and (iv) in the
statement of Theorem \ref{thm3.4}.
That will be the purpose of the next two ``steps''.
\vskip 0.1 in

   {\bf Step 6.  Proof (in our context here) of Statement (iii)
(including eq.\ (\ref{eq3.23})) in Theorem \ref{thm3.4}.}
The argument is exactly the same as in Step 6 in Section \ref{sc7}, 
in the proof of Theorem \ref{thm3.4}.
One can employ the same two-column chart there.
In particular, with regard to the right hand column of that chart, 
the ``second half'' of Sub-step 1D(g) here in Section \ref{sc8} is identical to  
Sub-step 1D(g) in Section \ref{sc7}, and
``Sub-steps 2A and 2B'' and ``(\ref{eq7.2C1})'' for Section \ref{sc8} here
are a direct repeat of Sub-steps 2A and 2B and (\ref{eq7.2C1}) in Section \ref{sc7}. 
\medskip

   {\bf Step 7.  Proof (in our context here) of 
Statement (iv) in Theorem \ref{thm3.4}.}    
Most of the argument will be a repeat of Step 7 in Section \ref{sc7}, in the proof
of Theorem \ref{thm3.4}.
As in Step 7 there, the argument in Step 7 here will be divided into four
``sub-steps''.
Only in the third sub-steps there and here will there be a real difference in the
arguments.
\smallskip

References to Sub-step 1D((a)-(g)) there in Section \ref{sc7} 
are simply replaced by references to the corresponding identical facts
in Sub-step 1D((a)-(g)) here in Section \ref{sc8}.  
\vskip 0.1 in   

   {\bf Sub-step 7A.}\ \ The argument is the same as in
Sub-step 7A in the proof (in Section \ref{sc7})    
of Theorem \ref{thm3.4}.
Thereby one (again) obtains here eq.\ (\ref{eq7.7A2}),
which we repeat for reference:
\begin{equation}
\label{eq8.7A1}
\frac {\theta_j} {h_j} \sum_{k=1}^{I(j)} h_j X^{(j)}_k\ \to\ \mu_{P1sL}\ 
{\rm in\ distribution\ as}\ j \to \infty. 
\end{equation}
\smallskip

   {\bf Sub-step 7B.}\ \ The argument is the same as in
Sub-step 7B in the proof (in Section \ref{sc7})    
of Theorem \ref{thm3.4}.
Thereby one (again) obtains here eq.\ (\ref{eq7.7B1}),
which we repeat for reference:
\smallskip
\begin{equation}
\label{eq8.7B1}
\frac {\theta_j} {h_j} 
\sum_{u=1}^{j-1} \sum_{k=1}^{I(j)} h_u X^{(u)}_k\
\to\ 0\ \, {\rm in\ probability\ as}\ j \to \infty.
\end{equation}

   {\bf Sub-step 7C.}\ \ The argument here differs from
that in Sub-step 7C in the proof of Theorem \ref{thm3.4}.
It is identical to the argument in Sub-step 7C of Section 7 of
\cite {ref-journal-Bradley2025},
and is repeated here directly for convenience:
\smallskip

   If $j$ and $u$ are any positive integers such that $u \geq j+1$,
then $I_{u-1} \geq I_j$ by Sub-step 1D(c), and hence by Sub-step 1D(g)
and the first sentence of Sub-step 2B (which was carried over directly 
from Section \ref{sc7} to here in Section \ref{sc8}) one has that 
$P(X^{(u)}_k \neq 0) = \zeps_u 
\leq 9^{-u}/I_{u-1} \leq 9^{-u}/I_j$.
Hence for any positive integer $j$, 
\begin{align}
\nonumber
P\biggl( \frac {\theta_j} {h_j}
\sum_{u=j+1}^\infty \sum_{k=1}^{I(j)} 
h_uX^{(u)}_k \neq 0 \biggl)\
&\leq\ \sum_{u=j+1}^\infty \sum_{k=1}^{I(j)}
P(X^{(u)}_k \neq 0)\
\leq\ \sum_{u=j+1}^\infty \sum_{k=1}^{I(j)} 
(9^{-u}/I_j)\  \\
\nonumber
&\leq\ \sum_{u=j+1}^\infty 9^{-u}\
=\ 9^{-(j+1)}/(1 - 9^{-1})\ <\ 9^{-j}.   
\end{align}
Hence 
$
P\bigl((\theta_j/h_j) \sum_{u=j+1}^\infty \sum_{k=1}^{I(j)} 
h_uX^{(u)}_k \neq 0\bigl) \to 0
$ 
as $n \to \infty$.
As a consequence,
\begin{equation}
\label{eq8.7C1}
\frac {\theta_j} {h_j} \sum_{u=j+1}^\infty \sum_{k=1}^{I(j)} h_uX^{(u)}_k 
\to\ 0\, \ {\rm in\ probability\ as}\  j \to \infty.
\end{equation} 

   {\bf Sub-step 7D.}\ \ The argument is the same as in
Sub-step 7D in the proof (in Section {\ref{sc7}) of
Theorem \ref{thm3.4} — simply referring here to
eqs.\ (\ref{eq8.7A1}), (\ref{eq8.7B1}), and (\ref{eq8.7C1}), 
in place of the respectively identical equations  
(\ref{eq7.7A2}), (\ref{eq7.7B1}), and (\ref{eq7.7C1}) there.
That completes the proof here of conclusion (iv) in Theorem \ref{thm3.4}. 
\medskip

     {\bf Step 8.}\ \ Our final task is to prove conclusion (iv) in the statement of
Theorem \ref{thm4.4}.  
From the formulation of that conclusion (iv), recall for each positive integer $n$ 
the definition ${\bf h}(n) := n^{-1} \var(\sum_{k=1}^n  X_k)$.
In Step 8 here the proof of eq.\ (\ref{eq4.41}) will be given.
In Step 9 below, the proof of eq.\ (\ref{eq4.42}) will be given, and the proof of
Theorem \ref{thm4.4} will then be complete.   
\medskip     
     
    {\bf Proof of eq.\ (\ref{eq4.41}).}\ \ This proof will be divided into
two ``sub-steps'', of which the first is labeled ``Sub-claim 8A''.
\medskip

     {\bf Sub-claim 8A.}\ \ {\it Suppose $j$ and $n$ are each a positive integer.
Then} 
\begin{equation}
\label{eq8.8A1}
n^{-1}\, \var\biggl(\, \sum_{k=1}^n X_k^{(j)} \biggl)\ \leq\ 2^{-j} g_n.
\end{equation}

    {\bf Proof of Sub-claim 8A.}\ \ First a preliminary observation is in order:
By (\ref{eq8.1B4c}), ${\bf L}_j(0)\ =\ 1 + \log(1/t_j)$. 
Hence by both ``halves''  of (\ref{eq8.1B5}), 
\begin{align}
\nonumber
\zeps_j\ &=\ \exp \bigl( {\bf L}_j(0) \bigl)\, \cdot\, \exp \bigl( B_j - j - 2 \bigl)\
=\ \exp \bigl( 1 + \log(1/t_j) \bigl)\, \cdot\, \exp \bigl( B_j - j - 2 \bigl)\\
\nonumber 
&=\ e \cdot (1/t_j)\, \cdot\, \exp \bigl( B_j - j - 2 \bigl)\
=\ \theta_j\, \cdot\, \exp \bigl( B_j - j - 1 \bigl).    
\end{align}
Hence
\begin{equation}
\label{eq8.8A21}
\zeps_j / \theta_j\ =\ \exp(B_j - j - 1).
\end{equation}

   Now the rest of the argument for Sub-claim 8A will be divided into two cases
according to whether $g_n < 2^{j+1} \exp(B_j)$ or $g_n \geq 2^{j+1} \exp(B_j)$. 
\medskip 

   {\bf Case 1.}\ \ $g_n < 2^{j+1} \exp(B_j)$.\ \ 
By (\ref{eq8.1B3}) and (\ref{eq8.1B5}), and then (\ref{eq8.1B1a}), 
$\theta_j \leq {\bf c}_j \leq 2^{-j} \cdot (g_n/n) \cdot  \exp(-B_j)$.
Hence by (\ref{eq8.8A21}),
\begin{equation}
\nonumber 
\zeps_j\ =\ (\zeps_j/\theta_j) \cdot \theta_j\
\leq\ \exp(B_j-j-1) \cdot 2^{-j} \cdot (g_n/n) \cdot \exp(-B_j)\
<\ 2^{-j} \cdot (g_n/n). 
\end{equation} 
Hence by Minkowski's inequality and Lemma \ref{lem6.5}(4),
\begin{equation}
\nonumber
(1/n) \cdot  \var\biggl(\, \sum_{k=1}^n X_k^{(j)} \biggl)\ \leq\ (1/n) \cdot n^2\, \var(X^{(j)}_0)\
=\ n \zeps_j\ \leq\ 2^{-j} g_n.
\end{equation} 
Thus (\ref{eq8.8A1}) holds for Case 1.
\medskip

     {\bf Case 2.}\ \ $g_n \geq 2^{j+1} \exp(B_j)$.\ \ 
By Lemma \ref{lem6.5}(9) and eq.\ (\ref{eq8.8A21}), 
\begin{equation}
\nonumber 
n^{-1} \var \biggl(\, \sum_{k=1}^n X_k^{(j)} \biggl)\ \leq\ 2 \zeps_j / \theta_j\
=\ 2\, \exp(B_j - j - 1)\ <\ 2g_n/2^{j+1}\ =\ 2^{-j} g_n. 
\end{equation} 
Thus (\ref{eq8.8A1}) holds for Case 2.
That completes the proof of Sub-claim 8A.
\medskip

   {\bf Sub-step 8B.}\ \
Recall (again) the definition of the numbers ${\bf h}(n)$, $n \in N$ in the statement of 
Theorem \ref{thm4.4}.
By (\ref{eq7.2C1}) (which carries over to Section 8 here) and then 
Lemma \ref{lem6.8}, eq.\ (\ref{eq8.101}), and Sub-claim 8A, for each positive integer $n$,
\begin{align}
\nonumber
{\bf h}(n)\ &=\ n^{-1} \var \biggl(\, \sum_{k=1}^n X_k \biggl)\ 
=\ n^{-1} \sum_{j=1}^\infty \var \biggl(\, \sum_{k=1}^n h_j X_k^{(j)} \biggl)\
=\ n^{-1} \sum_{j=1}^\infty h_j^2\, \var \biggl(\, \sum_{k=1}^n X_k^{(j)} \biggl)\\
\nonumber
&\leq\ n^{-1} \sum_{j=1}^\infty \var \biggl(\, \sum_{k=1}^n X_k^{(j)} \biggl)\
\leq\ \sum_{j=1}^\infty 2^{-j}g_n\ =\ g_n.    
\end{align}   
Thus eq.\ (\ref{eq4.41}) holds.  Step 8 is complete.
\medskip

    {\bf Step 9.\ \ Proof of eq.\ (\ref{eq4.42}).}\ \
The argument will be divided into three ``sub-steps'', of which one is a ``sub-claim''.
\medskip
      
     {\bf Sub-step 9A.}\ \ 
For each positive integer $j$, define the positive integer  $M_j$ [also denoted $M(j)$] 
as follows:
\begin{equation}
\label{eq8.9A1}
M_j\ =\ \biggl[\, \frac {1} {j^{1/2} \theta_j} \biggl], 
\end{equation}
that is, the greatest integer that is $\leq 1/(j^{1/2} \theta _j)$.
For each positive integer $j$, by trivial arithmetic and the second sentence 
after (\ref{eq8.1B5}), one has that $\theta_j \leq 1/9^j < 1/j$
and hence $j^{1/2} \theta_j \leq j\theta_j < 1$.
Hence the integer $M_j$ is in fact positive for each $j \in \N$;
and further (again refer to Notations \ref{nt2.1}), 
\begin{equation}
\label{eq8.9A1a}
j^{1/2} \theta_j \to 0\ \ \ {\rm and}\ \ \ M_j \sim 1/(j^{1/2} \theta_j) \to \infty\ \ \ {\rm as}\  
j \to \infty. 
\end{equation}

     Recall the first sentence of Remark \ref{rem4.5}(D).
To show that (\ref{eq4.42}) holds, it suffices to show that for every positive integer $N$,
\begin{equation}
\label{eq8.9A2}
\var\biggl(\, \sum_{k=1}^{N \cdot M(j)} X_k\biggl)\ 
\sim\ N^2 \cdot \var\biggl(\, \sum_{k=1}^{M(j)} X_k\biggl)\ \ \ 
{\rm as}\ j \to \infty. 
\end{equation} 
The following lemma is intended to help accomplish that:
\medskip

   {\bf Sub-claim 9B.}\ \ {\it For every positive integer $K$, the following statements hold:}
\begin{align}
\label{eq8.9B1a}
\var\biggl(\, \sum_{k=1}^{K \cdot M(j)} h_j X_k^{(j)} \biggl)\ 
&\sim\ K^2 M_j^2 h_j^2 \zeps_j\ \ \ {\rm as}\ j \to \infty.\\
\label{eq8.9B1b}
\sum_{u=1}^{j-1} \var \biggl(\, \sum_{k=1}^{K \cdot M(j)} h_u X_k^{(u)} \biggl)\
&=\ o(M_j^2 h_j^2 \zeps_j) \ \ \ {\rm as}\ j \to \infty.\\
\label{eq8.9B1c}
\sum_{u=j+1}^\infty \var \biggl(\, \sum_{k=1}^{K \cdot M(j)} h_u X_k^{(u)} \biggl)\
&=\ o(M_j^2 h_j^2 \zeps_j)\ \ \ {\rm as}\ \ j \to \infty.\\
\label{eq8.9B1d}
\var\biggl(\, \sum_{k=1}^{K \cdot M(j)} X_k\biggl)\ 
&\sim\ K^2 M_j^2 h_j^2 \zeps_j\ \ \ {\rm as}\ j \to \infty.   
\end{align}

   {\bf Proof of (\ref{eq8.9B1a}).}\ \ Suppose $K \in \N$.  One has that for each $j \in \N$,
by Lemma \ref{lem6.5}(4)(8a),   
\begin{align}   
\nonumber  
\var\biggl(\, \sum_{k=1}^{K M(j)} &h_j X_k^{(j)} \biggl)\ 
=\ \sum_{k=1}^{K M(j)} \sum_{\ell = 1}^{K M(j)} 
\cov\Bigl(h_j X_k^{(j)},\, h_j X_\ell^{(j)}\Bigl)\
=\ h_j^2\,  \sum_{k=1}^{K M(j)} \sum_{\ell = 1}^{K M(j)}  \cov\Bigl(X_k^{(j)},\, X_\ell^{(j)}\Bigl)\\
\label{eq8.9B2}
&=\ h_j^2\,  \sum_{k=1}^{K M(j)} \sum_{\ell = 1}^{K M(j)}  \zeps_j (1 - \theta_j)^{|k-\ell |}\
\leq\ h_j^2\,  \sum_{k=1}^{K M(j)} \sum_{\ell = 1}^{K M(j)}  \zeps_j \cdot 1\
=\ K^2 M_j^2 h_j^2 \zeps_j.    
\end{align}
Also, after repeating the first three equalities of (\ref{eq8.9B2}), 
one has that for each $j \in \N$, by Lemma \ref{lem6.5}(4)(8a), 
\begin{align}
\nonumber
\var\biggl(\, \sum_{k=1}^{K M(j)} &h_j X_k^{(j)} \biggl)\ 
=\ h_j^2\,  \sum_{k=1}^{K M(j)} \sum_{\ell = 1}^{K M(j)}  \zeps_j (1 - \theta_j)^{|k-\ell |}\
\geq\ h_j^2\,  \sum_{k=1}^{K M(j)} \sum_{\ell = 1}^{K M(j)}  \zeps_j (1 - \theta_j)^{KM(j)}\\ 
\label{eq8.9B3}
&\geq\ h_j^2\,  \sum_{k=1}^{K M(j)} \sum_{\ell = 1}^{K M_j}  
\zeps_j \Bigl(1\, -\, \theta_j  K M_j \Bigl)\
=\ K^2 M_j^2 h_j^2 \zeps_j  \Bigl(1\, -\, \theta_j K M_j \Bigl).     
\end{align}
By (\ref{eq8.9A1}), $M_j K \theta_j \to 0$ as $j \to \infty$.
Now (\ref{eq8.9B1a}) follows from (\ref{eq8.9B2}) and (\ref{eq8.9B3}).
\medskip

   {\bf Proof of (\ref{eq8.9B1b}).}\ \ Suppose $K \in \N$.  One has that for each $j \geq 2$
and each $u \in \{1,2,\dots, j-1\}$, by Lemma \ref{lem6.5}(9), 
\begin{equation}
\nonumber
\var \biggl(\, \sum_{k=1}^{KM(j)} h_u X_k^{(u)} \biggl)\ 
=\ h_u^2\, \var\biggl(\, \sum_{k=1}^{KM(j)} X_k^{(u)}\biggl)\
<\ h_u^2 K M_j \cdot 2 \zeps_u / \theta_u. 
\end{equation}
Hence considering integers $j \geq 2$, one has by Sub-step 1D(g) and (\ref{eq8.9A1a}) that
\begin{align}
\nonumber
\sum_{u=1}^{j-1}&\, \var \biggl(\, \sum_{k=1}^{KM(j)} h_u X_k^{(u)} \biggl)\
\leq\ 2KM_j \cdot \sum_{u=1}^{j-1} h_u^2 \zeps_u / \theta_u\
\leq\ 2KM_j h_j^2 \zeps_j / (j \theta_j)\\
\nonumber
&=\ 2KM_j^2 h_j^2 \zeps_j \cdot \bigl(j^{-1/2} \bigl) \cdot 1/\bigl(M_j \cdot j^{1/2} \theta_j\bigl)\
\sim\ 2KM_j^2 h_j^2 \zeps_j \cdot \bigl(j^{-1/2} \bigl) \cdot 1
=\ o(M_j^2 h_j^2 \zeps_j) \ \ \ {\rm as}\ j \to \infty. \quad   
\end{align}
Thus (\ref{eq8.9B1b}) holds.
\medskip

   {\bf Proof of (\ref{eq8.9B1c}).}\ \ Suppose $K \in \N$.  One has that for each $j \in \N$
and each $u\geq j+1$, by Minkowski's inequality and Lemma \ref{lem6.5}(4),  
\begin{equation}
\label{eq8.9B4}
\var \biggl(\, \sum_{k=1}^{K M(j)} h_u X_k^{(u)} \biggl)\ 
=\ h_u^2\, \var \biggl(\, \sum_{k=1}^{K M(j)} X_k^{(u)} \biggl)\
\leq\ h_u^2 (KM_j)^2\, \var(X_0^{(u)})\ =\ h_u^2 K^2 M_j^2 \zeps_u.  
\end{equation}     
Also, for each $j \in \N$ and each $u \geq j+1$, one trivially has that
$h_u < h_j$ by (\ref{eq8.101}), and $\zeps_u \leq (2^{-j})^{u-j} \zeps_j$ 
by the first sentence of Sub-step 1C (here in Section \ref{sc8}) 
and a crude induction argument.
Hence considering $j \in \N$, one has by (\ref{eq8.9B4}) that 
\begin{align}
\nonumber
\sum_{u=j+1}^\infty \var \biggl(\, \sum_{k=1}^{KM(j)} h_u X_k^{(u)} \biggl)\
&\leq\ \sum_{u= j+1}^\infty h_j^2 K^2 M_j^2 \cdot (2^{-j})^{u-j} \zeps_j\
=\ [2^{-j}/(1 - 2^{-j})] \cdot K^2 M_j^2 h_j^2 \zeps_j\\
\label{eq8.9B5}
&=\ o(M_j^2 h_j^2 \zeps_j)\ \ \ {\rm as}\ \ j \to \infty.   
\end{align}
Thus (\ref{eq8.9B1c}) holds.
\medskip

    {\bf Proof of (\ref{eq8.9B1d}).}\ \ Suppose $K \in \N$.  Eq.\ (\ref{eq8.9B1d}) holds by
eqs.\ (\ref{eq8.9B1a}), (\ref{eq8.9B1b}), and (\ref{eq8.9B1c}), and Lemma \ref{lem6.8},
the independence of the Markov chains $X^{(u)}$ for $u \in \N$, and the fact
that the random variables $X_k^{(u)}$ are centered.
That completes the proof of Sub-claim 9B.
\smallskip

     {\bf Sub-step 9C.}\ \ Now we are ready to prove (\ref{eq4.42}).
Suppose $N \in \N$.
By Sub-claim 9B [its eq.\ (\ref{eq8.9B1d})] with $K=1$, one has that
$\var(\sum_{k=1}^{M(j)} X_k) \allowbreak \sim M_j^2 h_j^2 \zeps_j$ 
as  $j \to \infty$.
Hence by trivial multiplication, 
$N^2\, \var(\sum_{k=1}^{M(j)} X_k) \allowbreak \sim N^2 M_j^2 h_j^2 \zeps_j$ 
as $j \to \infty$.   
Combining that with   
Sub-claim 9B [its eq.\ (\ref{eq8.9B1d})] for $K = N$, one has that
eq.\ (\ref{eq8.9A2}) holds.
Thus eq.\ (\ref{eq4.42}) holds.
That completes the proof of Theorem \ref{thm4.4}.

\section{Proof of Theorem \ref{thm5.5}}
\label{sc9}
The proof will be divided into a series of ``steps'',
starting with ``Step 0''.
\medskip

     {\bf Step 0.}\ \ This step will be divided into several ``sub-steps''.
\smallskip
     
     {\bf Sub-step 0A.}\ \ Referring to the function $f$ in the statement 
of Theorem \ref{thm5.5}, define the function 
$\phi: [1, \infty) \to (-\infty, 0]$ as follows:
\begin{equation}
\label{eq9.0A1}
{\rm For\ each}\ x \in [1, \infty), \indent
\phi(x)\ :=\ \log f(x).  
\end{equation}
By the hypothesis (of Theorem \ref{thm5.5}),
this function $\phi$ is continuous, strictly decreasing, and
convex on $[1, \infty)$, and  
\begin{equation}
\label{eq9.0A2}
\lim_{x \to \infty} \phi(x)\ =\ -\infty 
\indent {\rm and} \indent
\lim_{x \to \infty} \frac {\phi(x)} {x}\ =\ 0.  
\end{equation}
In fact by the hypothesis (of the theorem) and elementary 
arguments [involving convexity and (\ref{eq9.0A2})],
the function $\phi$ is twice continuously differentiable on the open
half line $(1, \infty)$, its derivative $\phi'(x)$ is nondecreasing 
on that open half line, and      
\begin{equation}
\label{eq9.0A3}
\frac {f'(x)} {f(x)}\ =\ 
\phi'(x)\ <\ 0\ \ {\rm for\ every}\ x \in (1, \infty), 
\indent {\rm and} \indent
\lim_{x \to \infty} \phi'(x)\ =\ 0.
\end{equation}

   From the preceding paragraph, one has that
\begin{align}
\label{eq9.0A4a}
&-\phi'(x)\ >\ 0\ \ {\rm for\ all}\ x \in (1,\infty);\ \ {\rm and} \\
\label{eq9.014b}
&-\phi'(x)\ {\rm is\ nonincreasing\ as}\ x\ {\rm increases\ in}\ (1,\infty);\ \  
{\rm and} \lim_{x \to \infty} (-\phi'(x))\ =\ 0. 
\end{align}
Accordingly, referring to the function $g$ in the statement of
Theorem \ref{thm5.5} and applying Lemma \ref{lem6.13} 
[with the function $f$ there being $x \mapsto -\phi'(x)$ here],
we assume without loss of generality that
\begin{equation}
\label{eq9.0A5} 
{\rm the\ product}\ \ g(x) \cdot (-\phi'(x))\ \  
{\rm is\ nonincreasing\ as}\ x\ {\rm increases\ in}\ (1,\infty). 
\end{equation}
Eq.\ (\ref{eq9.0A5}) will come into play somewhat later in the proof.
\medskip

   {\bf Note.}\ \ {\it The inequality $-\phi'(x) > 0$ for $x \in (1,\infty)$  
from (\ref{eq9.0A4a}) should be kept in mind in everything that follows.}\ \ 
In particular, the function $x \mapsto \log (-\phi'(x))$
for $x \in (1,\infty)$ will shortly come into play in calculations below.
\medskip

   {\bf Sub-step 0B.}\ \ 
For each $y \in (1, \infty)$ let $L^{(y)}: \R \to \R$ 
denote the affine function whose graph is tangent to
the graph of the (convex) function $\phi$ at the point
$(y, \phi(y))$.
Thus for each $y \in (1, \infty)$, one has that (see (\ref{eq9.0A3})) 
\begin{align}
\label{eq9.0A6a}
&L^{(y)}(y) = \phi(y) 
\indent {\rm and} \indent
L^{(y)}(x)\ \leq\ \phi(x)\ \ {\rm for\ all}\ x \in [1, \infty);\ \quad {\rm and} \\
\label{eq9.0A6b}
&[{\rm slope\ of}\ L^{(y))}]\ =\ \phi'(y)\ <\ 0.
\end{align}

     By (\ref{eq9.0A6a})--(\ref{eq9.0A6b}),  
for each $y \in (1, \infty)$,
\begin{equation}
\label{eq9.0B1}
L^{(y)}(0)\ =\ 
L^{(y)}(y) - y \cdot [{\rm slope\ of}\ L^{(y)}]\
=\ \phi(y)\ +\ y \cdot (-\phi'(y)).    
\end{equation} 
Sub-steps 0C and 0D below will examine the right hand side of
(\ref{eq9.0B1}) in more detail, for later use.
\medskip 

   {\bf Sub-step 0C.}\ \ Let ${\bf w} \in (1,\infty)$ be fixed such that
it satisfies the condition in hypothesis (e) in the theorem.  
For each $x \in [{\bf w}, \infty)$, by 
hypotheses (b), (c), and (e) in the theorem, one has that
\begin{equation}
\nonumber
f(x)\ +\ x \cdot f'(x)\ =\ 
\frac {d} {dx} [x \cdot f(x)]\ \leq\ 0
\end{equation}
and hence $f(x) \leq x \cdot (-f'(x))$.
Hence by (\ref{eq9.0A3}) (and hypothesis (a) in the theorem),
\begin{equation}
\label{eq9.0C1}
{\rm for\ each}\ x \in [{\bf w}, \infty),\ \ \   
1\ \leq\ x \cdot (-f'(x))/f(x)\ =\ x \cdot (-\phi'(x)).
\end{equation}

   {\bf Sub-step 0D.}\ \
Referring to (\ref{eq9.0B1}), define the function
$\xi: (1,\infty) \to \R$ as follows: 
\begin{equation}
\label{eq9.0D1} 
{\rm For\ each}\ x \in (1, \infty), \indent
\xi(x)\ :=\ \exp \bigl( L^{(x)}(0) \bigl)\ =\  
\exp \Bigl( \phi(x)\ +\ x\cdot (-\phi'(x)) \Bigl). 
\end{equation}
For every $x \in (1, \infty)$, by (\ref{eq9.0D1}) and (\ref{eq9.0A4a}),
one has that 
$\xi(x) > 0$ and $-\phi'(x) > 0$ and hence $\xi(x)/(-\phi'(x)) > 0$.
We shall now examine that latter ratio in more detail.
\smallskip

By (\ref{eq9.0D1}), one has that
for every $x \in (1, \infty)$, 
$\log \xi(x) = \phi(x) + x \cdot (-\phi'(x))$
and hence
\begin{equation}
\nonumber
\frac {d}{dx} \Bigl( \log \xi(x) \Bigl)\ =\ -x \cdot \phi''(x). 
\end{equation}
Hence for every $x \in (1, \infty)$,
\begin{align}
\label{eq9.0D2}
\frac {d} {dx} \Biggl[ \log \Biggl(
\frac {\xi(x)} {-\phi'(x)} \Biggl) \Biggl]\
=\ 
\frac {d} {dx} \Bigl[ \bigl(\log \xi(x) \bigl) 
- \log \bigl(-\phi'(x) \bigl) \Bigl]\
&=\ -x \cdot \phi''(x)\ -\ \frac {1} {-\phi'(x)} 
\cdot \bigl(-\phi''(x) \bigl) \nonumber\\
&=\ \biggl[ -x\ +\ \frac {1} {-\phi'(x)} \biggl]\ \cdot\ \phi''(x).
\end{align}
By (\ref{eq9.0C1}), for each $x \in [{\bf w}, \infty)$, 
one has that $1/(-\phi'(x)) \leq x$ and hence
$-x + 1/(-\phi'(x)) \leq 0$.
Since the function $\phi$ is convex, one also has that
$\phi''(x) \geq 0$ for each $x \in (1, \infty)$.
Hence for each $x \in [{\bf w}, \infty)$, one has that
$[-x + 1/(-\phi'(x))] \cdot \phi''(x) \leq 0$.
Hence by (\ref{eq9.0D2}), the expression
$\log[\xi(x)/(-\phi'(x))]$ is nonincreasing    
as $x$ increases in $[{\bf w}, \infty)$.
Hence 
\begin{equation}
\label{eq9.0D3}
{\rm the\ fraction}\ \ \xi(x)/(-\phi'(x))\ \ 
{\rm is\ nonincreasing\ as}\ x\ 
{\rm increases\ in}\ [{\bf w}, \infty).
\end{equation}
Accordingly, as a point of reference,
define the positive number $\Delta$ by
\begin{equation}
\label{eq9.0D4}
\Delta\ :=\ \frac {\xi({\bf w})} {-\phi'({\bf w})}\
=\ \sup_{x \in [{\bf w}, \infty)} \frac {\xi(x)} {-\phi'(x)}.
\end{equation}
(Of course $\Delta > 0$, by the sentence after (\ref{eq9.0D1}).)
\smallskip

That completes Step 0.
The quantity $\Delta$ in (\ref{eq9.0D4}) will come into play shortly.
\medskip

   {\bf Step 1.  The parameters.}\ \
In this step, we shall recursively define for each 
positive integer $j$ ten parameters
(one of those ``parameters'' being an affine function). 
They are listed here, along with their most basic restrictions: 
\medskip

    {\bf Key List:} \zhfb
(i) a positive number $B_j$; \zhfb
(ii) a positive number ${\bf c}_j$; \zhfb
(iii) a number $\zeps_j^* \in (0, 1/9]$; \zhfb 
(iv) a number $t_j \in [{\bf w}, \infty)$; \zhfb
(v) an affine function ${\bf L}_j: \R \to \R$; \zhfb
(vi) a number $\zeps_j \in (0, 1/9]$; \zhfb
(vii) a number $\theta_j \in (0, 1/9]$; \zhfb
(viii) a number $\theta^*_j \in (0, 1/8]$; \zhfb  
(ix) a positive integer $I_j$; and \zhfb
(x) a positive number $h_j$.
\medskip

   Such parameters $B_j$, $\zeps_j$, $\theta_j$, $I_j$, and $h_j$ 
in (i) and (vi)-(vii) and (ix)-(x) will also be defined for $j=0$.
Those five parameters for $j=0$ will provide a convenient ``springboard'' 
to start the recursion (but they will not otherwise play a significant
role later on in the construction itself for Theorem \ref{thm5.5}). 
\smallskip

     In superscripts, the number $t_j$ will be denoted $t(j)$ for typographical 
convenience.
 \medskip

   {\bf Sub-step 1A.\ \  The initial step.}\ \
We start with $j=0$, for just five parameters as mentioned above.
Define the following five parameters:
\begin{equation}
\label{eq9.1A1}
B_0 := 1 \quad {\rm and} \quad \zeps_0\ :=\ \theta_0\ :=\ 1/9  \quad {\rm and} \quad
I_0\ :=\ h_0:\ =\ 1.
\end{equation}
Obviously in the Key List above, the relevant basic restrictions 
in (i) and (vi)-(vii) and (ix)-(x) for the five parameters defined in 
(\ref{eq9.1A1}) are satisfied.
That completes the initial step.
\medskip

   {\bf Sub-step 1B.\ \  The recursion step.}\ \ 
Now suppose $j \geq 1$ is an integer, and that for each 
integer $u$ such that $0 \leq u \leq j-1$, the five parameters 
$B_u$, $\zeps_u$, $\theta_u$, $I_u$, and $h_u$ 
have been defined, meeting the basic requirements in (i) and (vi)-(vii) and (ix)-(x) 
in the Key List.
Our task now is to define, for the given $j$ here, all ten parameters in
the Key List, meeting the basic requirements there.
\medskip   

   To start off, referring to (\ref{eq9.0D4}), define the number $B_j$ by
\begin{equation}
\label{eq9.1B1}
B_j\ :=\ \max \biggl\{ 9h_{j-1}^2 \Delta\, ,\
 j \sum_{u=0}^{j-1} h_u^2 \zeps_u / \theta_u 
\biggl\}. 
\end{equation}
Of course $B_j > 0$.
[Recall that $\Delta > 0$ from the sentence after (\ref{eq9.0D4}), and
that from the recursive assumptions (from the Key List), the numbers
$\zeps_u$, $\theta_u$, and $h_u$ are positive for each $u \in \{0, \dots, j-1\}$.]
\smallskip

   Now refer to the assumptions on the function $g(x)$, $x \in [1, \infty)$ in the statement of Theorem \ref{thm5.5}, and the embellishment in eq.\ (\ref{eq9.0A5}).
For each $j \in \Z$, define the number ${\bf c}_j$ by
\begin{equation}
\label{eq9.1B1a}
{\bf c}_j\ :=\ \min_{n \in {\bf U}(j)} (2^{-j} \cdot n^{-1} g(n) \cdot B_j^{-1})
\end{equation}
where ${\bf U}(j) := \{1\} \cup \{i \in \N: g(i) < 2^{j+1} B_j\}$.
That set ${\bf U}_j$ is nonempty (as ensured by the ``$\{1\} \cup \dots$''), and it is finite.
The number ${\bf c}_j$ is well defined, and it is positive.
\smallskip  
 
   Next, let $\zeps_j^* \in (0, 1/9]$ be such that
\begin{equation}
\label{eq9.1B2} 
\zeps_j^*\ \leq\ \min \{ \zeps_{j-1}/3,\  9^{-j}/I_{j-1} \}.
\end{equation}

   Next, we apply Lemma \ref{lem6.11}(II) to the (convex) function $\phi$ here.
Note that $\log \zeps_j^* < 0$ by (\ref{eq9.1B2}).    
Referring to the number ${\bf w}$ in Sub-steps 0C and 0D, 
the notation $T(\dots)$ in Lemma \ref{lem6.11}(II), 
and the assumption in the theorem that $g(x) \to \infty$ as $x \to \infty$,  
let $t_j \in [{\bf w}, \infty)$ be a positive number with the following two properties:  
\begin{align}
\label{eq9.1B3a}
&t_j\ \geq\ T\bigl(\phi,\ \log \zeps_j^*,\ \min\{9^{-j},\ \theta_{j-1}/2,\ 
2^{-j} B_{j-1} \theta_{j-1} / B_j,\ {\bf c}_j \}\bigl); 
\indent {\rm and}\\
\label{eq9.1B3b}
&{\rm for\ all}\ x \in [t_j, \infty),\ \ g(x)\ \geq\ 2^{j+6} e^{j+2} B_j\ .
\end{align}

   Now referring to the entire paragraph of  (\ref{eq9.0A6a})-(\ref{eq9.0A6b}), 
define the affine function ${\bf L}_j: \R \to \R$ as follows:
\begin{equation}
\label{eq9.1B4}
{\rm For\ all}\ x \in \R,\ \ \ {\bf L}_j(x)\ =\ -(j+2)\ +\  L^{(t(j))}(x). 
\end{equation} 

   Note that ${\bf L}_j(0) < L^{(t(j))}(0) \leq \log \zeps^*_j$ from (\ref{eq9.1B4}),
(\ref{eq9.1B3a}), and Lemma \ref{lem6.11}(II).
It follows that  
$\exp({\bf L}_j(0)) < \zeps^*_j \leq 9^{-j}/I_{j-1} \leq 9^{-j} \leq 1/9$ by (\ref{eq9.1B2}).
Note also that from (\ref{eq9.1B3a}) and (again) Lemma \ref{lem6.11}(II),
$-1/9 \leq  \phi'(t_j) < 0$.
Accordingly, define the positive numbers 
$\zeps_j \in (0, 1/9]$, $\theta_j \in (0, 1/9]$, and
$\theta_j^*$ [\`a la (\ref{eq6.312})], the positive integer $I_j$, 
and the positive number $h_j$, as follows:    
\begin{align}
\label{eq9.1B5}    
&\zeps_j\ :=\ \exp\bigl({\bf L}_j(0)\bigl) \quad {\rm and} \quad
\theta_j\ := -\phi'(t_j) \quad {\rm and} \quad
\theta_j^*\ :=\ \frac {\theta_j} {1 - \zeps_j}\ ; 
\quad {\rm and}\\
\label{eq9.1B6}
&I_j\ :=\ \Biggl[ \frac {1} {\theta_j^* \zeps_j} \Biggl] 
\indent {\rm and} \indent
h_j\ :=\ \Bigl( B_j \theta_j / \zeps_j \Bigl)^{1/2}\ ,     
\end{align}
where in the first equality in (\ref{eq9.1B6}), the big brackets indicate the
greatest integer that is $\leq 1/(\theta_j^* \zeps_j)$.
By simple arithmetic, $\theta^*_j \in (0, 1/8]$ and $I_j \geq 72$. 
And since $B_j > 0$ (as mentioned above), one has that $h_j > 0$. 
The parameters defined here in Sub-step 1B meet all of the specified 
requirements in the Key List for the given integer $j$ in this recursion.
That completes the recursion step.
\smallskip

   That completes the recursive definition of the parameters. 
\medskip

   {\bf Sub-step 1C.\ \  Miscellaneous observations.}\ \ 
Here is a list of several observations that 
come from this recursion and will be employed later on.
\medskip

   First, from (\ref{eq9.1A1}), followed by (\ref{eq9.1B5}), 
the (first or) second sentence after (\ref{eq9.1B4}), and then 
(\ref{eq9.1B2}), one has that  
\begin{equation}
\label{eq9.1C1}
\zeps_0 = 1/9 \indent 
{\rm and\ for\ every}\ j  \in \N,\ \ 
\zeps_j\, =\, \exp({\bf L}_j(0))\, <\, \zeps_j^*\, \leq\, 
\min\{\zeps_{j-1}/3,\ 9^{-j}/I_{j-1}\}.
\end{equation}

   Next, from (\ref{eq9.1B3a}) and Lemma \ref{lem6.11}(II), for each $j \in \N$, 
$-\min\{9^{-j},\ \theta_{j-1}/2,\ {\bf c}_j\} \leq \phi'(t_j) < 0$.   Hence by (\ref{eq9.1B5}), 
\begin{equation}
\label{eq9.1C2}
{\rm for\ every}\ j \in \N,\ \ \ \theta_j\, \leq\, \min\{9^{-j},\ \theta_{j-1}/2,\ {\bf c}_j\}.
\end{equation}

   Next, by (\ref{eq9.1B6}) and (\ref{eq9.1B1}) and trivial arithmetic,
\begin{equation}
\label{eq9.1C3}
{\rm for\ every}\ j \in \N, \quad 
h_j^2 \zeps_j/\theta_j\ =\ B_j\ \geq\
j \cdot \sum_{u=0}^{j-1} h_u^2 \zeps_u/\theta_u.
\end{equation}

   Next, by (\ref{eq9.0D1}), (\ref{eq9.1B4}), and (\ref{eq9.1B5}), 
one has that
\begin{equation}
\label{eq9.1C4}
{\rm for\ every}\ j \in \N, \indent
\xi(t_j)\ =\ \exp\bigl(L^{(t(j))}(0)\bigl)\ =\  e^{j+2} \exp\bigl({\bf L}_j(0) \bigl)\ 
=\ e^{j+2} \zeps_j\ >\ \zeps_j. 
\end{equation}  
Hence [since $t_j \geq {\bf w}$ for each $j \in \N$ by the entire sentence of 
(\ref{eq9.1B3a})-(\ref{eq9.1B3b})], 
one has by (\ref{eq9.1C4}), (\ref{eq9.1B5}), and (\ref{eq9.0D4}) that
\begin{equation}  
\label{eq9.1C5}
{\rm for\ every}\ j \in \N,\ \ \
\zeps_j / \theta_j\ <\ \frac {\xi(t_j)} {-\phi'(t_j)}\ \leq\ \Delta. 
\end{equation}    
Now for every $j \in \N$, by (\ref{eq9.1C5}), 
(\ref{eq9.1C3}), and (\ref{eq9.1B1}),
one has that 
$h_j^2 \Delta \geq h_j^2 \zeps_j / \theta_j
= B_j \geq 9h_{j-1}^2 \Delta$, and hence
$h_j^2 \geq 9h_{j-1}^2$. 
Hence [recall (\ref{eq9.1A1})]
\begin{equation}
\label{eq9.1C6}
h_0 = 1\ \ \ 
{\rm and\ for\ every}\ j \in \N,\ \ \
h_j\ \geq\ 3h_{j-1}\ .
\end{equation}

     Finally, from (\ref{eq9.1B3a}) and Lemma \ref{lem6.11}(II), for each $j \in \N$, 
$-2^{-j} B_{j-1} \theta_{j-1} / B_j \leq \phi'(t_j) < 0$, and hence by (\ref{eq9.1B5}),
$0 < \theta_j \leq 2^{-j} B_{j-1} \theta_{j-1} / B_j$.
Hence 
\begin{equation}
\label{eq9.1C7}
{\rm for\ each}\ j \in \N,\ \ \ B_j \theta_j \leq 2^{-j} B_{j-1} \theta_{j-1}.
\end{equation}
In particular, the sequence of numbers $(B_0\theta_0,\ B_1 \theta_1,\ B_2 \theta_2,\ \dots)$
is strictly decreasing.  
Thus by (\ref{eq9.1C3}) (its equality) and (\ref{eq9.1C7}) and then (\ref{eq9.1A1}), one has that
\begin{equation}
\label{eq9.1C8}
\sum_{j=1}^\infty h_j^2 \zeps_j\ =\ \sum_{J=1}^\infty B_j \theta_j\ 
\leq\ \sum_{j=1}^\infty 2^{-j} B_{j-1} \theta_{j-1}\ 
\leq\ \sum_{j=1}^\infty 2^{-j} B_0 \theta_0\ =\ B_0 \theta_0\ =\ 1/9\ <\ \infty.        
\end{equation}  

    {\bf Sub-step 1D.\ \ A summary of some key properties 
of the parameters:}  \zhfb
(a) For each $j \in \N$, $\zeps_j \in (0, 1/9]$, $\theta_j \in (0, 1/9]$,
$\theta^*_j \in (0, 1/8]$, $I_j \geq 72$, and $h_j > 0$.  \zhfb
(b) As $j \to \infty$, $\zeps_j \to 0$ monotonically
and $\theta_j \to 0$ monotonically. \zhfb
(c)  As $j \to \infty$, $\theta^*_j \to 0$ monotonically
and $I_j \to \infty$ monotonically. \zhfb     
(d) For each $j \in \N$, one has that $\theta_j/\theta^*_j < 1$ and
$\theta_j^*\zeps_j I_j \leq 1$.   
\zhfb
(e) As $j \to \infty$, one has that $\theta_j/\theta^*_j \to 1$ and 
$\theta_j^*\zeps_j I_j \to 1$ and $\theta_j/h_j \to 0$. \zhfb
(f) One has that $\sum_{j=1}^\infty \zeps_j \leq \sum_{j=1}^\infty h_j^2 \zeps_j \leq 1/9$. 
\zhfb
(g) For each $j \geq 2$, one has that $\zeps_j \leq 9^{-j}/I_{j-1}$ and that
$h_j^2 \zeps_j/\theta_j = B_j 
\geq j \cdot [1 + \sum_{u=1}^{j-1} (h_u^2 \zeps_u/\theta_u)] > j$.  
\medskip

   Property (a) holds by the entire paragraph of (\ref{eq9.1B5})--(\ref{eq9.1B6}),
along with (\ref{eq9.1C6}).
Property (b) holds by (\ref{eq9.1C1}) and (\ref{eq9.1C2}). 
Property (c) holds by property (a) and (\ref{eq9.1B5})--(\ref{eq9.1B6})
and simple arithmetic.
Property (d) holds by property (a) and 
eqs.\ ({\ref{eq9.1B5})-(\ref{eq9.1B6}).
In property (e), the first two limits hold by properties (b) and (c) and
eqs.\ (\ref{eq9.1B5})--(\ref{eq9.1B6}); and the third limit holds by
property (a) ($0 < \theta \leq 1/9$ for all $j \in \N$) and 
eq.\ (\ref{eq9.1C6}) (which gives $h_j \to \infty$ as $j \to \infty$).
Property (f) holds by (a) and (\ref{eq9.1C6}) and (\ref{eq9.1C8}).
Finally, property (g) holds by (\ref{eq9.1C1}) and then
(\ref{eq9.1C3}) and (\ref{eq9.1A1}) (which gives $h_0^2 \zeps_0/ \theta_0 = 1$).
\medskip

   {\bf Step 2.}\ \  Here we repeat essentially verbatim the 
entire ``Step 2'' in the proof given in Section \ref{sc7} for 
Theorem \ref{thm3.4}.
Where ``Sub-step 1D(f)'' is cited in that context there, 
we simply cite the corresponding fact in Sub-step 1D(f) here. 
\medskip

   {\bf Step 3.}\ \
Here we shall essentially repeat Step 4 in the proof (in Section \ref{sc7})
of Theorem \ref{thm3.4}.
Because of one key difference (the random variables $X_k$ there were
bounded; here in general they are unbounded), we shall spell out the main details
for our context here.
\smallskip 
   
   In accordance with Step 2 above, refer to
Sub-steps 2A, 2B, and 2C in Section \ref{sc7}, in the proof of
Theorem \ref{thm3.4}.
By Sub-step 2A there, one has the following:
For each $j \in \N$, $E(X^{(j)}_0) = 0$ and
$\var(X^{(j)}_0) = \zeps_j$, and hence
$E(h_j X^{(j)}_0) = 0$ and 
$E[(h_j X^{(j)}_0)^2]  = \var(h_j X^{(j)}_0) = h_j^2 \zeps_j$.
Hence by Sub-step 2A there (the independence condition there), 
(\ref{eq7.2C1}), Lemma \ref{lem6.8}, and (\ref{eq9.1C8}) above, 
the random variable $X_0$ is square-integrable and satisfies 
$EX_0 = 0$ and 
$E(X_0^2) = \var(X_0) = \sum_{j=1}^\infty h_j^2\zeps_j \leq 1/9 < \infty$.
\smallskip

    By a trivial application of strict stationarity and Minkowski's inequality,
for each $j \in \N$ and each $n \in \N$,
$E[(\sum_{k=1}^n h_j X^{(j)}_k)^2] \leq n^2 E[(h_j X^{(j)}_0)^2] 
= n^2h_j^2 \zeps_j$.
Of course for any $n \in \N$, 
$\sum_{j=1}^\infty n^2h_j^2 \zeps_j < \infty$ 
by (\ref{eq9.1C8}) above (again).
Hence by Lemma \ref{lem6.8} (again),
for any positive integer $n$ and any nonempty set $S \subset \N$
(including $\N$ itself),
the random variable $\sum_{j \in S} \sum_{k=1}^n h_j X^{(j)}_k$ 
is square-integrable and has mean 0 and second moment 
(variance) $\sum_{j \in S} \var (\sum_{k=1}^n h_j X^{(j)}_k)$. 
\medskip 
 
{\bf Step 4.\ \ Proof of eq.\ (\ref{eq5.51}) in Theorem \ref{thm5.5}.}\ \
This step will be divided into six ``sub-steps'', including a technical lemma that will
be referred to as ``Sub-claim 4B''.
\medskip
    
     {\bf Sub-step 4A.}\ \ For each positive integer $j$, 
define the function $q_j: (0,1) \to [0, \infty)$ as follows:
\begin{equation}
\label{eq9.4A1}
q_j(u)\ :=\ 
\begin{cases}
 h_j\ \  {\rm if}\ \ 0 < u < \zeps_j \\
 0\ \  {\rm if}\ \ \zeps_j \leq u < 1.    
\end{cases}
\end{equation}

   {\bf Sub-claim 4B.}\ \ {\it For every\/} $j \in \N$ {\it and every\/} $x \in (1, \infty)$,  
\begin{equation}
\label{eq9.4B1}
\int_0^{f(x)} q_j^2(u) du\ \leq\ 2^{-(j+6)} g(x) \cdot (-\phi'(x)).
\end{equation}

     {\bf Proof of Sub-claim 4B.}\ \ Suppose $j$ is a positive integer.
Trivially by (\ref{eq9.4A1}), one has that 
\begin{equation}
\label{eq9.4B2}
{\rm for\ every}\ y \in (0,1), \quad
\int_0^y q_j^2(u)du\ =\ h_j^2 \cdot \min\{y, \zeps_j\}.
\end{equation}
Using just crude inequalities based on (\ref{eq9.4B2}), 
we shall divide the argument for (\ref{eq9.4B1}) into two cases
according to whether $1 < x \leq t_j$ or $t_j < x < \infty$.
\smallskip

   {\bf Case 1.\ \ $1 < x \leq t_j$ .}\ \ By (\ref{eq9.4B2}), (\ref{eq9.1C3}),
(\ref{eq9.1B5}), (\ref{eq9.1B3b}) [giving $g(t_j) \geq 2^{j+6} e^{j+2} B_j$], 
and (\ref{eq9.0A5}),
\begin{equation}
\nonumber
\int_0^{f(x)} q_j^2(u)du\ \leq\ h_j^2 \zeps_j\
=\ B_j \theta_j\ =\ B_j \cdot (-\phi'(t_j))\ 
<\ 2^{-(j+6)} g(t_j) \cdot (-\phi'(t_j))\ \leq\ 2^{-(j+6)} g(x) \cdot (-\phi'(x)).  
\end{equation}
Thus (\ref{eq9.4B1}) holds for Case 1.
\smallskip 

   {\bf Case 2.\ \ $t_j < x< \infty$.}\ \ 
Then $f(x) < f(t_j)$.
Recall from the entire sentence of (\ref{eq9.1B3a})-(\ref{eq9.1B3b})
that $t_j \geq {\bf w}$. 
By (\ref{eq9.4B2}), (\ref{eq9.0A1}), (\ref{eq9.0D1}) and (\ref{eq9.0A4a}),
then (\ref{eq9.0D3}), (\ref{eq9.1C4}) and (\ref{eq9.1B5}), then (\ref{eq9.1C3}), 
and finally (\ref{eq9.1B3b}),   
\begin{align}
\nonumber
\int_0^{f(x)} &q_j^2(u)du\ \leq\ h_j^2\, f(x)\ 
=\ h^2 \exp(\phi(x))\ 
\leq\ h_j^2\, \xi(x)\ 
=\ h_j^2 \cdot \frac {\xi(x)} {-\phi'(x)}  \cdot (-\phi'(x)) \\
\nonumber
&\leq\ h_j^2 \cdot \frac {\xi(t_j)} {-\phi'(t_j)}  \cdot (-\phi'(x))\
=\ h_j^2 \cdot \frac {e^{j+2}\, \zeps_j} {\theta_j}  \cdot (-\phi'(x))\ 
=\ B_j \cdot e^{j+2} \cdot  (-\phi'(x)) \\
& \leq\ 2^{-(j+6)} g(x) \cdot (-\phi'(x)). \nonumber   
\end{align}
Thus (\ref{eq9.4B1}) holds for Case 2. 
That completes the proof of  Sub-claim 4B.
\bigskip

     {\bf Sub-step 4C.}\ \ We shall return to Sub-step 4A and Sub-claim 4B later on. 
For now, in accordance with Step 2 above, refer in Section \ref{sc7}
to Sub-steps 2A, 2B, and 2C, including eq.\ (\ref{eq7.2C1}) 
and the sentence right after it. 
Also refer to (\ref{eq5.11}) and also to (\ref{eq9.0A1}).
To complete the proof of (\ref{eq5.51}), it obviously suffices to prove that
\begin{equation}
\label{eq9.4C1} 
{\rm for\ every}\ x \in (1, \infty),\ \ \ 
\int_0^{f(x)} Q_{|X(0)|}^2(u) du\ \leq\ g(x) \cdot (-\phi'(x)).
\end{equation}
\smallskip

     Define the nonnegative random variable $Y$ by
\begin{equation}
\label{eq9.4C2}
Y\ :=\ \sum_{j=1}^\infty h_j |X^{(j)}_0|\, .
\end{equation}
Then by (\ref{eq7.2C1}), 
$|X_0(\omega)| \leq Y(\omega)$ for every $\omega \in \Omega$. 
Hence trivially for any $t \geq 0$, $P(|X_0| > t) \leq P(Y > t)$;
and hence for any given $u \in (0,1)$,
$\{t \geq 0: P(|X_0| > t) \leq u \}  \supset \{t \geq 0: P(Y > t) \leq u \}$.
Hence by (\ref{eq5.11}), for any $u \in (0,1)$,
$Q_{|X(0)|}(u) \leq Q_Y(u)$.
Hence to prove (\ref{eq9.4C1}) and thereby complete the proof of (\ref{eq5.51}), 
it suffices to prove that
\begin{equation}
\label{eq9.4C3}
{\rm for\ every}\ x \in (1,\infty), \indent
\int_0^{f(x)} Q_Y^2(u)du\ \leq\ g(x) \cdot (-\phi'(x)). 
\end{equation}  

   {\bf Sub-step 4D.}\ \
For each $j \in \N$, define the positive number $a_j$ and the 
event $A_j$ as follows:
\begin{equation}
\label{eq9.4D1}
a_j\ :=\ 2\zeps_j \indent {\rm and} \indent A_j\ :=\ \{|X_0^{(j)}|=1\}. 
\end{equation}
From (\ref{eq9.4D1}) and Sub-step 2B (repeated directly from 
Sub-step 2B in Section 7), one has that
for each $j \in \N$, the random variable $X_0^{(j)}$ takes
only the values $-1$, 0, and 1, and that $P(A_j) = P(|X^{(j)}_0| = 1) = \zeps_j$.
(That equality is tacitly in play more than once below.)\ \   
By (\ref{eq9.4D1}), for each $j \in \N$, 
one now has that $a_j = 2P(A_j)$,
and also one has the equality of random
variables $I(A_j) =   |X_0^{(j)}|$, where $I(\dots)$ denotes the 
indicator function.
By (\ref{eq9.4C2}), one now has the equality of random variables
\begin{equation}
\label{eq9.4D2}
Y\ :=\ \sum_{j=1}^\infty \bigl(h_j \cdot I(A_j) \bigl).
\end{equation}
Further, by Sub-step 1D(a) and eq.\ (\ref{eq9.1C1}), one now has that
$P(A_1) < 1/2$ and $P(A_j) \leq (1/2)P(A_{j-1})$ for every $j \geq 2$.    
And now by (\ref{eq9.1C6}), (\ref{eq9.4D2}), (\ref{eq9.4D1}), 
and Lemma \ref{lem6.12},
one has that $0 < \dots < a_3 < a_2 < a_1 < 1$ 
and also $\lim_{j \to \infty} a_j = 0$, and that
\begin{align}
\label{eq9.4D3a}
&Q_Y(u) = 0\ \ {\rm for\ all}\ u \in [a_1,1); \quad {\rm and} \\
\label{eq9.4D3b}
&{\rm for\ each}\ j \in \N,\ {\rm one\ has\ that}\ 
Q_Y(u) \leq 2h_j\ {\rm for\ all}\ u \in [a_{j+1}, a_j).
\end{align}

   {\bf Sub-step 4E.}\ \ 
Now for a given $j \in \N$, we shall use here the function 
$q_j: (0,1) \to [0, \infty)$ defined in eq.\ (\ref{eq9.4A1}). 
First let us define, for each $j \in \N$, two other closely related functions 
$q_j^*: (0,1) \to [0, \infty)$ and $q_j^{**}: (0,1) \to [0, \infty)$, as follows: 
\begin{equation}
\label{eq9.4E1}
q_j^*(u)\ :=\ h_j\, I_{(0,a(j))}(u)  \indent {\rm and} \indent
q_j^{**}(u)\ :=\ h_j\, I_{[a(j+1),a(j))}(u),
\end{equation}
where $I_{…}$ denotes the indicator function, and for each $i \in \N$,
the number $a_i$ is written as $a(i)$ for typographical convenience.
For a given $u \in [a_1, 1)$, one has that
$\sum_{i=1}^\infty [q_i^{**}(u)]^2 = \sum_{i=1}^\infty 0 = 0$; and for a given
$j \in \N$ and a given $u \in [a_{j+1}, a_j)$, one has that
$\sum_{i=1}^\infty [q_i^{**}(u)]^2 = h_j^2 + \sum_{i \in \N-\{j\}} 0 = h_j^2$.
Hence by (\ref{eq9.4D3a})-({\ref{eq9.4D3b}), one has that
\begin{equation}
\label{eq9.4E2}
{\rm for\ each}\ u \in (0,1), \indent 
Q_Y^2(u)\ \leq\ 4 \cdot \sum_{i=1}^\infty [q_i^{**}(u)]^2. 
\end{equation}

   {\bf Sub-step 4F.}\ \ 
For any given $j \in N$, one has by (\ref{eq9.4E1})} that
$0 \leq q^{**}(u) \leq q^*(u)$ for all $u \in (0,1)$.
Hence for any given $c \in (0,1)$, by (\ref{eq9.4E2}),
(\ref{eq9.4D1}), and (\ref{eq9.4A1}) [again (\ref{eq9.4B2})],
\begin{align}
\nonumber
\int_0^c& Q_Y^2(u)du\ \leq\ \int_0^c 4 \cdot \sum_{j=1}^\infty [q_j^{**}(u)]^2 du\
=\ 4 \cdot \sum_{j=1}^\infty \int_0^c [q_j^{**}(u)]^2 du\
\leq\ 4 \cdot \sum_{j=1}^\infty \int_0^c [q_j^*(u)]^2 du\\
\nonumber
&=\ 4 \cdot \sum_{j=1}^\infty (h_j^2 \cdot \min\{c, a_j\})\
=\ 8 \cdot \sum_{j=1}^\infty (h_j^2 \cdot \min\{c/2, \zeps_j\})\
=\ 8 \cdot \sum_{j=1}^\infty \int_0^{c/2} q_j^2(u)du\  
\leq\ 8 \cdot \sum_{j=1}^\infty \int_0^c q_j^2(u)du.    
\end{align} 
Applying that and Sub-claim 4B, one has that for any $x \in (1, \infty)$, 
\begin{align}
\nonumber
\int_0^{f(x)} Q_Y^2(u)du\ 
&\leq\ 8 \cdot \sum_{j=1}^\infty \int_0^{f(x)} q_j^2(u)du\ 
\leq\ 8 \cdot \sum_{j=1}^\infty \Bigl[ 2^{-(j+6)} \cdot g(x) \cdot \bigl(-\phi'(x)\bigl) \Bigl]\\
\nonumber
&<\  \sum_{j=1}^\infty \Bigl[ 2^{-j} \cdot g(x) \cdot \bigl(-\phi'(x)\bigl) \Bigl]\
=\ g(x) \cdot \bigl(-\phi'(x)\bigl).  
\end{align}
Thus (\ref{eq9.4C3}) holds.
That completes the proof of eq.\ (\ref{eq5.51}).
\medskip

    {\bf Step 5.  Proof of eq.\ (\ref{eq5.52}) in Theorem 5.5.}\ \
For each $j  \in \N$,
one has the following:  
First, by (\ref{eq9.1B4}), (\ref{eq9.0A6b}), and (\ref{eq9.1B5}), 
$[{\rm slope\ of}\ {\bf L}_j] = [{\rm slope\ of}\ L^{(t(j))}] = \phi'(t_j) = -\theta_j$,  
and hence for every $x \in \R$, 
${\bf L}_j(x) = {\bf L}_j(0) - \theta_jx$.
Hence for every $x \in \R$, by (\ref{eq9.1B5}), 
\begin{equation}
\label{eq9.51}
\exp\bigl({\bf L}_j(x)\bigl)\ =\ 
\exp\bigl({\bf L}_j(0)\bigl)\, \cdot\, \exp(-\theta_j x)\
=\ \zeps_j \cdot \exp(-\theta_j x).  
\end{equation}
Second, by (\ref{eq9.1B4}), for every $x \in \R$,
$\exp({\bf L}_j(x)) 
= e^{-(j+2)} \cdot \exp \bigl(L^{(t(j))}(x) \bigl)$;
and hence by (\ref{eq9.0A6a}) and (\ref{eq9.0A1}), 
for every $x \in [1, \infty)$,
\begin{equation}
\label{eq9.52}
\exp\bigl({\bf L}_j(x)\bigl)\ 
\leq\ e^{-(j+2)} \cdot \exp \bigl(\phi(x) \bigl)\
=\ e^{-j} \cdot e^{-2} \cdot f(x) .
\end{equation}
By (\ref{eq9.51}) and (\ref{eq9.52}), one has that for every $j \in \N$
and every $x \in [1, \infty)$,
\begin{equation}
\label{eq9.53}
\zeps_j \cdot \exp(-\theta_j x)\ \leq\ e^{-j} \cdot e^{-2} \cdot f(x).
\end{equation}

   Now by Lemma \ref{lem6.5}(6) and Sub-step 2A in Step 2 (in Section \ref{sc7})
of the proof of Theorem {\ref{thm3.4},  and then (\ref{eq9.53}) 
(and the inequality $6 < e^2$),
for every $j \in \N$ and every $n \in \N$,
\begin{equation}
\nonumber
\beta_{X(j)}(n)\ \leq\ 6 \zeps_j (1 - \theta_j)^n\
\leq\ 6\zeps_j \cdot \exp(-\theta_j n)\
<\ e^{-j} \cdot f(n).
\end{equation}
[Here the notation $X(j)$ means the Markov chain $X^{(j)}$, 
and is employed in subscripts for typographical convenience.]\ \ 
Hence by Lemma \ref{lem6.7}, for every $n \in \N$,
$\beta_X(n) \leq \sum_{j=1}^\infty \beta_{X(j)}(n) 
\leq \sum_{j=1}^\infty e^{-j} f(n) < f(n)$. 
Referring to (\ref{eq2.33}), one now has that eq.\ (\ref{eq5.52}) holds.
\medskip

   {\bf Step 6/7.}\ \ Our next task is to prove that under the assumptions in 
Theorem \ref{thm5.5}, conclusions (iii) and (iv) of Theorem \ref{thm3.4} hold.
To accomplish this, we first note that Sub-step 1D((a)-(g)) here in Section \ref{sc9}
gives all essential information in the same format as in 
Sub-step 1D((a)-(g)) in either Section \ref{sc7} or Section \ref{sc8}.
\smallskip

   The proof (in our context here) of statement (iii) of Theorem \ref{thm3.4} mimics 
the argument in Step 6 of Section \ref{sc7} (in the proof of Theorem \ref{thm3.4} itself),
in exactly the same way as was done in Section \ref{sc8} (in the proof of
Theorem \ref{thm4.4}).
\smallskip

   The proof (in our context here) of statement (iv) of Theorem \ref{thm3.4}
is the same as the argument in Step 7 of Section \ref{sc8} (in the proof of
Theorem \ref{thm4.4}), three quarters of which was the same as in
Step 7 of Section \ref{sc7} (in the proof of Theorem \ref{thm3.4} itself).
\medskip    

        {\bf Step 8.}\ \ Our final task is to prove that in our context here,
conclusion (iv) [with $g_n := g(n)$ in (\ref{eq4.41})] 
in the statement of Theorem \ref{thm4.4} holds.  
From the formulation of that conclusion in that theorem, 
recall for each positive integer $n$ the definition 
${\bf h}(n) := n^{-1} \var(\sum_{k=1}^n  X_k)$.
In Step 8 here the proof of eq.\ (\ref{eq4.41}) will be given.
In Step 9 below, the proof of eq.\ (\ref{eq4.42}) will be given, and the proof of
Theorem \ref{thm5.5} will then be complete.   
\medskip     
     
    {\bf Proof of eq.\ (\ref{eq4.41}) with $g_n = g(n)$.}\ \ In its overall structure, 
this proof will be like that given in Step 8 in Section \ref{sc8} for the same equation.
However, because of several differences in technical details, we shall spell it out in
detail here.
\smallskip

This proof will be divided into two ``sub-steps'', of which the first is labeled ``Sub-claim 8A''.
First refer again to the positive function $g(x),\ x \in [1, \infty)$ in the statement of 
Theorem \ref{thm5.5}.
\medskip

     {\bf Sub-claim 8A.}\ \ {\it Suppose $j$ and $n$ are each a positive integer.
Then} 
\begin{equation}
\label{eq9.8A1}
n^{-1}\, \var\biggl(\, \sum_{k=1}^n h_j X_k^{(j)} \biggl)\ \leq\ 2^{-j} g(n).
\end{equation}

    {\bf Proof of Sub-claim 8A.}\ \ 
The argument for Sub-claim 8A will be divided into two cases
according to whether $g(n) < 2^{j+1} B_j$ or $g(n) \geq 2^{j+1} B_j$. 
\medskip 

   {\bf Case 1.}\ \ $g(n) < 2^{j+1} B_j$.\ \ 
By (\ref{eq9.1C2}) and (\ref{eq9.1B1a}), 
$\theta_j \leq {\bf c}_j \leq 2^{-j} \cdot n^{-1} g(n) \cdot  B_j^{-1}$.
Hence by Lemma \ref{lem6.5}(4) and eq.\ (\ref{eq9.1C3}),
\begin{align}
\nonumber
(1/n) \cdot  \var\biggl(\, \sum_{k=1}^n h_j X_k^{(j)} \biggl)\ 
&\leq\ (1/n) \cdot n^2\, \var(h_j X_0^{(j)})\
=\ n h_j^2 \zeps_j\  =\ n B_j \theta_j\\
\nonumber
&\leq\ n B_j \cdot \bigl(2^{-j} \cdot n^{-1} g(n) \cdot B_j^{-1} \bigl)\
=\ 2^{-j} g(n).  
\end{align} 
Thus (\ref{eq9.8A1}) holds for Case 1.
\medskip

     {\bf Case 2.}\ \ $g(n) \geq 2^{j+1} B_j$.\ \ 
By Lemma \ref{lem6.5}(9) and eq.\ (\ref{eq9.1C3}), 
\begin{equation}
\nonumber 
n^{-1} \var \biggl(\, \sum_{k=1}^n h_j X_k^{(j)} \biggl)\ \leq\ h_j^2 \cdot 2 \zeps_j / \theta_j\
=\ 2 B_j\  \leq\ 2^{-j} g(n). 
\end{equation} 
Thus (\ref{eq9.8A1}) holds for Case 2.
That completes the proof of Sub-claim 8A.
\medskip

    {\bf Sub-step 8B.}\ \
Recall again the definition of the numbers ${\bf h}(n)$, $n \in \N$ from the statement of 
Theorem \ref{thm4.4}.
By Lemma \ref{lem6.8} and Sub-claim 8A, for each positive integer $n$,
\begin{equation}
\nonumber
{\bf h}(n)\ =\ n^{-1} \var \biggl(\, \sum_{k=1}^n X_k \biggl)\ 
=\ n^{-1} \sum_{j=1}^\infty \var \biggl(\, \sum_{k=1}^n h_j X_k^{(j)} \biggl)\
\leq\ \sum_{j=1}^\infty 2^{-j} g(n)\ =\ g(n).    
\end{equation}   
Thus eq.\ (\ref{eq4.41}) [with $g_n = g(n)$] holds.  Step 8 is complete.
\zhfb

     {\bf Step 9.  Proof (in our context here) of eq.\ (\ref{eq4.42}).}\ \ 
The argument is essentially the same as in Step 9 in Section {\sc8}.
The only (minor) differences are as follows:
\medskip

   In the sentence after that of (\ref{eq8.9A1}), the reference to ``the second sentence after
(\ref{eq8.1B5})'' needs to be replaced by a reference to (\ref{eq9.1C2}). 
\medskip

   In the proof of (\ref{eq8.9B1b}), the reference to Sub-step 1D(g) (there in Section \ref{sc8})
needs to be replaced by a reference to (the corresponding identical fact in) 
Sub-step 1D(g) here in Section \ref{sc9}.
\medskip

   In the proof of (\ref{eq8.9B1c}), the material right after (\ref{eq8.9B4}) 
[including eq.\ (\ref{eq8.9B5})] needs to be replaced by the following:  
For each $j \in \N$ and each $u \geq j+1$, one has by (\ref{eq9.1C7}) (and the sentence 
right after it) that $B_u \theta_u \leq 2^{-u} B_{u-1} \theta_{u-1} \leq 2^{-u} B_j \theta_j$. 
Hence considering $j \in \N$, one has by (\ref{eq9.1C3}) (its equality, twice) that
\begin{equation}
\nonumber
\sum_{u=j+1}^\infty h_u^2 \zeps_u\ =\ \sum_{u=j+1}^\infty B_u \theta_u\ 
\leq\ B_j \theta_j \sum_{u=j+1}^\infty 2^{-u}\ =\ 2^{-j} h_j^2 \zeps_j\ 
=\ o(h_j^2 \zeps_j)\ \ {\rm as}\ j \to \infty.
\end{equation}
Hence considering $j \in \N$, one has by Lemma \ref{lem6.5}(4) that
\begin{equation}
\nonumber
\sum_{u = j+1}^\infty \var \biggl(\, \sum_{k=1}^{K M(j)} h_u X^{(u)}_k \biggl)\
\leq\ \sum_{u=j+1}^\infty K^2 M_j^2\, \var (h_u X^{(u)}_0)\ 
=\ K^2 M_j^2 \sum_{u=j+1}^\infty h_u^2 \zeps_u\
=\ o(M_j^2 h_j^2 \zeps_j)\ \ {\rm as}\ j \to \infty. 
\end{equation}
Thus (\ref{eq8.9B1c}) holds in our context here.
\medskip

Thus (\ref{eq4.42}) holds in our context here (with no other changes needed in the proof).
That completes the proof of Theorem \ref{thm5.5}.

\section{Appendix}
\label{sc10}

     Recall the discussion in Remark \ref{rem4.6} regarding the refined assumption 
[\cite{ref-journal-MP2}, eq.\ (3.2)] in the CLT in [\cite{ref-journal-MP2}, Theorem 4],  
which sharpened Theorem \ref{thm5.3} and (as a special case) Theorem \ref{thm4.2}.
As noted there, for strictly stationary sequences with finite second moments, 
it was shown in \cite{ref-journal-MP2} that that assumption implies eq.\ (\ref{eq4.61}) 
(in Section 4 here).
Hence by eq.\ (\ref{eq4.42}) in Theorem \ref{thm4.4}, the construction in Section \ref{sc8}
for that theorem fails to satisfy [\cite{ref-journal-MP2}, eq.\ (3.2)].
\medskip
 
That argument is somewhat ``indirect''.
Here in this section [as indicated in Remark \ref{rem4.6}(D)], 
to give some more insight, we shall provide explicitly the 
formulation of [\cite{ref-journal-MP2}, eq.\ (3.2)], and then provide in one
nice context  --- the construction in Section \ref{sc8} for Theorem \ref{thm4.4}, 
involving bounded random variables --- 
a ``direct'' illustration of how it is that that construction fails to satisfy that 
equation [\cite{ref-journal-MP2}, eq.\ (3.2)].
A common link between this ``direct'' calculation below and the proof in Section \ref{sc8}
of eq.\ (\ref{eq4.42}) is eq.\ (\ref{eq8.9B1d}) in Sub-claim 9B in Section \ref{sc8}
[see eq.\ (\ref{eq8.9A1})], involving some sparse (in some sense) but 
ever-recurring manifestations of ``approximate quadratic growth'' of variances of the partial sums.        
\medskip

   Suppose $(X_k, k \in \Z)$ is a strictly stationary (not necessarily Markovian) 
sequence of random variables
such that $E[X_0^2] < \infty$ and $E[X_0] = 0$.
For each positive integer $n$, define the number
$\sigma_n^2 := \var (\sum_{k=1}^n X_k)$.   
Essentially, that assumption [\cite{ref-journal-MP2}, eq.\ (3.2)] is as follows:
\begin{equation}
\label{eq10.01}
\lim_{n \to \infty} \frac {1} {\sigma_n^2} \sum_{i=1}^n
 i \int_0^{\alpha(i)} Q_{|X(0)|}^2(u)\, du\ =\ 0. 
\end{equation}
In the formulation in [\cite{ref-journal-MP2}, eq.\ (3.2)] itself, in the upper boundary
of the integrand, the $\alpha(i)$ was replaced by the corresponding dependence 
coefficient for a mixing condition that is similar to but weaker (hence more general,
more versatile) than $\alpha$-mixing.
However, for strictly stationary Markov chains, our context here, 
the mixing condition that was used there coincides with $\alpha$-mixing 
(with the corresponding dependence coefficients being equal), 
as an elementary consequence of eq.\ (\ref{eq2.41}). 
Accordingly, in the integrand in (\ref{eq10.01}), we shall for simplicity leave the 
$\alpha(i)$ intact.    
\zhfb

     {\bf A look at the left side of (\ref{eq10.01}) for the construction in Section 8.}\ \ 
Our task here is to show with a relatively simple calculation how the (strictly stationary,
countable-state, reversible) Markov chain constructed in Section \ref{sc8} 
for Theorem \ref{thm4.4} fails to satisfy (\ref{eq10.01}).
All notations and equations in Section \ref{sc8} 
(as well as any material in Section \ref{sc7} that is carried over to Section \ref{sc8})
will be taken for granted.
The notation $\alpha(n)$ in (\ref{eq10.01}) of course refers to $\alpha_X(n)$,
for the (strictly stationary, countable-sate, reversible) Markov chain $X$ constructed in 
Section \ref{sc8} for Theorem \ref{thm4.4}.
The analysis below will be divided into several ``steps''
(including two ``claims'').
\medskip

   {\bf Step 1.}\ \ Refer to the first sentence of Background \ref{bckg3.1} (or to 
property (ii) in Theorem \ref{thm4.4} itself, say with the natural restriction  
$g_n= o(n)$ as $n \to \infty$).
One automatically has that the Markov chain $X$ is $\alpha$-mixing;
that is, $\alpha(n) := \alpha_X(n) \to 0$ as $n \to \infty$. 
\medskip

  {\bf Step 2.\ \ More on the dependence coefficients $\alpha(n)$.}\ \ 
Recall again that for any $j \in \N$ and any $k \in \Z$, 
the random variable $X^{(j)}_k$ takes its values in the set $\{-1, 0, 1\}$.
By an elementary direct calculation, or trivially by an old classic inequality
(see e.g.\ [\cite{ref-journal-IbrLin}, Theorem 17.2.1] or
[\cite{ref-journal-Bradley2007}, Vol.\ 1, Theorem 1.11]), one has that
for any $j \in \N$ and any $k \in \Z$, 
\begin{equation}
\label{eq10.21}
\bigl| \cov(X^{(j)}_0, X^{(j)}_k) \bigl|\ 
\leq\ 4 \alpha \bigl(\sigma(X^{(j)}_0), \sigma(X^{(j)}_k) \bigl).
\end{equation}
For any positive integers $i$ and $j$, one has by (\ref{eq2.41}) and 
the paragraph right before that of (\ref{eq7.2D1}), followed by 
(\ref{eq10.21}) and Lemma \ref{lem6.5}(8a), that
\begin{align}
\nonumber
\alpha(i)\ &:=\ \alpha_X(i)\ =\ \alpha \bigl(\sigma(X_0), \sigma(X_i) \bigl)\ 
\geq\ \alpha \bigl(\sigma(X^{(j)}_0), \sigma(X^{(j)}_i) \bigl)\
\geq\ (1/4)\, |\cov(X^{(j)}_0, X^{(j)}_i)|\\
\label{eq10.22}
&=\ (1/4) \zeps_j(1 - \theta_j)^i\ \geq\ (1/4) \zeps_j(1 - i \theta_j).
\end{align}

   Now for each $j \in \N$, refer to the positive integer $M_j$ in (\ref{eq8.9A1}).
By that equation, $M_j \theta_j \leq 1/2$ for every $j \geq 4$.
Recall from (\ref{eq8.9A1a}) that $M_j \to \infty$ as $j \to \infty$.
By (\ref{eq10.22}), 
\begin{align}
\label{eq10.23}
\forall j \geq 4,\ \forall i \in \{1,2,\dots, M_j\}, \ \ \ 
\alpha(i)\ \geq\ (1/4) \zeps_j(1 - M_j \theta_j)\ \geq\ (1/4) \zeps_j \cdot (1 - 1/2)\ =\ (1/8) \zeps_j.
\end{align}
As a triviality, (\ref{eq10.23}) conveys that $\alpha(i) := \alpha_X(i) > 0$ for every $i \in \N$.
\medskip

   {\bf Step 3.\ \ The quantile function $Q_{|X(0)|}$.}\ \
Here for review, we shall include some concrete details that could be omitted with the use of
``generic'' properties of the quantile functions in Section \ref{sc5}. 
Recall that for any $j \in \N$ and any $k \in \Z$, the random 
variable $X^{(j)}_k$ takes its values in the set $\{-1, 0, 1\}$.
Step 3 here will be divided into two `` sub-steps'', of which the first is a ``sub-claim'': 
\medskip   

   {\bf Sub-claim 3A.}\ \ {\em For any $r \in (0,1/2)$, there exists $u_0 = u_0(r)$ such that
for all $u \in (0, u_0]$, $Q_{|X(0)|}(u) \geq r$.}

\medskip

   {\bf Proof of Sub-claim 3A.}\ \ Suppose $r \in (0,1/2)$.
By eq.\ (\ref{eq8.101}) and a routine calculation,  
$\lim_{K \to \infty} \sum_{j=1}^K h_j \allowbreak 
= \sum_{j=1}^\infty h_j = 1/2$.
Let $J \in \N$ be fixed sufficiently large that $\sum_{j=1}^J h_j > r$.
By (\ref{eq7.2D1}) (which carries over to Section \ref{sc8}),
with $z_j =1$ for $1 \leq j \leq J$ and $z_j = 0$ for $j \geq J+1$,  one has that
$P(|X_0| > r) \geq P(X_0 = \sum_{j=1}^J h_j) > 0$.  
Let $u_0 \in (0,1)$ be such that $P(|X_0| > r) > u_0$.
\smallskip

    Now suppose $u \in (0, u_0]$.   Then one has that $P(|X_0| > r) > u$.
Hence for all $t \in [0, r]$, $P(|X_0| > t) > u$.
Hence $\inf \{ t \geq 0: P(|X_0| > t) \leq u\} \geq r$.
That is, by (\ref{eq5.11}), $Q_{|X(0)|}(u) \geq r$.   
Thus Sub-claim 3A holds.
\medskip

   {\bf Sub-step 3B.}\ \ By eq.\ (\ref{eq8.101}) and eq.\ (\ref{eq7.2C1}) 
(which, again, carries over to Section \ref{sc8}), followed by a routine calculation (again),
for each $\omega \in \Omega$,
$|X_0(\omega)| \leq \sum_{j=1}^\infty h_j = 1/2$. 
Thus $P(|X_0| > 1/2) = 0$.
 Hence for any $u \in (0,1)$, again referring to (\ref{eq5.11}), 
 $Q_{|X(0)|}(u) = \inf \{t \geq 0: P(|X_0| > t) \leq u\} \leq 1/2$
(since for any such $u$, the number $1/2$ itself is an element of that set).
Combining that with Sub-claim 1A, one now has that
$\lim_{u \to 0+} Q_{|X(0)|}(u) = 1/2$.
As an elementary consequence (see Notations \ref{nt2.1} and the last sentence of 
each of Steps 1 and 2 above), one has that
$\int_0^{\alpha(n)} Q_{|X(0)|}^2 (u)\, du\ \sim\ (1/4)\, \alpha(n)$ as $n \to \infty$.
Accordingly, let ${\bf N}$ be a positive integer such that (say)  
\begin{equation}
\label{eq10.3B1}
 {\rm for\ all}\ n \geq {\bf N},\ \ \ \int_0^{\alpha(n)} Q_{|X(0)|}^2 (u)\, du\ \geq\ (1/8)\, \alpha(n).
\end{equation}

   {\bf Step 4.}\ \ Recall again from (\ref{eq8.9A1a}) that $M_j \to \infty$ as $j \to \infty$.
\smallskip

   For each positive integer $j$, define  the integer $K_j := [M_j/2]$ (the greatest integer
that is $\leq M_j/2$).
Then $K_j \to \infty$ as $j \to \infty$.
Referring to the integer ${\bf N}$ in (\ref {eq10.3B1}), let ${\bf J}$ be a positive integer 
such that 
\begin{equation}
\label{eq10.41}
{\bf J} \geq 4 \indent {\rm and} \indent \forall j \geq {\bf J},\ \ K_j \geq {\bf N}.
\end{equation} 

   {\bf Claim 5.}\ \ {\it If $j \geq {\bf J}$ and $i \in \{K_j, K_j+1, K_j+2, \dots, M_j\}$, then
$\int_0^{\alpha(i)} Q_{|X(0)|}^2 (u)\, du \geq (1/64) \zeps_j$.}
\medskip

   {\bf Proof of Claim 5.}\ \  Since $i \geq K_j \geq {\bf N}$ by (\ref{eq10.41}), one has by
(\ref{eq10.3B1}) that $\int_0^{\alpha(i)} Q_{|X(0)|}^2 (u)\, du \geq (1/8) \alpha(i)$.
Since $j \geq {\bf J} \geq 4$ by (\ref{eq10.41}), and $1 \leq i \leq M_j$, one has by   
(\ref{eq10.23}) that $\alpha(i) \geq (1/8) \zeps_j$.
From those two observations, Claim 5 follows.
\medskip

   {\bf Step 6.}\ \ In what follows, for a given $j \in \N$, the numbers $M_j$ and $K_j$ will 
also be written respectively as $M(j)$ and $K(j)$, for typographical convenience.
\smallskip

For each $j \in \N$, one has by the definition of $K_j$ in Step 4 that 
$K_j+1 > M_j/2 \geq K_j$ and (hence) $M_j - K_j \geq M_j/2$.
And hence also by Claim 5, one has that for each $j \geq {\bf J}$, 
\begin{align}
\nonumber 
\sum_{i=1}^{M(j)} & i \int_0^{\alpha(i)} Q_{|X(0)|}^2 (u)\, du\
\geq \sum_{i=K(j)+1}^{M(j)} i \int_0^{\alpha(i)} Q_{|X(0)|}^2 (u)\, du\ 
\geq \sum_{i=K(j)+1}^{M(j)} i \cdot (1/64) \zeps_j\\
\nonumber
&\geq \sum_{i=K(j)+1}^{M(j)} (K_j+1) \cdot (1/64) \zeps_j\ 
\geq \sum_{i=K(j)+1}^{M(j)} (M_j/2) \cdot (1/64) \zeps_j\
=\ (M_j - K_j) \cdot (M_j/2) \cdot (1/64) \zeps_j\\
\label{eq10.61}
&\geq\ (M_j/2) \cdot (M_j/2) \cdot (1/64) \zeps_j\
=\ (1/256) M_j^2 \zeps_j.    
\end{align}     
 
    {\bf Step 7.}\ \ 
Recall for $n \in \N$ the notation $\sigma_n^2 := \var (\sum_{k=1}^n X_k)$, from
the sentence right before that of eq.\ (\ref{eq10.01}). 
By eq.\ (\ref{eq8.9B1d}) (with $K=1$ there) in Sub-claim 9B in Section 8, one has that
$\sigma_{M(j)}^2 \sim M_j^2 h_j^2 \zeps_j$ as $j \to \infty$.
Now with the calculations below restricted to $j \geq {\bf J}$, 
one has by (\ref{eq10.61}) [its entire sentence] and (\ref{eq8.101}) that as $j \to \infty$,
\begin{equation}
\nonumber
\frac {1} {\sigma_{M(j)}^2} \sum_{i=1}^{M(j)} i \int_0^{\alpha(i)} Q_{|X(0)|}^2 (u)\, du\
\geq\ \frac {1} {\sigma_{M(j)}^2} \cdot \frac {1} {256} M_j^2 \zeps_j \
\sim\  \frac {1} {M_j^2 h_j^2 \zeps_j} \cdot \frac {1} {256} M_j^2 \zeps_j\
=\ \frac {1} {256\, h_j^2}\ \to\ \infty.  
\end{equation}
Since $M_j \to \infty$ as $j \to \infty$ (again see Step 4), it follows that
for the construction in Section \ref{sc8} for Theorem \ref {thm4.4},
eq.\ (\ref{eq10.01}) fails to hold.
\zhfb

   {\bf Acknowledgement.}\ \ The author thanks the referee for correcting several 
errors/typos, for calling attention to the references \cite {ref-journal-MP2} and
\cite {ref-journal-ZWV}, for suggesting the question posed in Remark \ref{rem4.6}(A)(D),  
and overall for various comments that led to a considerable improvement of this paper.

\end{document}